# Proof of Riemann's zeta-hypothesis

**Arne Bergstrom**

# Proof of Riemann's zeta-hypothesis

By ARNE BERGSTROM

### Abstract

Make an exponential transformation in the integral formulation of Riemann's zeta-function $\zeta(s)$ for $Re(s) > 0$. Separately, in addition make the substitution $s \to 1 - s$ and then transform back to $s$ again using the functional equation. Using residue calculus, we can in this way get two alternative, equivalent series expansions for $\zeta(s)$ of order $N$, both valid inside the "critical strip", i e for $0 < Re(s) < 1$. Together, these two expansions embody important characteristics of the zeta-function in this range, and their detailed behavior as $N$ tends to infinity can be used to prove Riemann's zeta-hypothesis that the nontrivial zeros of the zeta-function must all have real part ½.

## 1. Introduction

Riemann's zeta-hypothesis from 1859 [11] is expressed as follows:

CONJECTURE 1.1. *The nontrivial zeros of the Riemann zeta-function $\zeta(s)$ all have real part $Re(s) = ½$.*

The Riemann zeta-hypothesis is the most famous of the few still unsolved problems on Hilbert's list of twenty-three mathematical challenges, which he presented in 1900 at the dawn of the new century [12, 18]. It is also one of the seven Millennium Problems [19] named in 2000 by the Clay Mathematics Institute.

It can be shown (cf [15]) that the nontrivial zeros of the zeta-function must lie inside the "critical strip", i e for $0 < Re(s) < 1$, which is the range studied in this paper.

The Riemann zeta-hypothesis has been computationally verified for $|Im(s)|$ at least up to 2.4 trillion [17].

The intriguing possibility has been suggested that the Riemann zeta-function could correspond to a quantum-physical problem with its zeros corresponding to energy eigenvalues. The underlying physical problem would then correspond to a chaotic quantum system without time-reversal symmetry [4, 5].

---

With ($\sigma$ and $t$ are real)

$$s = \sigma + i\,t$$

Riemann's zeta-function $\zeta(s)$ can be defined as the following series, convergent for $\sigma > 1$,

$$\zeta(s) = \sum_{n=1}^{\infty} \frac{1}{n^s}$$

This Dirichlet series can also be expressed as follows (for $\sigma > 1$),

$$(1 - 2^{(-s)})\,\zeta(s) = \sum_{n=1}^{\infty} \frac{1}{(2\,n - 1)^s}$$

In Sects 2 through 5 below a modification of this latter series will be derived, giving the equivalent pair (9) and (11), which are valid also inside the critical strip. Although it will be shown that (9) and/or (11) are somewhat similar to previous results found in the literature, the approach described in the following permits a more detailed analysis, leading to a proof of Conjecture 1.1.

The proof of Riemann's zeta-hypothesis given in this paper is based on the following two fundamental properties of the Riemann zeta-function:

*the integral representation (1)*, valid for $Re(s) > 0$ [10, 14],

$$(1 - 2^{(1-s)})\,\Gamma(s)\,\zeta(s) = \int_0^{\infty} \frac{w^{(s-1)}}{e^w + 1}\,dw \tag{1}$$

*the functional equation (2)*, valid for all $s$ [7, 13],

$$\zeta(s) = 2^s\,\pi^{(s-1)}\,\sin\left(\frac{1}{2}\,s\,\pi\right)\Gamma(1-s)\,\zeta(1-s) \tag{2}$$

## 2. Variable transformation

We start by transforming the variable $w$ in (1) as follows

$$w = e^u$$

i e

$$(1 - 2^{(1-s)})\,\Gamma(s)\,\zeta(s) = \int_{-\infty}^{\infty} \frac{(e^u)^s}{e^{(e^u)} + 1}\,du \tag{3}$$

The integration variable $w$ in (1) being real, we can also set $u$ real. Then

$$(1-2^{(1-s)})\,\Gamma(s)\,\zeta(s)=\int_{-\infty}^{\infty}\frac{e^{(s\,u)}}{e^{(e^{u})}+1}\,du \tag{4}$$

Consider the integrand

$$F(u)=\frac{e^{(s\,u)}}{e^{(e^{u})}+1} \tag{5}$$

and extend $u$ to the entire complex plane,

$$u=x+i\,y$$

Extended over the complex plane, F($u$) is an analytic (meromorphic) function.

## 3. Poles and residues

We next calculate the poles of F($u$) above, i e we want to find all $u$ that satisfy the equation

$$e^{(e^{u})}+1=0$$

which can be verified to have the following solutions ($m$ and $n$ are integers, $n \geq 1$),

$$u=\ln(\pi\,(2\,n-1))+i\,\pi\left(\frac{1}{2}+m\right)$$

The poles are thus all situated in the half-plane $x > 0$, and are symmetric around the real axis in conjugate pairs at half-integer values of $\pi$ in the positive and negative imaginary directions. The residues of F($u$) corresponding to these poles are given by the following expression

$$\mathrm{Res}(n,m)=i\,(-1)^{m}\,(2\,n-1)^{(s-1)}\,\pi^{(s-1)}\,e^{(i\,(1/2+m)\,s\,\pi)}$$

Let $S_N$ be the sum of the residues in the strip $0 < y < 2\,\pi$ (i e for $m = 0$ and $m = 1$), and from $n = 1$ up to and including the pair of residues at $x = \ln((2N\text{-}1)\,\pi)$. Then

$$S_{N}=2\,\sin\left(\frac{1}{2}\,s\,\pi\right)e^{(i\,s\,\pi)}\,\pi^{(s-1)}\left(\sum_{n=1}^{N}(2\,n-1)^{(s-1)}\right) \tag{6}$$

(6), as well as (7) below, both tend to infinity with $N$ for $0 < Re(s) < 1$. Note however Remark 5.2 in Sect 5 below.

## 4. Contour integral

Consider now a closed contour $C_N$ in the complex plane (see Sect A1 in Appendix A), consisting of the real axis in the positive direction from $x = -\infty$ to $x = L$ just to the right of the pair of residues mentioned above at $x = \ln((2N\text{-}1)\,\pi)$, then a vertical connection from $y = 0$ to $y = 2\pi$ at $x = L$ up to a line from $x = L$ back to $x = -\infty$ in the negative direction parallel to the real axis and at a distance $2\pi$ above it, and then finally a vertical connection at negative infinity back down from $y = 2\pi$ to $y = 0$. This contour encloses the $N$ pairs of residues summed as $S_N$ in (6) above, and is here traversed in the positive direction.

THEOREM 4.1. *The integral* $I_N$ *of* (5) *around the contour* $C_N$ *as defined above is*

$$I_N = \frac{i\,2^{(s+1)}\,N^s\,\pi^s\,e^{(i s \pi)}\,\sin\left(\frac{1}{2}s\,\pi\right)E_N(s)}{s} - 2\,i\,e^{(i s \pi)}\sin(s\,\pi)\,(1-2^{(1-s)})\,\Gamma(s)\,\zeta(s) \qquad (7)$$

*where* $E_N(s)$ *is an error function incorporating truncation errors.*

*Proof.* See Appendix A.

## 5. Two equivalent expressions for ζ(*s*)

Now use Cauchy's theorem to equate the contour integral $I_N$ in (7) to the sum of residues $S_N$ in (6),

$$I_N = 2\,i\,\pi\,S_N$$

i e

$$\frac{i\,2^{(s+1)}\,N^s\,\pi^s\,e^{(i s \pi)}\,\sin\left(\frac{1}{2}s\,\pi\right)E_N(s)}{s} - 2\,i\,e^{(i s \pi)}\sin(s\,\pi)\,(1-2^{(1-s)})\,\Gamma(s)\,\zeta(s) = 4\,i\sin\left(\frac{1}{2}s\,\pi\right)\pi^s\left(\sum_{n=1}^{N}(2n-1)^{(s-1)}\right)e^{(i s \pi)} \qquad (8)$$

Solve for $\zeta(s)$,

$$\zeta(s) = \frac{\pi^s\left(-2^{(s-1)}\,N^s\,E_N(s) + \left(\sum_{n=1}^{N}(2\,n-1)^{(s-1)}\right)s\right)}{\cos\left(\frac{1}{2}s\,\pi\right)(-1+2^{(1-s)})\,\Gamma(s+1)} \qquad (9)$$

An equivalent expression can be obtained by making the substitution $s \to 1-s$ in (9),

$$\zeta(1-s) = \frac{\pi^{(1-s)}\left(-2^{(-s)}\,N^{(1-s)}\,E_N(1-s) + \left(\sum_{n=1}^{N}(2\,n-1)^{(-s)}\right)(1-s)\right)}{\sin\left(\frac{1}{2}s\,\pi\right)(-1+2^s)\,\Gamma(2-s)} \qquad (10)$$

and then transforming back to $\zeta(s)$ again by using the functional equation (2),

$$\zeta(s) = \frac{-N^{(1-s)}\,E_N(1-s) + 2^s\left(\sum_{n=1}^{N}(2n-1)^{(-s)}\right)(1-s)}{(-1+2^s)(1-s)} \qquad (11)$$

From (A10) in Appendix A, the error functions in (9) and (11) can be written as follows [regarding $O(1/N^{3.})$, see Sect A3.13, paragraph 2, in Appendix A],

$$E_N(s) = 1 + \frac{s\,(s-1)\,\varepsilon_\nu(s)}{N^2} + O\left(\frac{1}{N^{3.}}\right) \qquad (12a)$$

$$E_N(1-s) = 1 + \frac{s\,(s-1)\,\varepsilon_\nu(1-s)}{N^2} + O\left(\frac{1}{N^{3.}}\right) \qquad (12b)$$

*Remark 5.1.* The two equivalent expressions (9) and (11) above are somewhat analogous to the two equivalent expressions obtained by an integral and the same integral integrated by parts. In fact, performing the analogous operations as above on, e g, Euler's integral form [1] of the related gamma-function [in that case the substitution $s$ –> $s$ + 1 followed by the functional equation $\Gamma(s) = \Gamma(s+1)/s$ ] yields precisely the same result as integrating by parts.

*Remark 5.2.* It should be emphasized that (8), from which (9) and (11) were derived, is Cauchy's theorem, which thus rigorously connects the power $N^s$ in the first term to the zeta-function in the second term and to the sum over $N$ on the right-hand side. Since all functions involved are analytic also in the limit $N \to \infty$, this exact relationship between the terms is thus maintained to give finite results for $\zeta(s)$ also in the limit $N \to \infty$, even though the two contributions in (8) from (6) and (7) are then both divergent.

*Remark 5.3.* It is interesting to compare (11) above with the Dirichlet series valid for $\sigma > 1$ mentioned in the Introduction. Insert (12b) from above, and (A9) from Appendix A into (11),

$$(1-2^{(-s)})\,\zeta(s) = -\frac{N^{(1-s)}\left(1 - \frac{1}{24}\frac{s\,(s-1)}{N^2} + O\left(\frac{1}{N^{3.}}\right)\right)}{2^s\,(1-s)} + \left(\sum_{n=1}^{N}\frac{1}{(2n-1)^s}\right)$$

Comparing this with the Dirichlet series in the Introduction, we see that the last term on the right-hand side above is a finite form of the Dirichlet series. However, in contrast to that series, which is divergent for $\sigma < 1$, the present relationship for $\zeta(s)$ is derived from (1) and is thus valid also for $s$ in the critical strip, i e also for $0 < \sigma < 1$. This is a result of the rigorous derivation of (11) from Cauchy's theorem, and is effected by the first term on the right-hand side above tracking the behavior of the Dirichlet series as $N$ tends to infinity in order to give a correct rendering of the zeta-function.

*Remark 5.4.* By using alternative ways to extend the integrand in (3) to an analytic function on the complex plane, it is possible by the same technique as above to obtain variants of (9) and (11) [e g, by variations on the step from (3) to (4)]. Also other approaches lead to similar (but not identical) expressions for $\zeta(s)$, e g, the sum of the first $N$ terms of its Dirichlet series plus a power in $N$ as in (11) [16]. The particular variants (9) and (11) above are selected here since their properties turn out to be fortuitously well suited for the following proof of Conjecture 1.1, the Riemann zeta-hypothesis.

## 6. Functions $\zeta_N(s)'$ and $\zeta_N(s)''$

The two equivalent expressions for $\zeta(s)$ in (9) and (11) above should be understood as follows. For each $N$ there exists a particular error function $E_N(s)$ in (12a) within its Landau $O(1/N^3)$ such that (9) is exactly true. Thus for the right-hand side of (9) with this particular error function $E_N(s)$, the functional equation is exactly true also after the substitution in (10). This thus means that there exists a particular error function $E_N(1-s)$ in (12b) within its $O(1/N^3)$ such that also (11) holds exactly.

Now consider the following two functions in the range $0 < \sigma < 1$,

$$\zeta_N(s)' = \frac{\pi^s \left( -2^{(s-1)} N^s \,\Xi_N(s) + \left( \sum_{n=1}^{N} (2n-1)^{(s-1)} \right) s \right)}{\cos\left(\frac{1}{2} s \pi\right) (-1 + 2^{(1-s)}) \,\Gamma(s+1)} \tag{13}$$

$$\zeta_N(s)'' = \frac{-N^{(1-s)} \,\Xi_N(1-s) + 2^s \left( \sum_{n=1}^{N} (2n-1)^{(-s)} \right) (1-s)}{(-1 + 2^s)(1-s)} \tag{14}$$

where we let the functions $\zeta_N(s)'$ and $\zeta_N(s)''$ denote the two approximate functions around $\zeta(s)$, which we get if we set the remainders $O(1/N^3)$ in the error functions in (12a) and (12b) equal to zero instead of equal to the special remainders that correspond to the exact zeta-function as discussed above. To distinguish these two approximate error functions from the exact error functions $E_N(s)$ and $E_N(1-s)$ defined in the previous paragraph, we denote the approximate error functions we just defined as $\Xi_N(s)$ and $\Xi_N(1-s)$.

That we have the two expressions in (9) and (11) for the exact zeta-function $\zeta(s)$ is because we can write the exact zeta-function in two different but equivalent ways by using the functional equation combined with a variable transformation, as shown in Sect 5. Working with this function pair for $\zeta(s)$ is advantageous since in a simple way it automatically incorporates the functional equation with its symmetry properties into the derivation – and these symmetry properties are important for the proof.

For finite $N$, the approximate functions $\zeta_N(s)'$ and $\zeta_N(s)''$ above are of course now normally no longer equivalent, nor do they obey the functional equation.

As discussed above, the approximating functions $\zeta_N(s)'$ and $\zeta_N(s)''$ differ from the exact zeta-function $\zeta(s)$ by remainders $N^s\, O(1/N^3)$ and $N^{(1-s)}\, O(1/N^3)$, respectively. Due to the factors $N^s$ and $N^{(1-s)}$, the function $\zeta_N(s)'$ will then lie closer to $\zeta(s)$ than $\zeta_N(s)''$ does for $0 < \sigma < ½$, whereas for $½ < \sigma < 1$ the function $\zeta_N(s)''$ will lie closer to $\zeta(s)$ than $\zeta_N(s)'$ does. Nevertheless, it will be shown in Sect 7 below that the quotient $|\zeta_N(s)'/\zeta_N(s)''|$ in both cases tends to unity for all $s$ when $N \to \infty$. These two properties can be combined to a proof of Conjecture 1.1, as will be shown in Sect 9 below.

In the following analysis we shall study the quotient $\zeta_N(s)'/\zeta_N(s)''$ in the limit $N \to \infty$. As seen by writing $N^s$ as $N^\sigma e^{it\ln(N)}$, the argument of $N^s$ becomes indeterminate on the unit circle in the limit $N \to \infty$. Thus when limits of type $N^s$ are concerned, it is only relevant to consider their moduli, as we shall do in the following.

## 7. Quotient $|\zeta_N(s)'/\zeta_N(s)''|$

As discussed above, the approximating functions $\zeta_N(s)'$ and $\zeta_N(s)''$ differ from the exact zeta-function $\zeta(s)$ only by remainders $N^s\,\mathrm{O}(1/N^{3.})$ and $N^{(1\,-\,s)}\,\mathrm{O}(1/N^{3.})$, respectively. Hence we have

$$\left|\frac{\zeta_N(s)'}{\zeta_N(s)''}\right| = \left|\frac{\zeta(s)+\mathrm{O}\left(\frac{1}{N^{(3.-\sigma)}}\right)}{\zeta(s)+\mathrm{O}\left(\frac{1}{N^{(2.+\sigma)}}\right)}\right|$$

In the limit $N \to \infty$ when these remainders vanish, this quotient thus becomes equal to $|\zeta(s)/\zeta(s)| = 1$, except possibly in the case of such $s = s_0$ when $\zeta(s_0) = 0$, in which case this quotient in principle becomes an indeterminate expression of type 0/0 in the limit $N \to \infty$. However, in that case we can use l'Hôpital's rule on the above expression – possibly repeatedly if necessary. Since $\zeta(s)$ is an analytic function inside the critical strip, it can be expanded in a Taylor series, and since $\zeta(s)$ is nontrivial then at least some derivative at $\zeta(s_0) = 0$ must be nonzero. If the lowest order of such a nonzero derivative is $p$ (with $p \geq 1$), we then have (cf Appendix C for a further discussion of the double limit $N \to \infty$, $s \to s_0$)

$$\lim_{N\to\infty}\left|\frac{\zeta_N(s_0)'}{\zeta_N(s_0)''}\right| = \lim_{N\to\infty}\left|\frac{\zeta(s_0)+\mathrm{O}\left(\frac{1}{N^{(3.-\sigma)}}\right)}{\zeta(s_0)+\mathrm{O}\left(\frac{1}{N^{(2.+\sigma)}}\right)}\right| = \lim_{s\to s_0}\frac{\zeta(s)}{\zeta(s)} = \frac{\left(\frac{d^p\,\zeta(s)}{ds^p}\right)_{s=s_0}}{\left(\frac{d^p\,\zeta(s)}{ds^p}\right)_{s=s_0}} = 1 \qquad (15)$$

where the derivative is to be taken in the point $s = s_0$ where $\zeta(s_0) = 0$. We thus note that the quotient $|\zeta_N(s)'/\zeta_N(s)''|$ tends to unity for all $s$ when $N \to \infty$, as mentioned in Sect 6 above.

## 8. Quotient $|(\zeta_N(s)' - \zeta(s))/(\zeta_N(s)'' - \zeta(s))|$

As discussed above, the approximating functions $\zeta_N(s)'$ and $\zeta_N(s)''$ differ from the exact zeta-function $\zeta(s)$ by the remainders $\mathrm{O}(1/N^{(3.-\sigma)})$ and $\mathrm{O}(1/N^{(2.+\sigma)})$, respectively. Specifically, for finite $N$ the following differences between the approximate functions in (13), (14) and the exact zeta-function will then be nonzero (and unequal),

$$\zeta_N(s)' - \zeta(s) = \frac{\pi^s\left(-2^{(s-1)}N^s\,\Xi_N(s) + \left(\sum_{n=1}^{N}(2\,n-1)^{(s-1)}\right)s\right)}{\cos\left(\frac{1}{2}\,s\,\pi\right)(-1+2^{(1-s)})\,\Gamma(s+1)} - \zeta(s) \qquad (16)$$

$$\zeta_N(s)'' - \zeta(s) = \frac{-N^{(1-s)}\,\Xi_N(1-s) + 2^s\left(\sum_{n=1}^{N}(2\,n-1)^{(-s)}\right)(1-s)}{(-1+2^s)(1-s)} - \zeta(s) \qquad (17)$$

In Appendix B, the quotient of (16) and (17) is calculated in closed form to give

$$\frac{\zeta_N(s)' - \zeta(s)}{\zeta_N(s)'' - \zeta(s)} = \frac{1}{2}\frac{N^{(2s-1)}\,\pi^s\,(-4^s + 8^s)}{\cos\left(\frac{1}{2}s\,\pi\right)(-2+2^s)\,\Gamma(s-3)\,(s+2)\,(s+1)\,s} + O(N^{(2\sigma-2)}) \qquad (18)$$

For use in the proof of Conjecture 1.1 in Sect 9 below, we now give this quotient in the limit $N \to \infty$. As remarked at the end of Sect 6 above, it is then only relevant to consider its modulus, which thus becomes

$$\lim_{N\to\infty}\left|\frac{\zeta_N(s)' - \zeta(s)}{\zeta_N(s)'' - \zeta(s)}\right| = \lim_{N\to\infty}\frac{1}{2}\left|\frac{N^{(2s-1)}\,\pi^s\,(-4^s + 8^s)}{\cos\left(\frac{1}{2}s\,\pi\right)(-2+2^s)\,\Gamma(s-3)\,(s+2)\,(s+1)\,s}\right| \qquad (19)$$

The behavior of the expressions in (18) and (19) is dominated by the factor $N^{(2s-1)}$, which is central in the following proof.

## 9. Proof of Conjecture 1.1

The functions $\zeta_N(s)'$ and $\zeta_N(s)''$ above are functions of both $N$ and $s$. Rigorously, we thus need to study the quotients in Sects. 7 and 8 in the double limit when both $N \to \infty$ and $s \to s_0$, *i e* $\zeta(s) \to 0$. As discussed in Appendix C, we can choose to perform the passage to the limits as follows.

Consider first finite $N$ and let $s \to s_0$ so that $\zeta(s_0) = 0$. The values $\zeta_N(s_0)'$ and $\zeta_N(s_0)''$ differ from $\zeta(s_0)$ by the remainders, and are thus both nonvanishing for finite $N$. For $\zeta(s_0) = 0$, the quotient in (18) then becomes $\zeta_N(s_0)'/\zeta_N(s_0)''$.

Next let this quotient tend to the limit $N \to \infty$ as in (19), and equate the right-hand side of (19) to the expression (15) for $|\zeta_N(s_0)'/\zeta_N(s_0)''|$ in this limit. We then get

$$\lim_{N\to\infty}\frac{1}{2}\left|\frac{N^{(2s-1)}\,\pi^s\,(4^s - 8^s)}{\cos\left(\frac{1}{2}s\,\pi\right)(-2+2^s)\,(s+2)\,(s+1)\,s\,\Gamma(s-3)}\right| = 1 \qquad (20)$$

This equation can be true only if the modulus of $N^{(2s-1)}$ is equal to $N^0 = 1$, which thus requires that

$$\sigma = \frac{1}{2}$$

This thus proves Conjecture 1.1 that $Re(s)$ must be equal to ½ for all zeros of the Riemann zeta-function $\zeta(s)$ in the range $0 < Re(s) < 1$.

*Remark 9.1.* It should be noted that also other values of $\zeta(s)$ than $\zeta(s) = 0$ can make the quotient on the left-hand side in (19) become unity and give (20), and thus the above value of $\sigma$. However, among the $\zeta(s)$ that have this property, we must necessarily find also every nontrivial zero of the zeta-function, as was shown above (cf Remark B.2 in Appendix B).

*Remark 9.2.* Parenthetically, we note that for consistency the rest of the expression on the left-hand side of (20) should also become unity for $\sigma = ½$. Since for $s = ½ + i\,t$ we have $|(s+2)\,(s+1)\,s\,\Gamma(s-3)| = |\Gamma(s)|$, then the left-hand side of (20) can be calculated as follows for $s = ½ + i\,t$, where the last equality is a known property [2] of the gamma-function,

$$\frac{1}{2}\left|\frac{\pi^s\,(4^s-8^s)}{\cos\left(\frac{1}{2}\,s\,\pi\right)(-2+2^s)\,\Gamma(s)}\right| = \frac{\sqrt{\dfrac{\pi}{\cosh(\pi\,t)}}}{\left|\Gamma\left(\frac{1}{2}+i\,t\right)\right|} = 1$$

## Appendix A: Proof of Theorem 4.1

*A1. Contour integral $I_N$.* The integral $I_N$ of (5) around the closed contour $C_N$ defined in Sect 4 can be written as follows,

$$I_N = \int_{-\infty}^{L} \mathrm{F}(x)\,dx + i\int_{0}^{2\pi} \mathrm{F}(L+i\,y)\,dy - \int_{-\infty}^{L} \mathrm{F}(x+2\,i\,\pi)\,dx - i\int_{0}^{2\pi} \mathrm{F}(-\infty+i\,y)\,dy$$

Here the first term in $I_N$ is the (transformed) Riemann integral with finite upper limit $x = L$,

$$\int_{-\infty}^{L} \mathrm{F}(x)\,dx = \int_{-\infty}^{L} \frac{e^{(s\,x)}}{e^{(e^x)}+1}\,dx$$

The second term in $I_N$ is the integral of (5) from $y = 0$ to $y = 2\,\pi$ for $x = L$,

$$i\int_{0}^{2\pi} \mathrm{F}(L+i\,y)\,dy = i\int_{0}^{2\pi} \frac{e^{(s\,(L+i\,y))}}{e^{(e^{(L+i\,y)})}+1}\,dy$$

The third term in $I_N$, i e the integral of (5) from $x = L$ down to $x \to -\infty$ along $y = 2\,\pi$, can be shown to be a factor times the first term above,

$$-\int_{-\infty}^{L} \mathrm{F}(x+2\,i\,\pi)\,dx = -\int_{-\infty}^{L} \frac{e^{(s\,(x+2\,i\,\pi))}}{e^{(e^{(x+2\,i\,\pi)})}+1}\,dx$$

i e

$$-\int_{-\infty}^{L} \mathrm{F}(x+2\,i\,\pi)\,dx = -e^{(2\,i\,s\,\pi)}\int_{-\infty}^{L} \frac{e^{(s\,x)}}{e^{(e^x)}+1}\,dx$$

The fourth term in $I_N$, i e the vertical connection from $y = 2\pi$ down to $y = 0$ at $x \to -\infty$, can be shown to tend to zero (for $\sigma > 0$),

$$-i \int_0^{2\pi} \mathrm{F}(-\infty + i\,y)\,dy = 0$$

As shown in Sect 3, the integrand has poles at certain values of $x$ and $y$. In the following, $L$ is assumed to be chosen to stay clear of those poles, or specifically ($N$ is an integer, $N \geq 1$),

$$L = \ln(2\,N\,\pi)$$

Using the above results, the integral $I_N$ can be rewritten as

$$I_N = I_1 + I_2$$

where

$$I_1 = (1 - \mathrm{e}^{(2\,i\,s\,\pi)}) \int_{-\infty}^{\ln(2\,\pi\,N)} \frac{\mathrm{e}^{(s\,x)}}{\mathrm{e}^{(\mathrm{e}^x)} + 1}\,dx$$

$$I_2 = i \int_0^{2\pi} \frac{\mathrm{e}^{(s\,(\ln(2\,\pi\,N) + i\,y))}}{\mathrm{e}^{(\mathrm{e}^{(\ln(2\,\pi\,N) + i\,y)})} + 1}\,dy$$

*A2. Integral $I_1$.* The integral $I_1$ above can be rewritten as

$$I_1 = (1 - \mathrm{e}^{(2\,i\,s\,\pi)}) \int_{-\infty}^{\infty} \frac{\mathrm{e}^{(s\,x)}}{\mathrm{e}^{(\mathrm{e}^x)} + 1}\,dx - (1 - \mathrm{e}^{(2\,i\,s\,\pi)}) \int_{\ln(2\,\pi\,N)}^{\infty} \frac{\mathrm{e}^{(s\,x)}}{\mathrm{e}^{(\mathrm{e}^x)} + 1}\,dx$$

Use (4) to express the first integral in terms of the Riemann zeta-function, and transform the second integrand back to its original form,

$$I_1 = (1 - \mathrm{e}^{(2\,i\,s\,\pi)})\,(1 - 2^{(1-s)})\,\Gamma(s)\,\zeta(s) - (1 - \mathrm{e}^{(2\,i\,s\,\pi)}) \int_{2\,\pi\,N}^{\infty} \frac{w^{(s-1)}}{\mathrm{e}^w + 1}\,dw$$

For $0 < Re(s) < 1$, the last term can be estimated as follows,

$$\left| \int_{2\,N\,\pi}^{\infty} \frac{w^{(s-1)}}{\mathrm{e}^w + 1}\,dw \right| < \int_{2\,N\,\pi}^{\infty} \frac{w^0}{\mathrm{e}^w}\,dw$$

i e

$$\left| \int_{2N\pi}^{\infty} \frac{w^{(s-1)}}{e^{w}+1}\, dw \right| < e^{(-2N\pi)}$$

After converting the common factor in $I_1$ above to a sine, we thus have

$$I_1 = -2\, i\, e^{(i s \pi)} \sin(s\pi)\,(1 - 2^{(1-s)})\,\Gamma(s)\,\zeta(s) + O(e^{(-2\pi N)})$$

*A3. Integral $I_2$.* The integral $I_2$ in Sect A1 can be expanded as follows

$$I_2 = i\, 2^s \pi^s N^s \int_0^{2\pi} \frac{e^{(i s y)}}{e^{(2\pi N(\cos(y) + i \sin(y)))} + 1}\, dy \qquad (A1)$$

We will now calculate this integral by considering the corresponding integral $I_0$ of an approximating, piecewise function, which is zero for $0 \le y < \pi/2$ and for $3\pi/2 < y \le 2\pi$, whereas for $\pi/2 \le y \le 3\pi/2$ its integrand is given by the exponential in the numerator in (A1), i e

$$I_0 = i\, 2^s \pi^s N^s \int_{1/2\,\pi}^{3/2\,\pi} e^{(i s y)}\, dy \qquad (A2)$$

which integrates to

$$I_0 = \frac{2^s \pi^s N^s (e^{(3/2\, i s \pi)} - e^{(1/2\, i s \pi)})}{s} \qquad (A3)$$

For sufficiently large $N$, the integrands in $I_2$ and $I_0$ will differ appreciably only in the two regions where *cos*($y$) is close to zero, i e around $y = \pi/2$ and $y = 3\pi/2$, respectively, which points will be the centers for two corrections to the integrated result in (A3). These two corrections can be calculated by series expansions as follows around $y = \pi/2$ and $y = 3\pi/2$, respectively. The integral $I_2$ in (A1) can then be obtained by adding these two corrections (including remainders, see Sect A3.13) to the result in (A3), as will be shown in Sect A3.2 ending this Appendix.

*A3.1. Case $\pi/2$.* Study here first the behavior around $y = \pi/2$. After Taylor expansions around $\pi/2$ in the numerator and denominator in (A1), the integral $I_2$ can be written there as

$$I_2\!\left(\frac{1}{2}\pi, \delta\right) = i\, e^{(1/2\, i s \pi)}\, 2^s \pi^s N^s \int_{1/2\,\pi - \delta}^{1/2\,\pi + \delta} \frac{1 + i s\left(y - \frac{1}{2}\pi\right) - \frac{1}{2} s^2 \left(y - \frac{1}{2}\pi\right)^2 + \ldots}{e^{(2\pi N(i - y + 1/2\,\pi - 1/2\, i (y - 1/2\,\pi)^2 + 1/6 (y - 1/2\,\pi)^3 + \ldots))} + 1}\, dy$$

where the interval $\delta$ is chosen so that it (at least) covers the region of appreciable deviation from the piecewise integrand in $I_0$, as will be further discussed below.

After changing integration variable,

$$N\left(y-\frac{1}{2}\pi\right)=\tau$$

and noting that $\exp(2N\pi i)=1$, the expression above can be rewritten as follows (truncating at powers in $1/N$ of order two in the integrand, collecting the remainders into the numerator)

$$I_2\left(\frac{1}{2}\pi,\delta\right)=i\,e^{(1/2\,i\,s\,\pi)}\,2^s\,\pi^s\,N^s\int_{-N\delta}^{N\delta}\frac{\frac{1}{N}+\frac{i\,s\,\tau}{N^2}+O\left(\frac{1}{N^3}\right)}{e^{\left(-2\pi\tau-\frac{i\pi\tau^2}{N}\right)}+1}\,d\tau$$

*A3.11. Correction $\Delta I_2(\pi/2)$.* The correction $\Delta I_2(\pi/2)$ below is the correction necessary if we approximate the integral $I_2$ above around $y=\pi/2$ by the integral $I_0$ there. Similarly to $I_2$ versus $I_0$ above, it can be expressed as the above integral minus the corresponding integral with the exponential in the denominator and the lower limit both set to zero. The integration is to be performed over the interval $-N\delta<\tau<N\delta$ around $y=\pi/2$, defined as an interval that covers (at least) the region where the integrands below differ appreciably, i e more than $O(1/N^3)$.

$$\Delta I_2\left(\frac{1}{2}\pi\right)=i\,e^{(1/2\,i\,s\,\pi)}\,2^s\,\pi^s\,N^s\left(\int_{-N\delta}^{N\delta}\frac{\frac{1}{N}+\frac{i\,s\,\tau}{N^2}+O\left(\frac{1}{N^3}\right)}{e^{\left(-2\pi\tau-\frac{i\pi\tau^2}{N}\right)}+1}\,d\tau-\int_0^{N\delta}\frac{1}{N}+\frac{i\,s\,\tau}{N^2}+O\left(\frac{1}{N^3}\right)d\tau\right)$$

The size of the region where the integrands above differ appreciably as defined above, is determined by the exponential function in the denominator of the first integrand. If we set

$$N\delta=\nu=\ln(N)\qquad (A4)$$

then for $N\geq 5$ the interval $-N\delta\leq\tau\leq N\delta$ in the integrations above will with good margin include this region where the integrands differ appreciably, the margin becoming better and better as $N$ increases. At the same time, the relative proportion of the range of the integrations above compared to the original range $0<\delta<2\pi$ will become smaller and smaller as $N$ increases.

After expanding the first integrand above as a Taylor series in $1/N$, we get

$$\Delta I_2\left(\frac{1}{2}\pi\right)=i\,e^{(1/2\,i\,s\,\pi)}\,2^s\,\pi^s\,N^s\left(\int_{-N\delta}^{N\delta}\frac{a_1}{N}+\frac{a_2}{N^2}+O\left(\frac{1}{N^3}\right)d\tau-\int_0^{N\delta}\frac{1}{N}+\frac{i\,s\,\tau}{N^2}+O\left(\frac{1}{N^3}\right)d\tau\right)\qquad (A5)$$

where

$$a_1 = \frac{1}{e^{(-2\pi\tau)} + 1}$$

$$a_2 = \frac{i\,s\,\tau + \dfrac{i\,e^{(-2\pi\tau)}\,\pi\,\tau^2}{e^{(-2\pi\tau)} + 1}}{e^{(-2\pi\tau)} + 1}$$

The correction $\Delta I_2(\pi/2)$ in (A5) corresponds to the integral of a more or less sharp peak (depending on $N$) at $y = \pi/2$, i e at $\tau = 0$. This peak is well described by the functions of $\tau$ in $a_1$, $a_2$, and in the second integrand in (A5) above. The series expansions in $1/N$ in (A5) are required only to second order for the following calculations.

*A3.12. Factor $\varepsilon_\nu$.* Inserting the above expressions for $a_1$ and $a_2$ into (A5), and integrating for each power of $N$, we get

$$\int_{-N\delta}^{0} \frac{a_1}{N}\,d\tau + \int_{0}^{N\delta} \frac{a_1}{N} - \frac{1}{N}\,d\tau = 0 \qquad (A6)$$

$$\int_{-N\delta}^{0} \frac{a_2}{N^2}\,d\tau + \int_{0}^{N\delta} \frac{a_2}{N^2} - \frac{i\,s\,\tau}{N^2}\,d\tau = \frac{i\,(s-1)\,\varepsilon_\nu(s)}{N^2} \qquad (A7)$$

where the factor $\varepsilon_\nu(s) = \varepsilon(s, N\delta)$ is a function of $s$ and the limit $N\delta = \nu$, and can be expressed in exact, explicit form as follows using the dilogarithm function $Li_2(x)$ as defined in [3] (NB other definition in [9]),

$$\varepsilon_\nu(s) = \frac{1}{4}\,\frac{Li_2(e^{(2\nu\pi)}+1) - Li_2(e^{(-2\nu\pi)}+1)}{\pi^2} + \frac{\ln(e^{(2\nu\pi)}+1)\,\nu}{\pi} - \frac{3}{2}\nu^2 + \frac{\nu^2}{(e^{(2\nu\pi)}+1)(1-s)} \qquad (A8)$$

For large $\nu$, the function $\varepsilon_\nu(s)$ in (A8) can be written as the following series expansion,

$$\varepsilon_\nu(s) = -\frac{1}{24} + \left(\frac{1}{2}\,\frac{1}{\pi^2} + \frac{\nu}{\pi} + \frac{\nu^2}{1-s}\right) e^{(-2\nu\pi)} + \mathrm{O}(\nu^2\, e^{(-4\nu\pi)})$$

giving the following asymptotic value for $\nu >> 1$,

$$\varepsilon_\nu(s) = -\frac{1}{24} + \mathrm{O}(\nu^2\, e^{(-2\nu\pi)})$$

Inserting $\nu$ from (A4) we get

$$\varepsilon_\nu(s) = -\frac{1}{24} + \mathrm{O}\left(\frac{\ln(N)^2}{N^{(2\pi)}}\right) \qquad (A9)$$

*A3.13. **Remainders.*** The remainder in the factor $\varepsilon_\nu(s)$ in (A9) above is one order smaller than the remainders in the final results in Sect A3.2 below, and thus negligible in comparison. However, when we later calculate the series expansions with remainders $O(N^{(\sigma-6)})$ and $O(N^{(-\sigma-5)})$, respectively, in Appendix B, the terms of lower order vanish identically, so in principle the remainder in (A9) could then potentially be left over (with a factor $N^s$ and $N^{(1-s)}$, respectively, in front). But even in that case the remainder from (A9) would still be negligible compared to said remainders (the negative power in the denominator above being -2π compared to at the most -6 in the remainders in Appendix B). For these reasons, we can for simplicity disregard the remainder in (A9) in all the following calculations.

When integrating the remainders $O(1/N^3)$ in (A5) over the interval $2N\delta = 2\ln(N)$, the integrated remainder becomes of type $O(\ln(N)^p/N^3)$, where $p$ is a positive number, in this case equal to three. Note that any remainder of this type, although greater than $O(1/N^3)$, is for sufficiently large $N$ still always smaller than any remainder $O(1/N^{3-|\varepsilon|})$, no matter how small $|\varepsilon|$ may be. Throughout this paper we denote for simplicity a remainder of this approximate power type as $O(1/N^{3.})$, i e with a decimal point in the power to signify that the power is not an exact integer.

For completeness, we also need to estimate the error we make in $\Delta I_2(\pi/2)$ when we neglect the rest of the integral outside the interval $-N\delta \le \tau \le N\delta$, which error can be written

$$\Delta\left(\Delta I_2\left(\frac{1}{2}\pi\right)\right) = i\,(2\,\pi\,N)^s \left(\int_0^{1/2\,\pi-\delta} \mathrm{h}(y, s, N)\,dy + \int_{1/2\,\pi+\delta}^{\pi} \mathrm{h}(y, s, N) - \mathrm{g}(y, s)\,dy\right)$$

where

$$\mathrm{h}(y, s, N) = \frac{e^{(i\,s\,y)}}{e^{(2\,\pi\,N\,(\cos(y) + i\,\sin(y)))} + 1} \qquad\qquad \mathrm{g}(y, s) = e^{(i\,s\,y)}$$

and which can be shown to be

$$\Delta\left(\Delta I_2\left(\frac{1}{2}\pi\right)\right) < O\left(\frac{1}{N^4}\right)$$

i e it is thus negligible compared to the remainders $O(1/N^{3.})$ below [as expected, since this is the part of the integral in (A5) outside the region where the integrands differ appreciably as defined above].

*A3.2. **Final results.*** Inserting (A6) and (A7) into (A5) we thus finally get (note that there are no terms of zeroth and first order in $1/N$)

$$\Delta I_2\left(\frac{1}{2}\pi\right) = i\,e^{(1/2\,i\,s\,\pi)}\,2^s\,\pi^s\,N^s\left(\frac{i\,\varepsilon_\nu(s)\,(s-1)}{N^2} + O\left(\frac{1}{N^{3.}}\right)\right)$$

For the region around $3\pi/2$ we can similarly calculate the following correction (note again that there are no terms of zeroth and first order in $1/N$),

$$\Delta I_2\left(\frac{3}{2}\pi\right) = i\,e^{(3/2\,i\,s\,\pi)}\,2^s\,\pi^s\,N^s\left(-\frac{i\,\varepsilon_\nu(s)\,(s-1)}{N^2} + O\left(\frac{1}{N^{3.}}\right)\right)$$

Adding the two corrections above to the integral $I_0$ in (A3) gives

$$I_2 = \frac{2^s \pi^s N^s \left(-e^{(1/2\, i\, s\, \pi)} + e^{(3/2\, i\, s\, \pi)}\right)\left(1 + \frac{s\,(s-1)\,\varepsilon_v(s)}{N^2} + O\left(\frac{1}{N^{3.}}\right)\right)}{s}$$

Simplify the sum of exponentials to a sine function, and define the error function $E_N(s)$ as follows

$$I_2 = \frac{2\, i\, 2^s N^s \pi^s e^{(i\, s\, \pi)} \sin\left(\frac{1}{2} s\, \pi\right) E_N(s)}{s}$$

$$E_N(s) = 1 + \frac{s\,(s-1)\,\varepsilon_v(s)}{N^2} + O\left(\frac{1}{N^{3.}}\right) \qquad \textit{(A10)}$$

where $\varepsilon_v(s)$ is given in (A8) and (A9) above and the notation $O(1/N^{3.})$ is explained above in Remark A3.13, paragraph 2.

Combining the above result for the integral $I_2$ with the result in Sect A2 for the integral $I_1$, the contour integral $I_N$ can thus finally be written as follows [the remainder in $I_1$ is negligible compared to the remainder in $I_2$ from $E_N(s)$],

$$I_N = \frac{i\, 2^{(s+1)} N^s \pi^s e^{(i\, s\, \pi)} \sin\left(\frac{1}{2} s\, \pi\right) E_N(s)}{s} - 2\, i\, e^{(i\, s\, \pi)} \sin(s\, \pi)\,(1 - 2^{(1-s)})\,\Gamma(s)\,\zeta(s) \qquad \textit{(A11)}$$

(A11), (A10), and (A9) derived above thus prove Theorem 4.1 in Sect 4.

## Appendix B: Calculation of $\zeta_N$ - $\zeta$ in closed form

Eqs (13) and (14) should hold (with other error functions) also for $N$ –> $N+1$, i e

$$\zeta_{N+1}(s)' = \frac{\pi^s \left(-2^{(s-1)} (N+1)^s\, \Xi_{N+1}(s) + \left(\sum_{n=1}^{N+1} (2\,n-1)^{(s-1)}\right) s\right)}{\cos\left(\frac{1}{2} s\, \pi\right)(-1 + 2^{(1-s)})\,\Gamma(s+1)}$$

$$\zeta_{N+1}(s)'' = \frac{-(N+1)^{(1-s)}\, \Xi_{N+1}(1-s) + 2^s \left(\sum_{n=1}^{N+1} (2\,n-1)^{(-s)}\right)(1-s)}{(-1+2^s)\,(1-s)}$$

where the notation $\zeta_{N+1}(s)$ denotes that these expressions for the Riemann zeta-function correspond to truncation at order $N+1$ in the above series.

Calculate the difference between (13) and (14), respectively, and the above relationships,

$$\zeta_N(s)' - \zeta_{N+1}(s)' = \frac{\pi^s \left(-2^{(s-1)} N^s\, \Xi_N(s) + 2^{(s-1)} (N+1)^s\, \Xi_{N+1}(s) - s\,(1+2N)^{(s-1)}\right)}{\cos\left(\frac{1}{2} s\,\pi\right) (-1 + 2^{(1-s)})\, \Gamma(s+1)}$$

$$\zeta_N(s)'' - \zeta_{N+1}(s)'' = \frac{-N^{(1-s)}\, \Xi_N(1-s) + (N+1)^{(1-s)}\, \Xi_{N+1}(1-s) - 2^s (1-s)(1+2N)^{(-s)}}{(1-s)(-1+2^s)}$$

*Remark B.1.* It should be remembered that the error function $E_N$ was calculated in (A10) in Appendix A using the explicit calculation of the error $\varepsilon_\nu(s)$ given in (A8) and (A9), and was then used in (7) etc through to (12a) and (12b). The error function $E_N$ was then subsequently used to define $\Xi_N$ by setting the remainder $O(1/N^{3\cdot})$ equal to zero at (13) and (14) onwards, via (16) and (17), and through the rest of the present Appendix B.

The error functions $\Xi_N(s)$ and $\Xi_N(1\text{-}s)$ are thus given by, respectively, (12a) and (12b) in the special case when we set the remainders $O(1/N^{3\cdot})$ equal to zero there, so the above relationships can hence be written explicitly as follows

$$\zeta_N(s)' - \zeta_{N+1}(s)' =$$
$$= \frac{\pi^s \left(-2^{(s-1)} N^s \left(1 - \frac{1}{24}\frac{s\,(s-1)}{N^2}\right) + 2^{(s-1)} (N+1)^s \left(1 - \frac{1}{24}\frac{s\,(s-1)}{(N+1)^2}\right) - s\,(1+2N)^{(s-1)}\right)}{\cos\left(\frac{1}{2} s\,\pi\right)(-1+2^{(1-s)})\,\Gamma(s+1)}$$

$$\zeta_N(s)'' - \zeta_{N+1}(s)'' =$$
$$= \frac{-N^{(1-s)}\left(1 - \frac{1}{24}\frac{s\,(s-1)}{N^2}\right) + (N+1)^{(1-s)}\left(1 - \frac{1}{24}\frac{s\,(s-1)}{(N+1)^2}\right) - 2^s (1-s)(1+2N)^{(-s)}}{(1-s)(-1+2^s)}$$

In order to study these two relationships, we make Taylor expansions in $1/N$ of the factors $(N+1)$ and $(1+2N)$ with their different exponents. By their very nature, the two differences above are very small. As a consequence, the handling of the terms and the resulting truncation errors in the above two relationships require series expansions of high order to get non-vanishing results. [For the same reason, numerical calculations of the right-hand sides above have to be made using high accuracy (30 digits or more) in order to give correct results].

The leading terms and remainders in the Taylor expansions of the two relationships above can be shown to be as follows (after some calculation),

$$\zeta_N(s)' - \zeta_{N+1}(s)' = \frac{7}{11520}\frac{\pi^s\, 4^s\, N^{(s-5)}}{\cos\left(\frac{1}{2} s\,\pi\right)(-2+2^s)\,\Gamma(s-4)} + O(N^{(\sigma-6)})$$

$$\zeta_N(s)'' - \zeta_{N+1}(s)'' = -\frac{7}{5760}\frac{(s+3)(s+2)(s+1)\,N^{(-s-4)}\, s}{-1+2^s} + O(N^{(-\sigma-5)})$$

We now wish to reformulate these expressions in terms of the zeta-function itself rather than in $\zeta_{N+1}(s)$. We begin by considering first-order Taylor expansions in $1/N$ as follows

$$N^{(s-4)} - (N+1)^{(s-4)} = -N^{(s-5)}(s-4) + O(N^{(\sigma-6)})$$

$$N^{(-3-s)} - (N+1)^{(-3-s)} = -N^{(-s-4)}(-3-s) + O(N^{(-\sigma-5)})$$

and divide the left-hand sides of the two previous results by the left-hand sides of the respective Taylor expansion, and similarly for the right-hand sides,

$$\frac{\zeta_N(s)' - \zeta_{N+1}(s)'}{N^{(s-4)} - (N+1)^{(s-4)}} = -\frac{7}{11520}\frac{\pi^s 4^s}{\cos\left(\frac{1}{2}s\pi\right)(-2+2^s)\,\Gamma(s-3)} + O\left(\frac{1}{N}\right)$$

$$\frac{\zeta_N(s)'' - \zeta_{N+1}(s)''}{N^{(-3-s)} - (N+1)^{(-3-s)}} = -\frac{7}{5760}\frac{(s+2)(s+1)s}{-1+2^s} + O\left(\frac{1}{N}\right)$$

For use in Sect 8, we equate the quotient of the left-hand sides above to the quotient of the right-hand sides,

$$\frac{\dfrac{\zeta_N(s)' - \zeta_{N+1}(s)'}{N^{(s-4)} - (N+1)^{(s-4)}}}{\dfrac{\zeta_N(s)'' - \zeta_{N+1}(s)''}{N^{(-3-s)} - (N+1)^{(-3-s)}}} = \frac{\dfrac{7}{11520}\dfrac{\pi^s 4^s}{\cos\left(\frac{1}{2}s\pi\right)(-2+2^s)\,\Gamma(s-3)}}{\dfrac{7}{5760}\dfrac{(s+2)(s+1)s}{-1+2^s}} + O(1/N)$$

The right-hand side above can be simplified to

$$RHS = \frac{1}{2}\frac{\pi^s(-4^s+8^s)}{\cos\left(\frac{1}{2}s\pi\right)\Gamma(s-3)(s+2)(s+1)s(-2+2^s)} + O\left(\frac{1}{N}\right)$$

For sufficiently large $N$, the right-hand side ***RHS*** is thus a nonzero constant, i e independent of $N$. This thus means that the corresponding left-hand side (***LHS***) above is the same [within $O(1/N)$] even if $N+1$ is replaced by $N+k$, where $k$ is an arbitrary positive integer, and where furthermore nothing prevents us from letting $k$ tend to infinity. In the limit $k \to \infty$, we have $\zeta_{N+k}(s) = \zeta(s)$ as discussed in Sect 6, and also $(N+k)^{(s-4)} = 0$ and $(N+k)^{(-3-s)} = 0$ (since we assume $0 < \sigma < 1$). The following equalities thus hold within $O(1/N)$,

$$LHS = \frac{\dfrac{\zeta_N(s)' - \zeta_{N+1}(s)'}{N^{(s-4)} - (N+1)^{(s-4)}}}{\dfrac{\zeta_N(s)'' - \zeta_{N+1}(s)''}{N^{(-3-s)} - (N+1)^{(-3-s)}}} = \frac{\dfrac{\zeta_N(s)' - \zeta_{N+k}(s)'}{N^{(s-4)} - (N+k)^{(s-4)}}}{\dfrac{\zeta_N(s)'' - \zeta_{N+k}(s)''}{N^{(-3-s)} - (N+k)^{(-3-s)}}} = \frac{\dfrac{\zeta_N(s)' - \zeta(s)}{N^{(s-4)}}}{\dfrac{\zeta_N(s)'' - \zeta(s)}{N^{(-3-s)}}}$$

Equating the last member of these left-hand sides *LHS* to the right-hand side *RHS* given above, we obtain

$$\frac{\zeta_N(s)' - \zeta(s)}{\zeta_N(s)'' - \zeta(s)} = \frac{1}{2}\,\frac{N^{(2s-1)}\,\pi^s\,(-4^s + 8^s)}{\cos\left(\frac{1}{2}s\,\pi\right)(-2+2^s)\,\Gamma(s-3)\,(s+2)\,(s+1)\,s} + O(N^{(2\sigma-2)}) \quad (18)$$

which is used as (18) in the above proof of Conjecture 1.1.

*Remark B.2.* Note that in the limit $N \rightarrow \infty$, the right-hand side of (18) becomes zero for $\sigma < ½$ and tends to infinity for $\sigma > ½$. In the limit $N \rightarrow \infty$, it has a finite, nonzero value only for $\sigma = ½$, when the modulus of the right-hand side |*RHS*| of (18) becomes unity (cf Remark 9.2). Corresponding to this case when $\sigma = ½$ and |*RHS*| is unity, the modulus of the left-hand side |*LHS*| of (18) can become unity in the limit $N \rightarrow \infty$ in the following two different but partially overlapping cases:

(*a*) for $\sigma = ½$ and $\zeta(s) = 0$, in which case |*LHS*| becomes $|\zeta_N(s)'|/|\zeta_N(s)''| \rightarrow 1$,

(*b*) for $\sigma = ½$, in which case |*LHS*| becomes $|\zeta_N(s)' - \zeta(s)|/|\zeta_N(s)'' - \zeta(s)| \rightarrow 1$.

A necessary and sufficient condition for the modulus of the left-hand side of (18) to tend to unity in the limit $N \rightarrow \infty$ is thus that $\sigma = ½$. According to case (*b*) the modulus of the left-hand side of (18) tends to unity for <u>all</u> $t$ in $s = ½ + i\,t$, whereas according to case (*a*) it does so specifically for those $t$ in $s = ½ + i\,t$ for which $\zeta(s) = 0$. The latter, special case (*a*) thus reiterates the proof in Sect 6 that $Re(s) = ½$ for every (nontrivial) zero $\zeta(s) = 0$ of the Riemann zeta-function.

## Appendix C: Double limit $N \rightarrow \infty$, $\zeta(s) \rightarrow 0$

The functions $\zeta_N(s)'$ and $\zeta_N(s)''$ are functions of both $N$ and $s$. We thus need to study the quotients in Sects. 7 and 8 above in the double limit when both

(a) $N \rightarrow \infty$, and

(b) $s \rightarrow s_0$, so that $\zeta(s) \rightarrow 0$.

This double limit is governed by the following theorem on uniform convergence of double limits [8]

THEOREM C: *If the limit lim $a_{nm} = l_m$ for $n \to \infty$ exists uniformly with respect to m, and if further the limit lim $l_m = l$ for $m \to \infty$ exists, then the double limit lim $a_{nm}$ for $n \to \infty$ and $m \to \infty$ exists and has the value l. We can then reverse the order of the passages to the limit, provided that lim $a_{nm} = \lambda_n$ for $m \to \infty$ exists.*

Below in Theorems C.I and C.II are formulations of the above theorem for the two cases we study above in Sects. 7 and 8, respectively. Note that if both limits in a double limit exist, then only one of them needs to be uniform [6] according to Theorem C.

THEOREM C.I (re Sect. 7): *If the limit*

$$L_{\zeta(s)} = \lim_{N \to \infty} \left| \frac{\zeta_N(s)'}{\zeta_N(s)''} \right| = \lim_{N \to \infty} \left| \frac{\zeta(s) + O(N^{(-3. + \sigma)})}{\zeta(s) + O(N^{(-2. - \sigma)})} \right| \qquad (\alpha)$$

*i e*

$$L_{\zeta(s)} = \left| \frac{\zeta(s)}{\zeta(s)} \right| \qquad (\beta)$$

*exists uniformly with respect to* $\zeta(s)$, *and if further the limit*

$$\lim_{\zeta(s) \to 0} L_{\zeta(s)} = \lim_{\zeta(s) \to 0} \left| \frac{\zeta(s)}{\zeta(s)} \right| \qquad (\gamma)$$

*exists, and has the value*

$$\lim_{\zeta(s) \to 0} \left| \frac{\zeta(s)}{\zeta(s)} \right| = 1 \qquad (\delta)$$

*then the double limit*

$$\lim_{(\zeta(s), N) \to (0, \infty)} \left| \frac{\zeta(s) + O(N^{(-3. + \sigma)})}{\zeta(s) + O(N^{(-2. - \sigma)})} \right| = \lim_{(\zeta(s), N) \to (0, \infty)} \left| \frac{\zeta_N(s)'}{\zeta_N(s)''} \right| \qquad (\varepsilon)$$

*exists, and has the value*

$$\lim_{(\zeta(s), N) \to (0, \infty)} \left| \frac{\zeta_N(s)'}{\zeta_N(s)''} \right| = 1 \qquad (\zeta)$$

*We can the reverse the order of the passages to the limit, provided that the limit*

$$L_N = \lim_{\zeta(s) \to 0} \left| \frac{\zeta_N(s)'}{\zeta_N(s)''} \right| \qquad (\alpha')$$

*exists.*

***Remark C.I:*** The step from (α) and (β) above to (γ) is valid at least for $\zeta(s) \neq 0$ as discussed in Sect. 7 above, where in (15) also the limit in (γ) is shown to exist and have the value in (δ). Hence Theorem C.I is applicable, and thus the double limit in (ε) exists and has the value in (ζ), i e unity. We can then reverse the order of the passages to the limit since from Sects. 6 and 7 follows that the limit in (α′) exists.

THEOREM C.II (re Sect 8): *Since the limit*

$$\Lambda_N = \lim_{\zeta(s)\to 0} \left| \frac{\zeta_N(s)' - \zeta(s)}{\zeta_N(s)'' - \zeta(s)} \right| \qquad (\theta)$$

*i e for s –> $s_0$ where $\zeta(s_0)$ = 0 (and with the right-hand side here taken, e g, from (18) in Sect. 8),*

$$\Lambda_N = \left| \frac{\zeta_N(s_0)'}{\zeta_N(s_0)''} \right| = \lim_{s\to s_0} \left| \frac{1}{2} \frac{N^{(2s-1)}\,\pi^s\,(-4^s+8^s)}{\cos\left(\frac{1}{2}s\,\pi\right)(-2+2^s)\,\Gamma(s-3)\,(s+2)\,(s+1)\,s} + O(N^{(2\sigma-2)}) \right| \qquad (\iota)$$

*exists uniformly with respect to N for $0 < \sigma \leq ½$ (for $½ < \sigma < 1$ see Remark C.II below), and since further the limit*

$$\lim_{N\to\infty} \Lambda_N = \lim_{N\to\infty} \left| \frac{\zeta_N(s_0)'}{\zeta_N(s_0)''} \right| \qquad (\kappa)$$

*according to* ($\varepsilon$) *and* ($\zeta$) *in Theorem C.I exists and has the value*

$$\lim_{N\to\infty} \left| \frac{\zeta_N(s_0)'}{\zeta_N(s_0)''} \right| = 1 \qquad (\lambda)$$

*then the double limit*

$$\lim_{(\zeta(s),\,N)\to(0,\,\infty)} \left| \frac{\zeta_N(s)' - \zeta(s)}{\zeta_N(s)'' - \zeta(s)} \right| \qquad (\mu)$$

*exists, and has the value*

$$\lim_{(\zeta(s),\,N)\to(0,\,\infty)} \left| \frac{\zeta_N(s)' - \zeta(s)}{\zeta_N(s)'' - \zeta(s)} \right| = 1 \qquad (\nu)$$

***Remark C.II:*** The step from ($\theta$) to ($\iota$) is here shown using (18). The right-hand side, as given in ($\iota$), in this case exists uniformly with respect to $N$ for $0 < \sigma \leq ½$. Hence the step ($\kappa$) to ($\lambda$) then follows directly from Theorem C.I above. For $½ < \sigma < 1$ we use the substitution $s$ -> 1 - $s$ to convert ($\iota$) to an equation involving $\zeta_N(1 - s_0)'/\zeta_N(1 - s_0)''$ on the left hand side, and on the right-hand side $N^{(1-2s)}$ times a factor that is a mirror image of the right-hand side around $\sigma = ½$ (and now with a remainder $O(N^{-2\sigma})$). Hence this right-hand side then exists uniformly with respect to $N$ for $½ \leq \sigma < 1$. Since $|\zeta_N(s)'| = |\zeta_N(s)''| = |\zeta(s)|$ in the limit $N \to \infty$, the quotient on the left-hand side then becomes $|\zeta_N(1 - s_0)'/\zeta_N(1 - s_0)''| = |\zeta_N(s_0)'/\zeta_N(s_0)''| = 1$ as in ($\lambda$) according to Theorem C.I. Thus Theorem C.II is applicable for all $\sigma$, and the limit in ($\mu$) exists and has the value in ($\nu$), *i e* unity, as required in the proof of RH in Sect 9.

Theorems C.I and C.II thus validate the passages to the limits in Sects. 7 - 9 in the proof of Conjecture 1.1.

*Acknowledgments.* The author wishes to express deep gratitude to Professor Sir Michael Berry for his friendly and encouraging criticism with regard to a very early version of the manuscript, and to Dr Hans-Olov Zetterström for many invaluable, clarifying discussions in the initial stages of this project. My profound thanks also go to the many readers who have critically studied the successive updates of my manuscript arXiv:0809.5120 [math.GM] on http://arXiv.org/ and asked all the relevant questions that are now summarized in the extensive collection of FAQ also filed there, with special thanks to Mr Hisanobu Shinya for meticulously checking my calculations over a period of many months, leading to several improvements in the text. My special thanks also go to Professor Yasuyuki Kawahigashi and Professor Takayuki Oda for kindly and helpfully making me aware of the need for the detailed discussion now made in the present Sects 7 and 8. I am also grateful to Dr Sergey K Sekatskii for interesting comments in connection with his work on applying my exponential transformation in Sect 2 to other transcendental functions (cf arXiv:0912.4116 [math.FA]).

B&E Scientific Ltd, Seaford, BN25 4PA, United Kingdom
*E-mail address:* arne.bergstrom@physics.org

## References

[1] M. Abramowitz and I. M. Stegun (Eds.), *Handbook of Mathematical Functions*, Dover (1972), eq. 6.1.1.
[2] ——, eq. 6.1.30.
[3] ——, p. 1004.
[4] M. V. Berry and J.P. Keating, SIAM Review 41, No. 2 (1999), 236-266.
[5] ——, *Supersymmetry and trace formulae: chaos and disorder*, Eds I. V. Lerner, J. P. Keating, Plenum, New York, 1999, pp. 355-367.
[6] K. G. Binmore, *Foundations of Analysis 2: Topological Ideas*, Cambridge University Press (2008), p.136.
[7] H. M. Edwards, *Riemann's Zeta Function*, Dover (2001), p. 13.
[8] R. Courant, *Differential and Integral Calculus II*, Wiley (1988), p 104-105.
[9] A. Erdélyi (Ed.), *Higher Transcendental Functions*, Vol. 1, McGraw-Hill (1953), p 31.
[10] E. Jahnke and F. Emde, *Tables of Functions*, 4th Ed, Dover (1945), p. 269.
[11] B. Riemann, Monatsberichte der Berliner Akademie, Nov 1859, http://www.claymath.org/millennium/Riemann_Hypothesis/1859_manuscript/
[12] M. du Sautoy, *The Music of the Primes*, Fourth Estate (2003), p. 1.
[13] E. C. Titchmarsh (revised by D. R. Heath-Brown), *The Theory of the Riemann Zeta-function*, 2nd Ed, Clarendon Press (2003), p. 16.
[14] ——, p. 21.
[15] ——, p. 45.
[16] ——, p. 77.
[17] E. W. Weisstein, http://mathworld.wolfram.com/RiemannHypothesis.html/ (2008).
[18] http://en.wikipedia.org/wiki/Hilbert's_problems/ (2008).
[19] http://www.claymath.org/millennium/ (2000).

# Frequently Asked Questions

(including some questions that should have been asked, but haven't)

**on "Proof of Riemann's zeta-hypothesis" (arXiv:0809.5120)**
**by Arne Bergstrom**

---

## FAQ #1

**"So many people have been trying in vain for so long to find a proof of the Riemann Hypothesis. Your proof seems to involve no new mathematics, only traditional complex analysis, and uses only methods known even to Riemann himself. So what makes you think you have found an approach that everyone has missed until now?"**

ANSWER: Actually some of the mathematics I use is due to Paul Bachmann/Edmund Landau, and is thus a shade more recent than Riemann's time, but otherwise you are right.

The proof is based on three key elements: 1) a somewhat less-studied formulation of the zeta-function, 2) a particular transformation, the potential of which may possibly have been overlooked, and 3) new, powerful tools in the form of algebraic software.

The particular formulation of the zeta-function I use in (1) in my preprint may perhaps be somewhat less studied than many of the other formulations of the zeta-function – it is not even mentioned in some standard works on the zeta-function such as H M Edwards' book from 1974.

The specific variable transformation (see Sect 2 in my preprint), which I then use to transform this formulation of the zeta-function, may possibly not have been studied with sufficient interest before - if at all. The reason for this might be that at first sight it just seems to complicate the problem by introducing lots of new poles. This, however, turns out to be a blessing in disguise, since suddenly much more structure is introduced into the problem, and which may be used to find a route to the proof.

In my work on the proof I have also benefitted greatly from the existence nowadays of algebraic computer software (e g Maple), which was not available to the old masters, and which permits making long, tedious algebraic calculations with very little effort (even though algebraic software normally needs to be held firmly by the hand so that it does not get lost), and it also permits checking the algebraic derivations numerically with any given high accuracy.

---

## FAQ #2

**"One conceivable road to a proof of RH is to show that the real part σ of all nontrivial zeros of the zeta-function must lie within a narrow strip 1/2 - δ < σ < 1/2 + δ around 1/2, and then prove that δ = 0. A proof of RH must somehow single out the case σ = 1/2 as being special, i e to have some unique, discrete property, and one which σ = 1/2 ± δ with δ ≠ 0 does not have, no matter how small δ is.**

**So what is this unique, discrete property in your proof that makes the case σ = 1/2 so radically different from a value of σ just an infinitesimal bit away?"**

ANSWER: The unique, discrete property you are looking for is the limit when $N$ tends to infinity of the factor $|N^{2s-1}| = N^{2\sigma-1}$ in the numerator in (20) in my preprint. As summarised in my Remark B.2 in Sect B in my preprint, this expression has exactly this property of singling out the case σ = 1/2 as being special, which you are looking for. For σ < 1/2 this factor vanishes when $N \to \infty$, whereas for σ > 1/2 it tends to infinity with $N$. When $N \to \infty$, it has a finite, nonvanishing value ($N^0 = 1$) only in the unique special case σ = 1/2; for any other value of σ, even if only an infinitesimal bit away, it is either zero or infinity.

---

## FAQ #2a

**"Can you summarize in a few words what is the essence of your proof?"**

ANSWER: Yes, I derive the two approximations $\zeta_N(s)'$ and $\zeta_N(s)''$ in (13) and (14) in the preprint of the zeta-function $\zeta(s)$. I then study the quotient $Q_1 = |\zeta_N(s)'/\zeta_N(s)''|$, and show that $Q_1 = \zeta(s)/\zeta(s) = 1$ in limit $N \to \infty$. Independently, I also derive an expression for the quotient $Q_2 = |(\zeta_N(s)' - \zeta(s))/(\zeta_N(s)'' - \zeta(s))|$. At zeros of the zeta-function $\zeta(s) = 0$, I then have $Q_1 = Q_2$. The expression $Q_2$ contains the factor $|N^{2s-1}|$, cf FAQ #2 above, which then requires that σ = ½ in order for $Q_2$ to be equal to $Q_1 = 1$ in the limit $N \to \infty$, thus proving RH.

---

## FAQ #3

**"I think your preprint would be more readable if it was structured better in the customary Theorem&Proof style. This is how readers of mathematical papers expect the material to be presented."**

ANSWER: I agree with you in principle. However, the present proof is rather intricate and I'm afraid that trying to impose a certain form on it would just make it longer without necessarily making it more readable – quite possibly instead having the opposite effect.

I would also like to stress that the detailed analytical calculations in the preprint are absolutely crucial for the proof. Some readers have given me various suggestions how to restructure the presentation. However, these suggestions have invariably only meant that the text would have become longer - and still left the question with the crucial analytical calculations unresolved.

So I think the proof might perhaps be best presented as it is, i e by a rather compact preprint, supplemented with a frequently updated online collection of FAQ, where particular points that readers have found need to be explained in the preprint can be further elaborated and explained in considerable detail, much more so than would have been possible in a journal article.

---

## FAQ #4a (page 3)

**"Does not (4) on page 3 follow immediately from (3) without any requirement on *u* being necessary?"**

ANSWER: No, in order for the step from (3) to (4) to be correct, the condition $-\pi < \mathrm{Im}(u) \leq \pi$ needs to be satisfied (see, e g, Abramotitz and Stegun, eq 4.2.19). Since $w$ is real, it is natural to set also $u$ real, which thus then automatically satisfies the above condition on $u$ – but any $u$ would thus not do so.

---

## FAQ #4b (page 3)

**"In Sect 3 in your preprint, you state that the equation**

$$e^{(e^u)} + 1 = 0$$

**can be verified to have the following solutions (*m* and *n* are integers, $n \geq 1$)**

$$u = \ln(\pi\,(2\,n - 1)) + i\,\pi\left(\frac{1}{2} + m\right) \tag{1}$$

**That is incorrect. The correct solution is given by Mathematica as**

$$u = \ln(i\,\pi\,(2\,n + 1)) + 2\,i\,\pi\,m \tag{2}$$

**where m and n are integers.**

**To try this yourself, enter the following line into Wolfram Alpha:**

**Reduce[Exp[Exp[u]] + 1 == 0, u] ."**

ANSWER: You state categorically that my expression (1) is incorrect since it does not agree with what you get if you ask Mathematica to give you the solution, in which case you get (2) instead

Now, getting something like (2) from Mathematica that looks different from my (1) may however not necessarily imply that my (1) must be wrong. Maybe we should allow for the possibility that they are just two equivalent expressions? For, if you express your (2) in real and imaginary parts, you get

$$u = \ln(\pi\,(2\,n + 1)) + i\,\pi\left(\frac{1}{2} + 2\,m\right) \tag{2a}$$

As you can see, this rewriting of your (2) as (2a) now makes it look rather like my (1), except for how $n$ and $m$ might be defined in the two cases. As defined in the preprint, my $n$ is an integer $n \geq 1$, and my $m$ is a positive or negative integer or zero.

---

## FAQ #4c (page 3)

**"The expression in (2a), derived from Mathematica for $u$ in your previous FAQ,**

$$u = \ln(\pi\,(2\,n + 1)) + i\,\pi\left(\frac{1}{2} + 2\,m\right) \tag{2a}$$

**still does not agree with the expression you give in your preprint,**

$$u = \ln(\pi\,(2\,n - 1)) + i\,\pi\left(\frac{1}{2} + m\right) \tag{1}$$

**So there is thus still a serious problem here."**

<u>ANSWER:</u> My (1) is correct, but this turns out to be a rather intricate question. We want to find the solutions to the equation

$$e^{(e^{u})} + 1 = 0 \,.$$

Rewrite this as

$$e^{(e^{u})} = -1 = e^{(-i\,\pi + 2\,\pi\,i\,n)}, \tag{3}$$

where $n$ is a positive or negative integer or zero.

Form logarithms for both sides in (3)

$$e^{u} = -\,i\,\pi + 2\,\pi\,i\,n = (2n - 1)\,\pi\,i\,. \tag{4}$$

After rewriting $i$, and with $m$ being a positive or negative integer or zero, we get

$$e^{u} = (2n - 1)\,\pi\,e^{(i\,\pi/2 + 2\,\pi\,i\,m)}. \tag{5}$$

Since $n$ is a positive or negative integer or zero, there are thus the following two cases for (5):

**Case A:** $2\,n - 1 > 0$, *i e* $n \geq 1$, in which case we have (with $n = 1, 2, 3 \ldots$),

$$u = \ln(\pi\,(2\,n - 1)) + i\,\pi\left(\frac{1}{2} + 2\,m\right) \tag{6A}$$

**Case B:** $2\,n - 1 < 0$, in which case we can set $2\,n - 1 > 0$ if at the same time we change sign of the exponential factor in (5), *i e* if we add a term $i\,\pi$ to the exponent (since $e^{i\pi} = -1$). The corresponding expression to (6A) thus becomes (where now $n = 1, 2, 3 \ldots$, and $m$ is a positive or negative integer or zero),

$$u = \ln(\pi\,(2\,n - 1)) + i\,\pi\left(\frac{1}{2} + 2\,m + 1\right) \tag{6B}$$

The zeroes in (6A) are thus supplemented by zeroes according to (6B). The zeroes can thus be summarised as in (1) above, and as is also stated in the preprint,

$$u = \ln(\pi\,(2\,n - 1)) + i\,\pi\left(\frac{1}{2} + m\right)$$

and where $n$ is a positive integer $n \geq 1$, and where $m$ is a positive or negative integer or zero.

It should be pointed out that the requirement $n \geq 1$ is what gives the term $m$ instead of $2\,m$ in the expression for $u$ above. If $n$ is instead allowed to take on all integer values (positive or negative or zero), then a negative argument ($i$ $e$ for $n \leq 0$) in the logarithm will give a contribution ln(-1) = $i\,\pi$, corresponding to changing $2\,m$ to $2\,m + 1$ in (2a). So in this respect the expression in (2a) given by Mathematica is thus formally correct, albeit in an easily misleading and not very transparent way.

---

## FAQ #4d (page 3)

**"The residue you calculate on page 3 in your preprint,**

$$\mathrm{Res}(n, m) = i\,(-1)^{m}\,(2\,n - 1)^{(s-1)}\,\pi^{(s-1)}\,\mathrm{e}^{(i\,(1/2+m)\,s\,\pi)}$$

**is wrong. There is a minus sign missing. The correct expression should be**

$$\mathrm{Res}(n, m) = -i\,(-1)^{m}\,(2n - 1)^{(s-1)}\ldots$$

**This sign affects the rest of paper. For example, your formula (11)**

$$\zeta(s) = \frac{-N^{(1-s)}\,\mathrm{E}_N(1-s) + 2^s\left(\sum_{n=1}^{N}(2\,n-1)^{(-s)}\right)(1-s)}{(-1+2^s)\,(1-s)} \tag{11}$$

**is wrong and should correctly be**

$$\zeta(s) = -\frac{N^{(1-s)}\,\mathrm{E}_N(1-s) + 2^s\left(\sum_{n=1}^{N}(2\,n-1)^{(-s)}\right)(1-s)}{(-1+2^s)\,(1-s)} \tag{a}$$

ANSWER: I welcome all questions and objections, but sometimes it would be helpful if the person asking the question could do some simple tests himself before asking. It may be very difficult to find a mathematical error in some long derivation, even if it is done on a computer (or maybe especially in that case, as FAQ #4b and #4c above show). But one can always do some numerical checks. If such a numerical check comes out right, it does of course not prove that something is right. But what a simple numerical check can do, however, is to give a clear indication when something may be basically wrong in a derivation.

Despite the questioner's rather categorical statement above that my expressions are wrong, I have of course continuously made sure that I have done all numerical tests on my derivations that I have been able to invent – even sometimes with as much as 500 digital digits numerical accuracy, as described in FAQ #5a below.

Thus some simple numerical checks (given here with 8 decimals accuracy, *cf* FAQ #5) for, *e g*, $s = 0.01+100\ i$ and $N = 10000$, give in this case

for the questioner's expression (a):

$$\zeta(.01 + 100.\ i) \stackrel{?}{=} 733.90352864 - 531.52064969\ i$$

and for my expression (11):

$$\zeta(.01 + 100.\ i) \stackrel{?}{=} 6.38166671 + .17431634\ i$$

to be compared to the exact value:

$$\zeta(.01 + 100.\ i) = 6.38166671 + .17431634\ i$$

Although of course only a numerical example, this nevertheless strongly suggests that the questioner's expression (a) is not correct, and thus consequently that his expression for the residues is also wrong (which he has also kindly conceded).

I thus maintain that my expression for the residues on line 18, page 3 in my manuscript is correct. I should mention that all derivations in my manuscript have been checked and crosschecked whenever at all possible – and also by many independent mathematicians, as I mention in the Acknowledgment on page 21 of my preprint.

---

## FAQ #4e (page 3)

**"I cannot immediately see how to get to (6) on page 3 in your preprint from the expression for the residues immediately before. Can you write out the intermediate steps for me?"**

ANSWER: Writing for simplicity

$$\phi(s) = i\ (2\ n - 1)^{(s-1)}\ \pi^{(s-1)}$$

the residues as given in Sect 3 in my preprint can be written

$$\mathrm{Res}(n, m) = (-1)^m\ \mathrm{e}^{(i\ (1/2 + m)\ s\ \pi)}\ \phi(s)$$

The sum for $m = 0$ and $m = 1$ then becomes

$$\mathrm{Res}(n, 0) + \mathrm{Res}(n, 1) = \mathrm{e}^{(1/2\ i\ s\ \pi)}\ \phi(s) - \mathrm{e}^{(3/2\ i\ s\ \pi)}\ \phi(s)$$

which can be rewritten as

$$\mathrm{Res}(n, 0) + \mathrm{Res}(n, 1) = \mathrm{e}^{(i\, s\, \pi)} \left(\mathrm{e}^{(-1/2\, i\, s\, \pi)} - \mathrm{e}^{(1/2\, i\, s\, \pi)}\right) \phi(s)$$

Converting the expression within the parenthesis to a sine then gives

$$\mathrm{Res}(n, 0) + \mathrm{Res}(n, 1) = -2\, i\, \mathrm{e}^{(i\, s\, \pi)} \sin\left(\frac{1}{2}\, s\, \pi\right) \phi(s)$$

Reinserting the expression for F(s) above, we get

$$\mathrm{Res}(n, 0) + \mathrm{Res}(n, 1) = 2\, \mathrm{e}^{(i\, s\, \pi)} \sin\left(\frac{1}{2}\, s\, \pi\right) (2\, n - 1)^{(s-1)}\, \pi^{(s-1)}$$

which after summation over $n$ from 1 to $N$ gives my equation (6) in the preprint.

---

## FAQ #5 (pages 4-5)

**"You derive (9) and (11) in your preprint to be two alternative formulations of the zeta-function in the critical strip. But consider (11) on page 5:**

$$\zeta(s) = \frac{-N^{(1-s)}\, \mathrm{E}_N(1-s) + 2^s \left(\sum_{n=1}^{N} (2\, n - 1)^{(-s)}\right) (1 - s)}{(-1 + 2^s)\, (1 - s)} \qquad (11)$$

**In the critical strip the first term in the numerator on the right-hand side is a positive power of *N*, which clearly tends to infinity with *N*. The second term is divergent for $\sigma < 1$.**

**I find it surprising that these two divergent terms can combine to something finite for all values of *s*, and even to the zeta-function. How can one somehow make this plausible theoretically?**

**For, e g, $N = 10^7$ and $s = 0.01 + 100\, i$, the first of those terms is of the order of $10^7$. The second term thus needs to calculate to almost the same $10^7$ in order to combine to the zeta-function, which is of the order of 1."**

<u>ANSWER:</u> The power of Cauchy's theorem is the reason why they do combine to the zeta-function as I derive in Sections 3 and 4 of my preprint, and which is also commented in Remarks 5.2 and 5.3 there. If we insert (12b) and (A9) into (11) in my preprint (and put the remainder with its proper power at the end), we get

$$\zeta(s) = \frac{-N^{(1-s)} \left(1 - \frac{1}{24} \frac{s\,(s-1)}{N^2}\right) + 2^s \left(\sum_{n=1}^{N} (2\, n - 1)^{(-s)}\right) (1 - s)}{(-1 + 2^s)\, (1 - s)} + \mathrm{O}\left(\frac{\ln(N)}{N^{(2+\sigma)}}\right)$$

Calculating an example as the one you suggest is illustrative. With $s = 0.01 + 100\, i$ and $N = 10^7$, the terms on the right-hand side above become (when calculated with 50 digits, rounding off the result to 20 decimals)

412229.02030205334476757607 - 104245.48488699243492823898 $i$ -
412222.63863533509177677516 + 104245.65920333443557465984 $i$ + O($10^{-14}$ )

The large numbers on the right-hand side evaluate to

6.38166671825299080091 + 0.17431634200064642086 $i$ + O($10^{-14}$ )

which should be compared to the corresponding exact value of the zeta-function on the left-hand side which is

6.38166671825299080590 + 0.17431634200064641950 $i$

i e the difference between the calculated value and the exact zeta-function is

$- 0.499\ 10^{-17} + 0.136\ 10^{-17}\ i$

which thus falls well within the accuracy of the order of $10^{-14}$ as defined by the remainder.

An example is of course only an example. But it clearly illustrates the strength of Cauchy's theorem when it requires the two divergent terms in (9) to match.

Inserting (12a) and (A9) similarly gives for the alternative expression for $\zeta(s)$ in (9) in my preprint

$$\zeta(s) = \frac{\pi^s \left( -2^{(s-1)} N^s \left( 1 - \frac{1}{24} \frac{s\,(s-1)}{N^2} \right) + \left( \sum_{n=1}^{N} (2\,n-1)^{(s-1)} \right) s \right)}{\cos\left(\frac{1}{2} s\,\pi\right) (-1 + 2^{(1-s)})\, \Gamma(s+1)} + O\left( \frac{\ln(N)}{N^{(3-\sigma)}} \right)$$

which with $s = 0.01 + 100\,i$ and $N = 10^7$ as above evaluates to

6.38166671825299080590 + 0.17431634200064641950 $i$ + O($10^{-21}$)

with a difference between the calculated value and the exact zeta-function of

$- 0.311\ 10^{-24} + 0.442\ 10^{-24}\ i$

and which thus again falls well within the accuracy of the order of $10^{-21}$ as defined by the remainder.

Omitting the remainders in the two expressions calculated above from (9) and (11) in the preprint, these results correspond to the approximations $\zeta_N(s)'$ and $\zeta_N(s)''$ in (13) and (14). According to the discussion on page 7 in the preprint, the quotient $|\zeta_N(s)'/\zeta_N(s)''|$ should become unity in the limit $N \to \infty$. With the values calculated above we get

$$\frac{\zeta_N(s)'}{\zeta_N(s)''} = \frac{6.38166671825299080590 + 0.17431634200064641950\ i}{6.38166671825299080091 + 0.17431634200064642086\ i}$$

which gives

$$|\zeta_N(s)'/\zeta_N(s)''| = 1.00000000000000000078$$

In agreement with Sect 7 in the preprint, the quotient $|\zeta_N(s)'/\zeta_N(s)''|$ is thus equal to unity within the truncation error in the calculations. We also see from the above calculations that the function $\zeta_N(s)'$ in the second calculation lies closer to $\zeta(s)$ than $\zeta_N(s)''$ in the first calculation does, as they should for $\sigma < ½$ according to Sect 8 in the preprint.

---

## FAQ #5a (pages 4-7)

**"It is all very well that your quotient $|\zeta_N(s)'/\zeta_N(s)''|$ seems to become unity in the limit $N \rightarrow \infty$ for some value of Re($s$) far away from ½, but the crux in your proof of RH is that the quotient $|\zeta_N(s)'/\zeta_N(s)''|$ should then be exactly unity also at zeros of the zeta-function. So why don't you repeat the above calculation at a zero of $\zeta_N(s)$? That would be a much tougher test on the validity of your derivation - not that any numerical example for an isolated case with ever so high accuracy proves anything of course."**

<u>ANSWER:</u> OK, as an example I calculate $\zeta_N(s)'$ and $\zeta_N(s)''$ for, say, the 25th zero of $\zeta(s)$ with an accuracy of 500 digits and with $N = 10^6$.

I then have the zero $\zeta(s_0) = 0$ at

$s_0$ = ½ + 88.809111207634465423682348079509378395444893409818675042199871
618814013559182198439520793279503933064153393514217920973698829552912796
435947430022616561789270621547005213034296606152586194041769538634530945
503364017906804361782732047293903104050652975462272566220454237002694748
322991711060120807226592762152718464656078715516747596277156935025449524
613402429805860214583563456820971738674177274265862494749169298610068752
635619844014549917115019165805602013934741385882124229542710375363168474
405876 *i*

Inserting (12a) and (A9) into (13) I get for $\zeta_N(s)'$ (which thus differs from 0 by the remainder)

$\zeta_N(s_0)'$ = - 0.5598210206848707731410994520741900667412976330788770765008036
300878349353174722551675241728810292446029785798674016633877871251749794
231093091210550197963151743912917094764849214918555239296048134152023473
557742481431882049753947643886769356133228450120886558883649195054476910
177526811281702345355209442148469017097976229715257797960900805875427126
545882612610694751695482473836763960478808056453953943858959422670831170
971226399452268820138742722202517264372269377967458059047445516856227970
9014170 $10^{-18}$
- 0.1429224283891745919434669738274127140045463690562680648278519
807015162252011180817196455403413985104568964029486877196116456363012784
273835631759510834040175739590184135951062546856519713996693021912468731
269347221890898741412491277289482532175044643417845045347197825098382576
148835382066365704537256975393093196783297046458637066016115679034796390
148399941779092971037971155727295584959289475506477211274603703461575780
262741766942825067961262458314703229340943081797695739650870185984615194
6020033 $10^{-18}$ *i*

Similarly inserting (12b) and (A9) into (14) I get for $\zeta_N(s)''$

$\zeta_N(s_0)''$ = - 0.118885488635737982032407088316718813593396998979140135378408
927185266850261594283017500975118423960470984586658853755884524257332688
428939642086728176189836187209169259991285501919906819480784526568559524
762337259290275322938313185822513706802048515784740461783016585802794345
277282236358529143718103220901167926334438503057410112729869431613713758
922001597951837819108646115918387331602871718921342250938827734512134516
313885287375541308106450076685153579250518589118607654925671472489883690
67817269 $10^{-18}$

\+ 0.565413686011535331422018295350295874961405420529835408493291
416287419825108067298691784936160554126957273806815046626804981164871179
068033043282715281559247568272057910968694790590494707543075853678476697
054473426531443922937723512436532054572368980461019259530103882954072754
665560321222593445767090972067061321584979710308249297971463416171316180
674215616730008236695692721622467840015605346317614991018622869125056488
583584423765936974414805142350116184222333215357488553591958976253618453
43344048 $10^{-18}$ ***i***

and hence I can finally calculate the quotient

$$|\zeta_N(s)'/\zeta_N(s)''| = 1 + 0.74\ 10^{-480}$$

which thus agrees with unity within the accuracy used in the calculations. But, as you say, a numerical example of course does not prove anything – there is a long way from $N = 10^6$ to infinity and from the 25th zero to every zero. However, what proves that $|\zeta_N(s)'/\zeta_N(s)''| = 1$ for all $s$ when $N \to \infty$ is the fact, as discussed on page 6 in the preprint and further in FAQ #7 below, that $|\zeta_N(s)'|$ and $|\zeta_N(s)''|$ differ from $|\zeta(s)|$ only by the remainders, and these remainders vanish when $N \to \infty$, so that in this limit we have $|\zeta_\infty(s)'| = |\zeta_\infty(s)''| = |\zeta(s)|$ for all $s$, i e $|\zeta_\infty(s)'/\zeta_\infty(s)''| = |\zeta(s)/\zeta(s)| = 1$.

---

## FAQ #6 (pages 4-6)

**"I think there could be a problem in connection with your discussion on top of page 6 in the preprint and the two expressions (9) and (11), which you have derived to give two variants of the (exact) zeta-function. If so, then one could solve the two error functions $E_N(s)$ and $E_N(1-s)$ from (9) and (11). These two expressions for the error functions should then be equal after setting $s \to 1-s$ in one of them. But judging from how different (9) and (11) are, this would seem unlikely, and there could hence be an inconsistency here.**

ANSWER: Expressions (9) and (11) in the preprint read

$$\zeta(s)=\frac{\pi^s\left(-2^{(s-1)}N^s\,\mathrm{E}_N(s)+\left(\sum_{n=1}^{N}(2n-1)^{(s-1)}\right)s\right)}{\cos\left(\frac{1}{2}s\pi\right)(-1+2^{(1-s)})\,\Gamma(s+1)} \tag{9}$$

$$\zeta(s)=\frac{-N^{(1-s)}\,\mathrm{E}_N(1-s)+2^s\left(\sum_{n=1}^{N}(2n-1)^{(-s)}\right)(1-s)}{(-1+2^s)(1-s)} \tag{11}$$

Solving for $\mathrm{E}_N(s)$ and $\mathrm{E}_N(1 - s)$, respectively, we get

$$\mathrm{E}_N(s)=2\,\frac{\zeta(s)\cos\left(\frac{1}{2}s\pi\right)\Gamma(s+1)-\zeta(s)\cos\left(\frac{1}{2}s\pi\right)\Gamma(s+1)\,2^{(1-s)}+\pi^s\left(\sum_{n=1}^{N}(2n-1)^{(s-1)}\right)s}{\pi^s\,2^s\,N^s}$$

$$\mathrm{E}_N(1-s)\ =$$

$$\frac{1}{2}\,\frac{-2\,\zeta(s)\,s+2\,\zeta(s)+2\,\zeta(s)\,2^s s-2\,\zeta(s)\,2^s+2^{(s+1)}\left(\sum_{n=1}^{N}(2n-1)^{(-s)}\right)-2^{(s+1)}\left(\sum_{n=1}^{N}(2n-1)^{(-s)}\right)s}{N^{(1-s)}}$$

Setting $s$ -> 1 - $s$ in the second equation above as you say, we get

$$\mathrm{E}_N(s)=\frac{1}{2}\Bigg(-2\,\zeta(1-s)(1-s)+2\,\zeta(1-s)+2\,\zeta(1-s)\,2^{(1-s)}(1-s)-2\,\zeta(1-s)\,2^{(1-s)}$$
$$+2^{(2-s)}\left(\sum_{n=1}^{N}(2n-1)^{(s-1)}\right)-2^{(2-s)}\left(\sum_{n=1}^{N}(2n-1)^{(s-1)}\right)(1-s)\Bigg)\Big/N^s$$

From (9) and (11) we have thus obtained two equations above for $\mathrm{E}_N(s)$. Forming the quotient of their left- and right-hand sides of these equations we get

$$1=4\left(\zeta(s)\cos\left(\frac{1}{2}s\pi\right)\Gamma(s+1)-\zeta(s)\cos\left(\frac{1}{2}s\pi\right)\Gamma(s+1)\,2^{(1-s)}+\pi^s\left(\sum_{n=1}^{N}(2n-1)^{(s-1)}\right)s\right)\Big/\Bigg($$
$$\pi^s\,2^s\Bigg(-2\,\zeta(1-s)(1-s)+2\,\zeta(1-s)+2\,\zeta(1-s)\,2^{(1-s)}(1-s)-2\,\zeta(1-s)\,2^{(1-s)}$$
$$+2^{(2-s)}\left(\sum_{n=1}^{N}(2n-1)^{(s-1)}\right)-2^{(2-s)}\left(\sum_{n=1}^{N}(2n-1)^{(s-1)}\right)(1-s)\Bigg)\Bigg)$$

After using the functional equation (2) in the preprint

$$\zeta(s)=2^s\,\pi^{(s-1)}\sin\left(\frac{1}{2}s\pi\right)\Gamma(1-s)\,\zeta(1-s)$$

the quotient above becomes

$$1=4\left(2^s\pi^{(s-1)}\sin\left(\frac{1}{2}s\pi\right)\Gamma(1-s)\zeta(1-s)\cos\left(\frac{1}{2}s\pi\right)\Gamma(s+1)\right.$$

$$\left.-2^s\pi^{(s-1)}\sin\left(\frac{1}{2}s\pi\right)\Gamma(1-s)\zeta(1-s)\cos\left(\frac{1}{2}s\pi\right)\Gamma(s+1)\,2^{(1-s)}+\pi^s\left(\sum_{n=1}^{N}(2n-1)^{(s-1)}\right)s\right)$$

$$\Big/\left(\pi^s 2^s\left(-2\zeta(1-s)(1-s)+2\zeta(1-s)+2\zeta(1-s)\,2^{(1-s)}(1-s)-2\zeta(1-s)\,2^{(1-s)}\right.\right.$$

$$\left.\left.+2^{(2-s)}\left(\sum_{n=1}^{N}(2n-1)^{(s-1)}\right)-2^{(2-s)}\left(\sum_{n=1}^{N}(2n-1)^{(s-1)}\right)(1-s)\right)\right)$$

which can be simplified to

$$1=-2\frac{2^s\sin\left(\frac{1}{2}s\pi\right)\zeta(1-s)\cos\left(\frac{1}{2}s\pi\right)-2\sin\left(\frac{1}{2}s\pi\right)\zeta(1-s)\cos\left(\frac{1}{2}s\pi\right)+\left(\sum_{n=1}^{N}(2n-1)^{(s-1)}\right)\sin(s\pi)}{\sin(\pi(s+1))\left(\zeta(1-s)2^s-2\zeta(1-s)+2\left(\sum_{n=1}^{N}(2n-1)^{(s-1)}\right)\right)}$$

or equivalently

$$\sin(\pi(s+1))\left(\zeta(1-s)2^s-2\zeta(1-s)+2\left(\sum_{n=1}^{N}(2n-1)^{(s-1)}\right)\right)+2\,2^s\sin\left(\frac{1}{2}s\pi\right)\zeta(1-s)\cos\left(\frac{1}{2}s\pi\right)$$

$$-4\sin\left(\frac{1}{2}s\pi\right)\zeta(1-s)\cos\left(\frac{1}{2}s\pi\right)+2\left(\sum_{n=1}^{N}(2n-1)^{(s-1)}\right)\sin(s\pi)=0$$

which simplifies to

$$-\sin(s\pi)\,\zeta(1-s)\,2^s+2\sin(s\pi)\,\zeta(1-s)+2\,2^s\sin\left(\frac{1}{2}s\pi\right)\zeta(1-s)\cos\left(\frac{1}{2}s\pi\right)$$

$$-4\sin\left(\frac{1}{2}s\pi\right)\zeta(1-s)\cos\left(\frac{1}{2}s\pi\right)=0$$

or

$$-\zeta(1-s)(-2+2^s)\left(\sin(s\pi)-2\cos\left(\frac{1}{2}s\pi\right)\sin\left(\frac{1}{2}s\pi\right)\right)=0$$

which evaluates to

$$0=0$$

Thus the two expressions for the error functions are equal, as they should be, and there is no inconsistency in the equations (9) and (11) for the zeta-function in the preprint.

## FAQ #6a (page 6)

**"The text on page 6 in your preprint seems somewhat vague. You seem to use a single notation, $\zeta_N(s)$, to denote many different functions."**

ANSWER: No, actually I don't. I use $\zeta_N(s)$ with prime and double-prime, respectively, to denote the approximations in (13) and (14)

$$\zeta_N(s)' = \frac{\pi^s \left( -2^{(s-1)} N^s \, \Xi_N(s) + \left( \sum_{n=1}^{N} (2n-1)^{(s-1)} \right) s \right)}{\cos\left(\frac{1}{2} s \pi\right) (-1 + 2^{(1-s)}) \, \Gamma(s+1)} \tag{13}$$

$$\zeta_N(s)'' = \frac{-N^{(1-s)} \, \Xi_N(1-s) + 2^s \left( \sum_{n=1}^{N} (2n-1)^{(-s)} \right) (1-s)}{(-1 + 2^s)(1-s)} \tag{14}$$

of the exact expressions for zeta-function $\zeta(s)$ in (9) and (11) in my preprint,

$$\zeta(s) = \frac{\pi^s \left( -2^{(s-1)} N^s \, \mathrm{E}_N(s) + \left( \sum_{n=1}^{N} (2n-1)^{(s-1)} \right) s \right)}{\cos\left(\frac{1}{2} s \pi\right) (-1 + 2^{(1-s)}) \, \Gamma(s+1)} \tag{9}$$

$$\zeta(s) = \frac{-N^{(1-s)} \, \mathrm{E}_N(1-s) + 2^s \left( \sum_{n=1}^{N} (2n-1)^{(-s)} \right) (1-s)}{(-1 + 2^s)(1-s)} \tag{11}$$

That I have the two expressions in (9) and (11) for the exact zeta-function $\zeta(s)$ is because one can always write the exact zeta-function in two different ways by using the functional equation combined with a variable transformation, as I show in Sect 5 in my preprint. As discussed in the preprint, working with this function pair for $\zeta(s)$ is advantageous since in a simple way it automatically incorporates the functional equation with its symmetry properties into the derivation – and these symmetry properties are important for the proof. But, for each $N$, the expressions in (9) and (11) give the exact zeta-function $\zeta(s)$ only for a very special choice of the error functions, denoted by $\mathrm{E}_N(s)$ and $\mathrm{E}_N(1\text{-}s)$, within their respective remainders $O(1/N^{3\cdot})$, as I discuss in the first paragraph in Sect 6 in the preprint.

The trick in the proof is then to consider also error functions around these two very special choices $\mathrm{E}_N(s)$ and $\mathrm{E}_N(1\text{-}s)$ of the error functions that give the exact expressions for $\zeta(s)$ in (9) and (11). Such modified error functions can be obtained by varying the error functions within the remainders around the error functions that give the exact zeta-function. If inserted into (9) and (11), all such modified error functions would of course destroy the accuracy of (9) and (11). Then those relationships would no longer give the exact zeta-function, nor would they obey the functional equation. In particular, I consider the error functions I get if I set the remainders $O(1/N^{3\cdot})$ in (12a) and (12b) equal to zero. I call those error functions $\Xi_N(s)$ and $\Xi_N(1\text{-}s)$, and denote the corresponding approximations of the zeta-function by $\zeta_N(s)'$ and $\zeta_N(s)''$, respectively.

So $\zeta(s)$ thus denotes the exact zeta-function in the two exact forms given in (9) and (11) with their special error functions $E_N(s)$ and $E_N(1\text{-}s)$, whereas $\zeta_N(s)'$ and $\zeta_N(s)''$ in (13) and (14) denote the two approximations of the zeta-function I get when I set the remainders in the error functions equal to zero instead of equal to the particular remainders that give the special error functions $E_N(s)$ and $E_N(1\text{-}s)$ in (9) and (11), corresponding to the exact zeta-function.

Of course we may then later perhaps want to make all possible error functions $E_N(s)$, $E_N(1\text{-}s)$ and $\Xi_N(s)$, $\Xi_N(1\text{-}s)$ equal (and equal to 1) by considering the case $N \to \infty$, in which case we would get $|\zeta_\infty(s)'| = |\zeta(s)|$ and $|\zeta_\infty(s)''| = |\zeta(s)|$, but that should not really cause any confusion.

---

## FAQ #6b (pages 4-6)

**"In (9) and (11) in the preprint you give what seems to be two new expressions for the zeta-function in the critical strip. But what about previously known expressions for the zeta-function? Could you not have used a similar pair of them?"**

ANSWER: I tried that first, but I could not get anywhere; that is why I ended up with the ones I have now derived.

For comparison I below give five different expressions for the zeta-function in the critical strip, which are obtained as follows. Expression (a) is equivalent to (1) in the preprint, from which the expressions (d) and (e) can be derived as in (9) and (11), and which are here given in the limit $N \to \infty$ when the remainders vanish. Expression (b) is a known alternating Dirichlet series for the zeta-function valid for Re($s$) > 0, see ref [11] in the preprint, and expression (c) is the result of an exactly analogous derivation as to how (11) in the preprint was derived from (9).

$$\zeta(s) = \frac{\displaystyle\int_0^\infty \frac{w^{(s-1)}}{e^w + 1}\, dw}{(1 - 2^{(1-s)})\,\Gamma(s)} \qquad (a)$$

$$\zeta(s) = \lim_{N \to \infty} \frac{\displaystyle\sum_{n=1}^{N} \frac{(-1)^{(n+1)}}{n^s}}{1 - 2^{(1-s)}} \qquad (b)$$

$$\zeta(s) = \lim_{N\to\infty} \frac{2^s \pi^{(s-1)} \sin\left(\frac{1}{2} s \pi\right) \Gamma(1-s) \left( \sum_{n=1}^{N} \frac{(-1)^{(n+1)}}{n^{(1-s)}} \right)}{1 - 2^s} \qquad (c)$$

$$\zeta(s) = \lim_{N\to\infty} \frac{\pi^s \left( -2^{(s-1)} N^s \left( 1 - \frac{1}{24} \frac{s(s-1)}{N^2} \right) + \left( \sum_{n=1}^{N} (2n-1)^{(s-1)} \right) s \right)}{\cos\left(\frac{1}{2} s \pi\right) (-1 + 2^{(1-s)}) \Gamma(s+1)} \qquad (d)$$

$$\zeta(s) = \lim_{N\to\infty} \frac{-N^{(1-s)} \left( 1 - \frac{1}{24} \frac{s(s-1)}{N^2} \right) + 2^s \left( \sum_{n=1}^{N} (2n-1)^{(-s)} \right) (1-s)}{(-1 + 2^s)(1-s)} \qquad (e)$$

---

## FAQ #6c (page 7)

**"Your assertion**

$$\lim_{N\to\infty} \left| \frac{\zeta_N(s_0)'}{\zeta_N(s_0)''} \right| = 1$$

**in (15) on page 7 in the preprint and further in Theorem C.I on page 19 is not justified, because from the formula before it on page 7, one can conclude that $\zeta_N(s)'$ and $\zeta_N(s)''$ can be written as follows**

$$\zeta_N(s)' = \zeta(s) + A_N(s) / N^{(3-\sigma)}$$

$$\zeta_N(s)'' = \zeta(s) + B_N(s) / N^{(2+\sigma)}$$

**where, for fixed $s$, $A_N(s)$ and $B_N(s)$ are bounded as functions of $N$. Thus**

$$\frac{\zeta_N(s)'}{\zeta_N(s)''} = 1 + \frac{N^{(2\sigma-1)} A_N(s) - B_N(s)}{N^{(2+\sigma)} \zeta(s) + B_N(s)}$$

**For fixed $s$ with $\zeta(s) \neq 0$ this indeed tends to 1 as $N \to \infty$. But I see no reason to believe that this convergence is uniform in $s$. That would require bounding the remainder term independently of $s$. Without uniform convergence in $s$ one cannot conclude that there is convergence."**

ANSWER: You consider the function

$$a_{n\,m} = \left| \frac{\zeta_N(s)'}{\zeta_N(s)''} \right| = \left| \frac{\zeta(s) + A_N(s)/N^{(3-\sigma)}}{\zeta(s) + B_N(s)/N^{(2+\sigma)}} \right|,$$

and rightly conclude that this function is not necessarily uniformly convergent in $s$ (cf figure below).

But that is not the point. Theorem C on uniform convergence on page 19 in the preprint does not say that $a_{nm}$ should converge uniformly. But what it says should converge uniformly is the limit $lim\; a_{nm}$ for $n \to \infty$ as given in (α) in my Appendix C.I. And as I show in my Remark C.1 there, this limit is unity for all $s$ in this case, and thus converges trivially.

I thus maintain that my expression (15) on page 7 and (ζ) in Appendix C.I is correct.

As illustration, the figure below (from FAQ #7a) shows the quotient $a_{nm}$ above plotted over the complex plane. This quotient has a pole and a zero close to a zero of the zeta-function, but still between the pole and the zero has the value unity for σ = ½. As $N$ increases, the extensions of the pole and the zero get more and more localised, and drift closer and closer together towards σ = ½, but the straight line $a_{nm} = 1$ for σ = ½ between them remains. This thus gives a numerical illustration of how the quotient of the two identical derivatives in (15) on page 7 in the preprint gives the result unity.

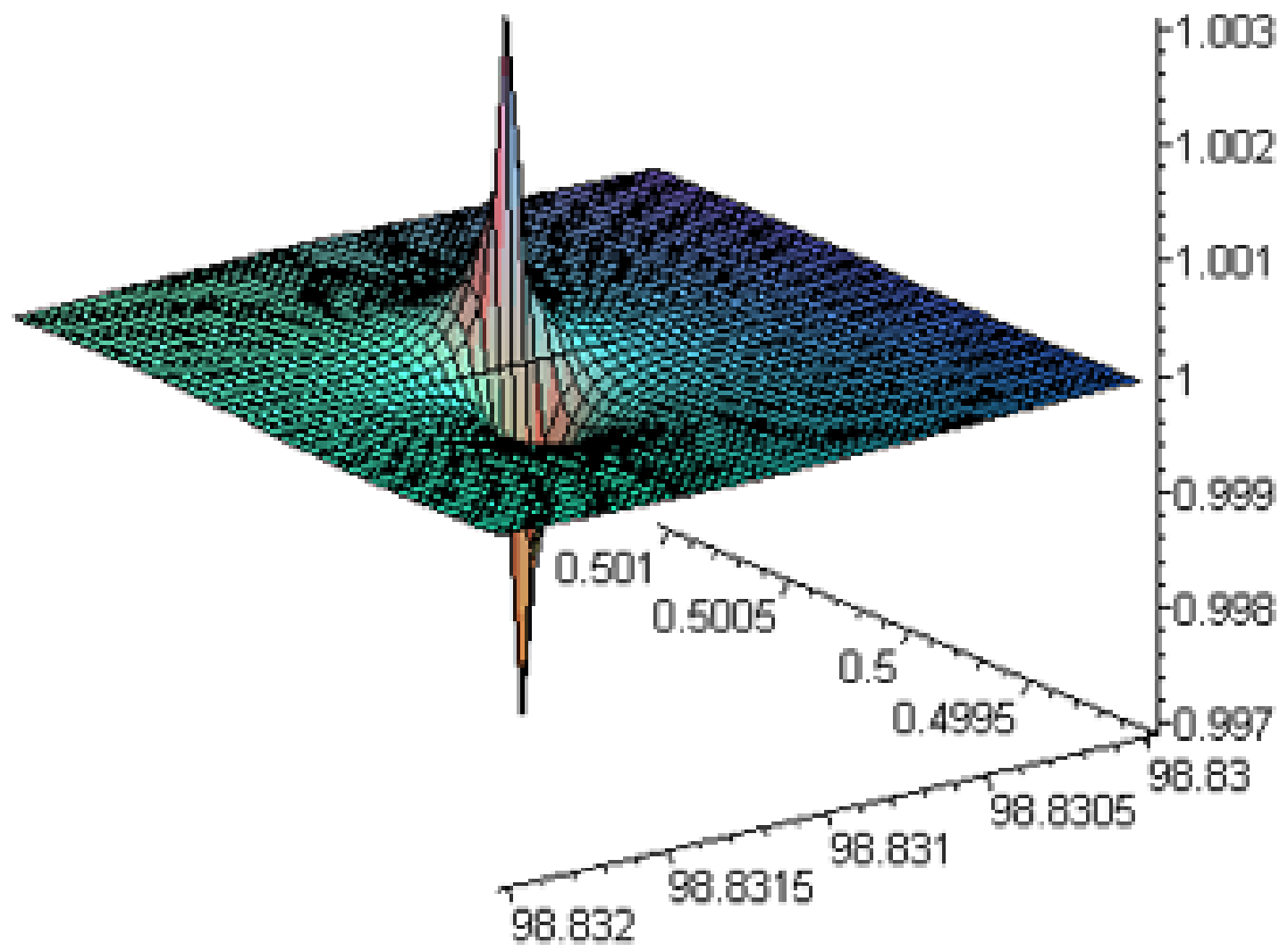

## FAQ #7 (page 8)

**"There may be an error in your argument on page 8 in the preprint. When $s = s_0$ is a zero of the zeta-function, then $\zeta(s_0)$, $\zeta_N(s_0)'$, and $\zeta_N(s_0)''$ are all zero in the limit $N \to \infty$. Thus if you set $\zeta(s) = 0$ in**

$$|\zeta_N(s)' - \zeta(s)|/|\zeta_N(s)'' - \zeta(s)| \qquad (1)$$

**as you do, then the quotient (1) in that case need not at all necessarily become equal to the quotient in (2),**

$$|\zeta_N(s)'|/|\zeta_N(s)''| \qquad (2)$$

**Consider, e g, the following counterexample when we set $\zeta(s) = \phi$, $\zeta_N(s)' = \sin(\phi)$, and $\zeta_N(\phi)'' = \tan(\phi)$,**

$$(\sin(\phi) - \phi)/(\tan(\phi) - \phi) \qquad (1')$$

**Setting $\phi = 0$ in (1') does not make the quotient (1') equal to (2')**

$$\sin(\phi)/\tan(\phi) \qquad (2')$$

**Instead (2') is equal to 1 for $\phi = 0$, whereas (1') then is equal to -½."**

ANSWER: As will be seen below, the problem illustrated by your counterexample is not completely applicable to the case studied in the preprint.

Since

$$\sin(\phi) = \phi - \frac{1}{6}\phi^3 + O(\phi^5) \qquad \tan(\phi) = \phi + \frac{1}{3}\phi^3 + O(\phi^5)$$

$$\sin(\phi) - \phi = -\frac{1}{6}\phi^3 + O(\phi^5) \qquad \tan(\phi) - \phi = \frac{1}{3}\phi^3 + O(\phi^5)$$

we have for your (1') and (2')

$$\frac{\sin(\phi) - \phi}{\tan(\phi) - \phi} = -\frac{1}{2} + O(\phi^2) \qquad (1')$$

$$\frac{\sin(\phi)}{\tan(\phi)} = 1 - \frac{1}{2}\phi^2 + O(\phi^4) \qquad (2')$$

so that the quotient on the left-hand side of (1') for $\phi = 0$ is thus indeed different from the quotient on the left-hand side of (2'), as you say.

In contrast, my quotients (1) and (2) above can be calculated in the limit $N \to \infty$ to be as follows. [The quotients in (1) and (2) need to be taken both in the limit $N \to \infty$ and in the limit $\zeta(s) \to 0$. We here first take the limits in the order $N \to \infty$, $\zeta(s) \to 0$, and then in FAQ #7a below in the order $\zeta(s) \to 0$, $N \to \infty$ as comparison, where both alternatives give the same result (cf Appendix C in the preprint).]

From (19) (and Remark 9.2) in the preprint, we get for (1) above in the limit $N \to \infty$

$$\lim_{N\to\infty}\left|\frac{\zeta_N(s)' - \zeta(s)}{\zeta_N(s)'' - \zeta(s)}\right| = \begin{cases} 0 & \sigma < \frac{1}{2} \\ 1 & \sigma = \frac{1}{2} \\ \infty & \frac{1}{2} < \sigma \end{cases} \qquad (1a)$$

Since $\zeta_N(s)'$ and $\zeta_N(s)''$ differ from $\zeta(s)$ only by the remainders, we have from (2) above in the limit $N \to \infty$, where these remainders vanish,

$$\lim_{N\to\infty}\left|\frac{\zeta_N(s)'}{\zeta_N(s)''}\right| = \left|\frac{\zeta(s)}{\zeta(s)}\right| = 1 \qquad 0 < \sigma < 1 \qquad (2a)$$

We see that for $\sigma \neq ½$ there is a problem between (1a) and (2a) of the type in your counterexample above. On the other hand, in the case $\sigma = ½$ the quotients in (1a) and (2a) are indeed equal. Comparing (1a) and (2a), we see that this corresponds to the case $\zeta(s) = 0$.

[Parenthetically, we note from (1a) (and discussed in Remark B.2 in the preprint) that also cases with $\zeta(s) \neq 0$ can make the right-hand side of (1a) equal to unity (namely for $s = ½ + i\,t$, with such $t$ that $\zeta(s) \neq 0$), but that does not affect the conclusion above.]

From the above we can thus conclude that the special case $\zeta(s) = 0$ can occur if and only if $\sigma = ½$, which thus proves RH in accordance with the discussion in the preprint.

---

## FAQ #7a (page 8)

**"In your proof you first calculate the quotients in the limit $N \to \infty$, and then you study this limit when $s$ tends to a zero of $\zeta(s)$. What happens if you go to the limits in reverse order, i e first calculate the limit in $s$ and then go to the limit in $N$?"**

ANSWER: This question is treated more generally in Appendix C in the preprint. Prior to calculating the quotient in (18) in Appendix B in the preprint, we have for the differences $\zeta_N(s)' - \zeta(s)$ and $\zeta_N(s)'' - \zeta(s)$ relative to the exact zeta-function $\zeta(s)$ the following two expressions

$$\zeta_N(s)' - \zeta(s) = -\frac{7}{11520}\,\frac{N^{(s-4)}\,\pi^s\,4^s}{\cos\left(\frac{1}{2}s\,\pi\right)(-2+2^s)\,\Gamma(s-3)} + O(N^{(s-5)})$$

$$\zeta_N(s)'' - \zeta(s) = -\frac{7}{5760}\,\frac{N^{(-3-s)}\,(s+2)\,(s+1)\,s}{-1+2^s} + O(N^{(-s-4)})$$

Starting in this case by first calculating the limit $\zeta(s) = 0$, the values of $\zeta_N(s)'$ and $\zeta_N(s)''$ in this limit can be obtained by setting $\zeta(s) = 0$ in the expressions above. We then get

$$\zeta_N(s)' = -\frac{7}{11520} \frac{N^{(s-4)} \pi^s 4^s}{\cos\left(\frac{1}{2} s \pi\right) (-2 + 2^s) \Gamma(s-3)} + O(N^{(s-5)})$$

$$\zeta_N(s)'' = -\frac{7}{5760} \frac{N^{(-s-3)} (s+2)(s+1) s}{-1 + 2^s} + O(N^{(-s-4)})$$

With the right-hand sides as above, the quotient $\zeta_N(s)'/\zeta_N(s)''$ thus evaluates to the right-hand side of (18) in the preprint. Studying next this expression for the quotient $\zeta_N(s)'/\zeta_N(s)''$ in the limit $N \rightarrow \infty$, where $|\zeta_\infty(s)'/\zeta_\infty(s)''| = |\zeta(s)/\zeta(s)| = 1$, this requires $\sigma$ on the right-hand side to be $\sigma = ½$ as discussed in Sect 6 in the preprint, and thus again proves RH.

As illustration, the loglog plot below shows $|\zeta_N(s)'|$ and $|\zeta_N(s)''|$ (coinciding, red) as functions of $N$, calculated from (13) and (14) (or from above), and compared to the function $10^3 N^{-3.5}$ (green) for values of $N$ from $10^2$ to $10^7$ and for the 25th zero of $\zeta(s)$ as calculated in FAQ #5a above.

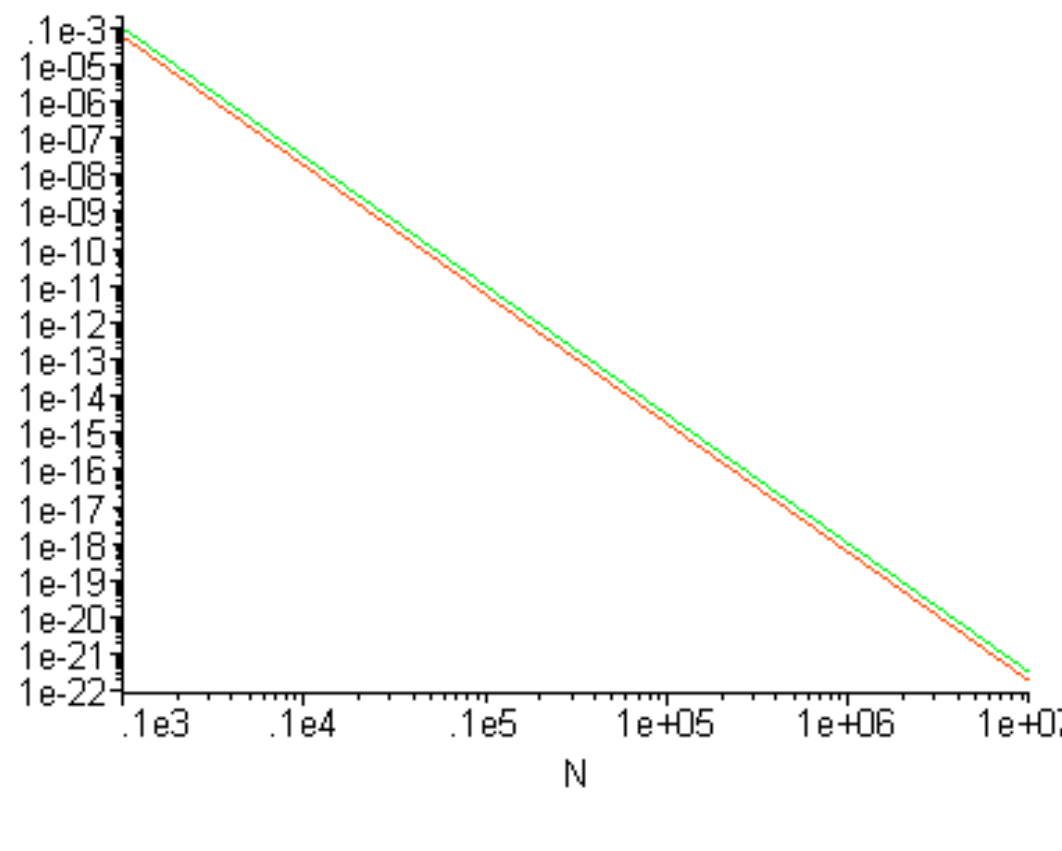

This plot illustrates the fact that the quotient $|\zeta_N(s)'/\zeta_N(s)''|$ is "normally" equal to unity. However, the zeros of the functions $\zeta_N(s)'$ and $\zeta_N(s)''$, which are approximations of $\zeta(s)$, do not completely coincide for finite $N$. Hence, close to a zero of $\zeta(s)$ the quotient $|\zeta_N(s)'/\zeta_N(s)''|$ will contain pairs of adjacent points where the functions $\zeta_N(s)'$ and $\zeta_N(s)''$ are zero independently, and where thus the quotient $|\zeta_N(s)'/\zeta_N(s)''|$ is zero and infinity, respectively, as illustrated by the following plot.

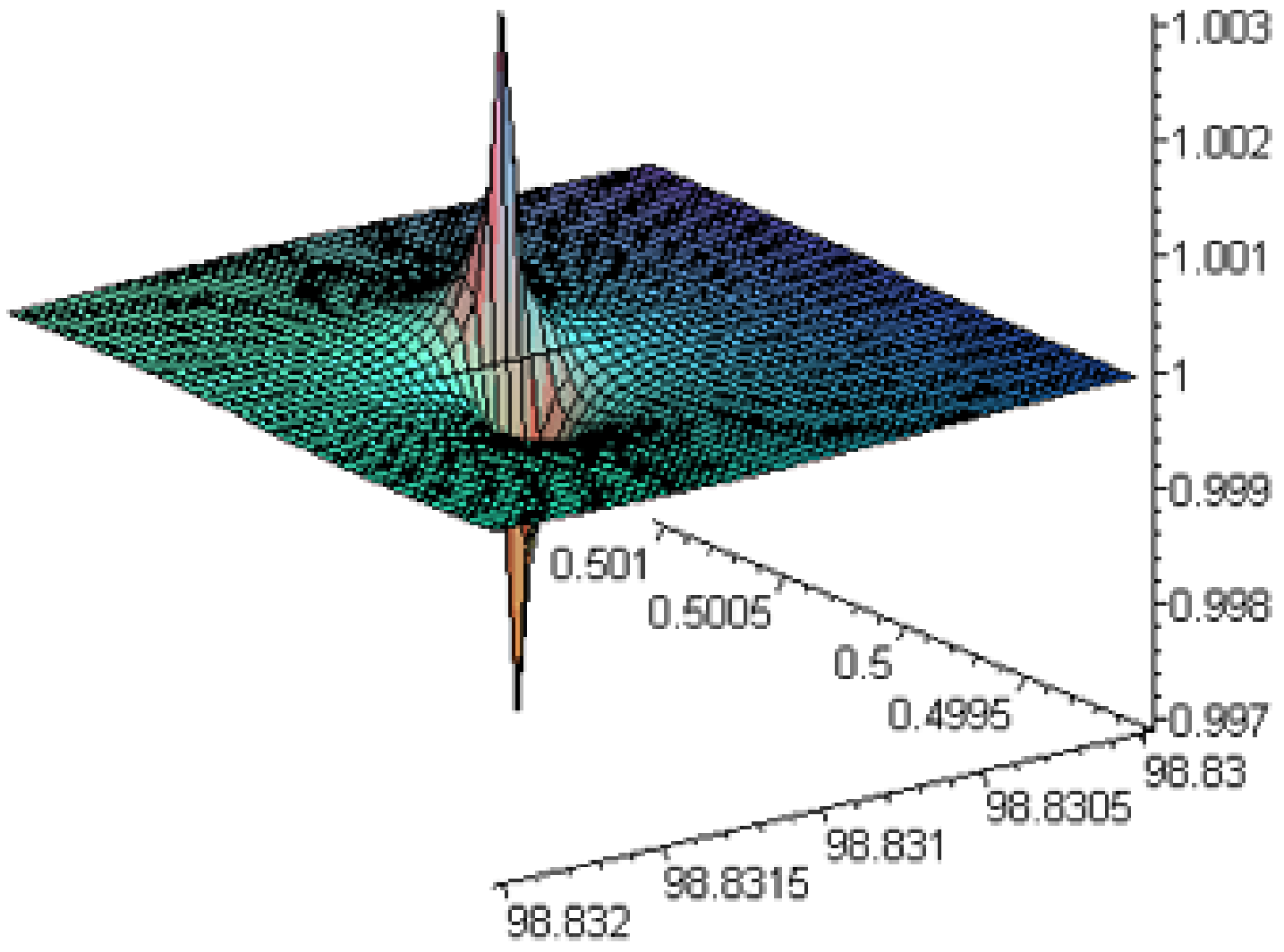

The above plot for the 29th zero of $\zeta(s)$ clearly illustrates the property discussed in remark 9.2 in Sect. 9 in the preprint that the quotient $|\zeta_N(s)'/\zeta_N(s)''|$ - in between the pole and the zero - is indeed unity for all $t$ in $s = \frac{1}{2} + i\,t$.

---

## FAQ #7b (page 8)

**"In your paper you take the limits in the order $N \to \infty$ , $\zeta(s) \to 0$, whereas in FAQ #7a above you take them in the reverse order $\zeta(s) \to 0$, $N \to \infty$. But it seems to me that in both cases you are doing the same analysis of expression (18) in the preprint."**

ANSWER: Yes, superficially it may look like that, because in both cases I end up using an expression of type (18) below - but in different contexts.

$$\frac{\zeta_N(s)' - \zeta(s)}{\zeta_N(s)'' - \zeta(s)} = \frac{1}{2}\,\frac{N^{(2s-1)}\,\pi^s\,(-4^s + 8^s)}{\cos\left(\frac{1}{2}\,s\,\pi\right)(-2+2^s)\,\Gamma(s-3)\,(s+2)\,(s+1)\,s} + O(N^{(2\sigma-2)}) \quad (18)$$

**In my paper,** I start by considering the following quotient in the limit $N \to \infty$.

$$\lim_{N\to\infty}\left|\frac{\zeta_N(s)' - \zeta(s)}{\zeta_N(s)'' - \zeta(s)}\right| \qquad (i)$$

Next I want to set $\zeta(s) = 0$ there and then use (18) to get an equation from which I can solve the values of $s$ which thus give a zero of $\zeta(s)$. Setting $\zeta(s) = 0$ in (*i*) above, I would expect to get the following limit,

$$\lim_{N \to \infty} \left| \frac{\zeta_N(s)'}{\zeta_N(s)''} \right| = 1 \qquad (ii)$$

which I know to have the value unity, since in the limit $N \to \infty$ we have $|\zeta_N(s)'/\zeta_N(s)''| = |\zeta_\infty(s)'/\zeta_\infty(s)''| = |\zeta(s)/\zeta(s)| = 1$ (cf also the numerical example in FAQ #7a above). In this way I would thus get unity on the left-hand side of (18) above and thus have an equation from which I can solve for the values of $s$ which give a zero of $\zeta(s)$.

Unfortunately, setting $\zeta(s) = 0$ in (*i*) does not necessarily give (*ii*) as shown by the counterexample in FAQ #7. Fortunately, however, there is as shown in FAQ #7 one case when this is true, namely when $\mathrm{Re}(s) = ½$. So this thus proves RH anyway.

**In FAQ #7a,** I <u>start</u> instead by setting $\zeta(s) = 0$ in the explicit expressions given there for the differences $\zeta_N(s)' - \zeta(s)$ and $\zeta_N(s)'' - \zeta(s)$. In this way I thus get explicit expressions for $\zeta_N(s)'$ and $\zeta_N(s)''$ in the case $\zeta(s) = 0$. These expressions give the following quotient

$$\frac{\zeta_N(s)'}{\zeta_N(s)''} = \frac{1}{2} \frac{N^{(2s-1)} \pi^s (-4^s + 8^s)}{\cos\left(\frac{1}{2} s \pi\right) (-2 + 2^s) \Gamma(s-3) (s+2) (s+1) s} + O(N^{(2\sigma - 2)})$$

<u>Next</u> letting $N$ tend to infinity, the left-hand side above becomes unity, since as discussed above $|\zeta_N(s)'/\zeta_N(s)''| = |\zeta_\infty(s)'/\zeta_\infty(s)''| = |\zeta(s)/\zeta(s)| = 1$, and this is compatible with the limit $N \to \infty$ of the right-hand side above only if $\mathrm{Re}(s) = ½$, which thus again proves RH.

---

## FAQ #7c (page 8)

**"I have difficulties visualizing the functions $\zeta_N(s)'$ and $\zeta_N(s)''$ relative to the 'normal' zeta-function. Can you give some plots that show how these functions behave around a zero of the zeta-function? Do they have minima or zeros there?"**

<u>ANSWER:</u> Figures 1 through 6 below show $\zeta_N(s)'$ (red), $\zeta_N(s)''$ (turquoise), and $\zeta(s)$ (green) as functions over the complex plane $s = \sigma + i\,t$ around the zero of the zeta-function at $s = ½ + 40.918719\,i$.

In a sufficiently small domain around a zero $s = s_0$ of the zeta-function $\zeta(s)$, the real parts of the functions $\zeta_N(s)'$, $\zeta_N(s)''$, and $\zeta(s)$, can be approximated by parallel planes over the complex plane $s = \sigma + i\,t$, separated from each other by distances corresponding to the real part of the remainders, as illustrated in Figure 1 below. The imaginary parts of the functions $\zeta_N(s)'$, $\zeta_N(s)''$, and $\zeta(s)$ can similarly also be approximated by parallel planes over the complex plane, separated by distances corresponding to the imaginary part of the remainders as illustrated in Figure 2 below.

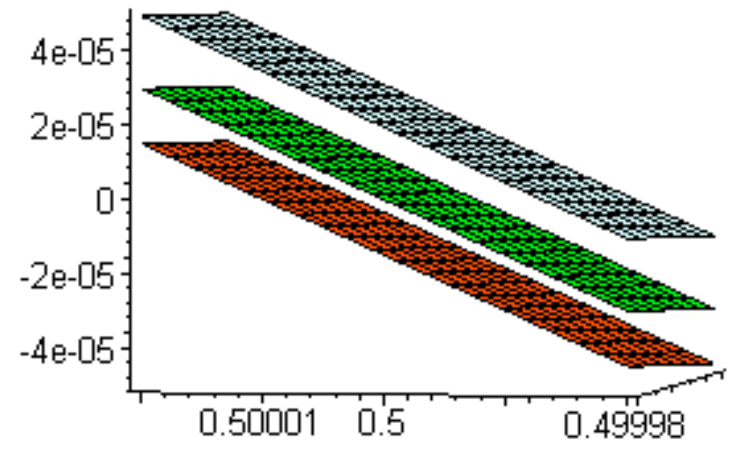

Fig 1. Real parts of $\zeta_N(s)'$, $\zeta_N(s)''$, and $\zeta(s)$

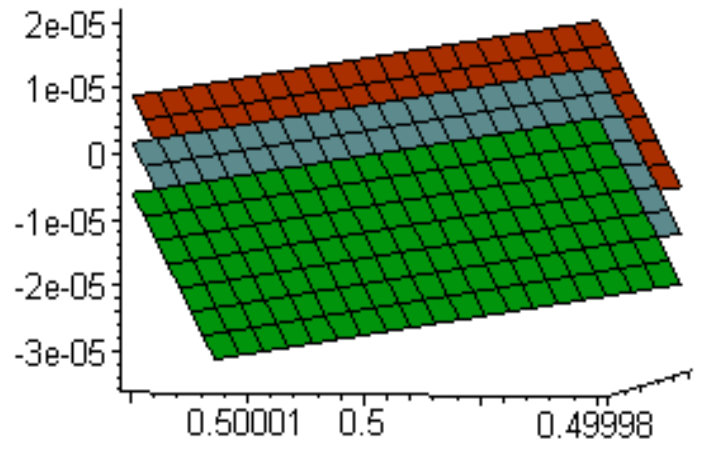

Fig 2. Imaginary parts of $\zeta_N(s)'$, $\zeta_N(s)''$, and $\zeta(s)$

This thus means that the real and imaginary parts of $\zeta_N(s)'$ will each intersect the complex plane along some lines at some distances away from $s_0$, as illustrated in Figure 3, and similarly will the real and imaginary parts of $\zeta_N(s)''$ also intersect the complex plane along some other lines at some distances away from $s_0$, as illustrated in Figure 4.

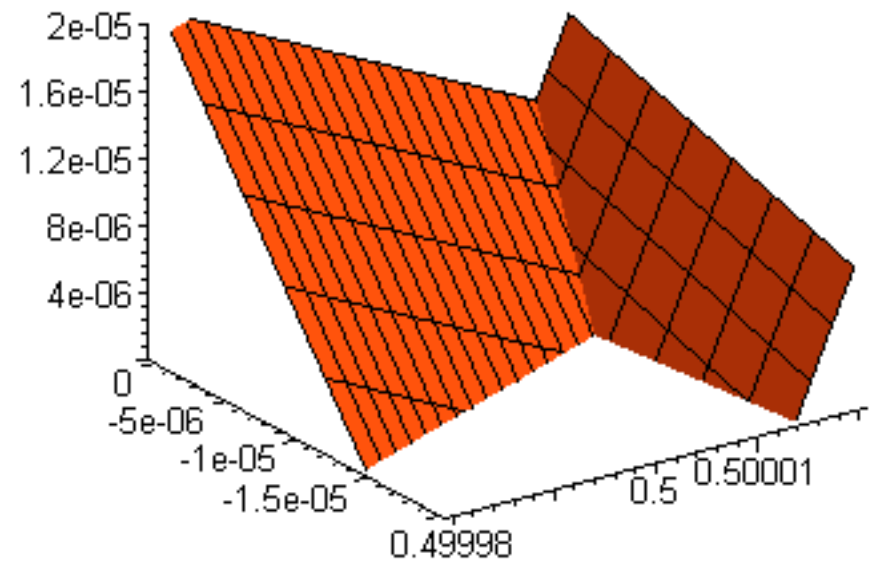

Fig 3. Real and imaginary parts of $\zeta_N(s)'$

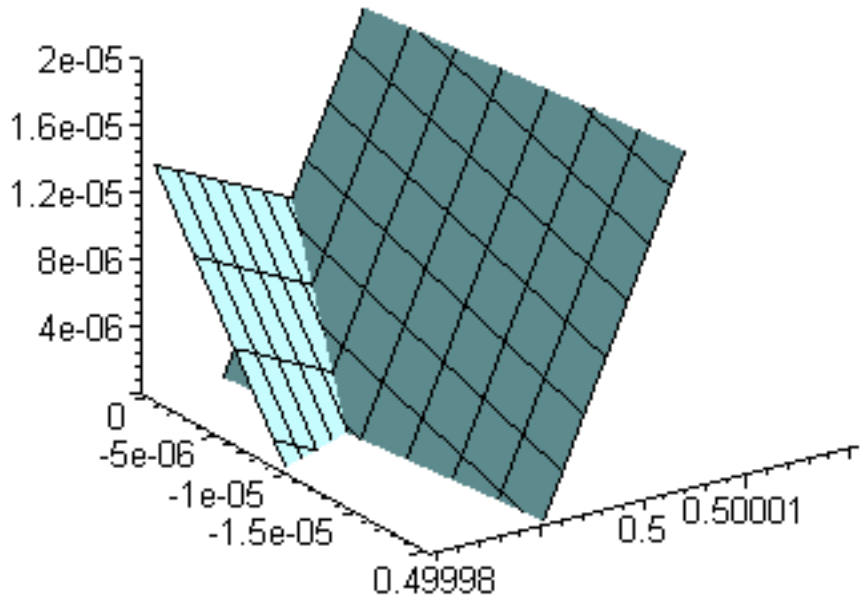

Fig 4. Real and imaginary parts of $\zeta_N(s)''$

In particular, there will as seen in Figure 3 above then be a point in the complex plane where the real and imaginary parts of $\zeta_N(s)'$ are both zero, and similarly some other point in the complex plane where the real and imaginary parts of $\zeta_N(s)''$ are both zero, as shown in Figure 4. At these two points at some distances away from the point $s = s_0$ where $|\zeta(s)| = 0$, we thus have that $|\zeta_N(s)'| = 0$ and $|\zeta_N(s)''| = 0$, as illustrated by the rather complicated structure in Figure 5 below, which depicts the absolute values of the functions $\zeta_N(s)'$ (red), $\zeta_N(s)''$ (turquoise), and $\zeta(s)$ (green).

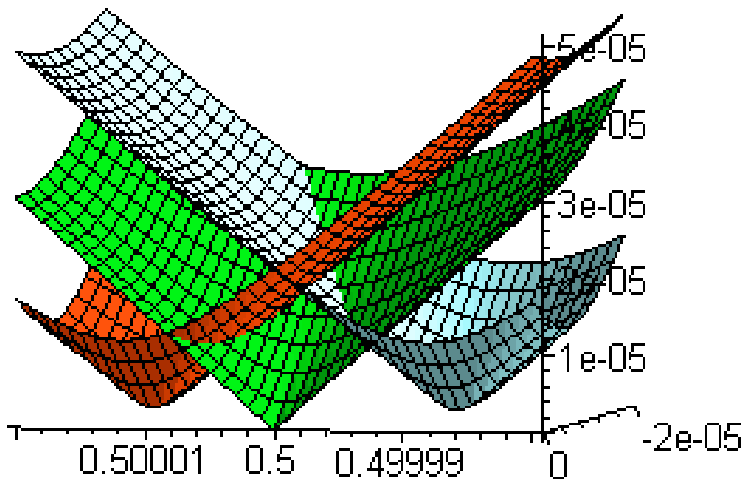

Fig 5. Absolute values of $\zeta_N(s)'$, $\zeta_N(s)''$, and $\zeta(s)$

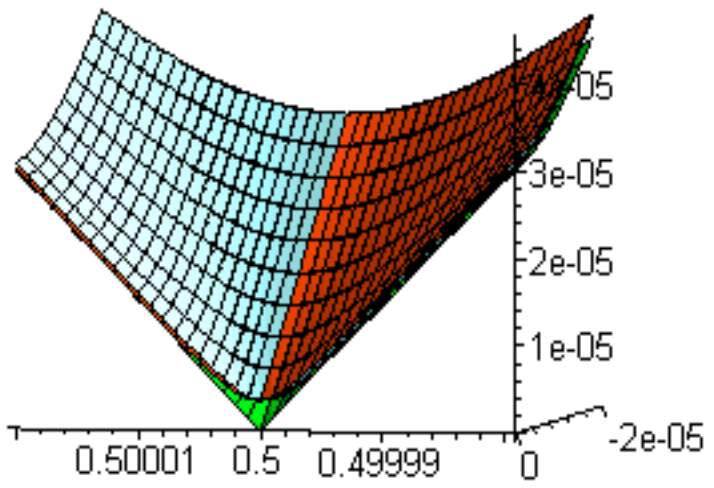

Fig 6. As Fig 5 but for larger $N$

Figure 6 shows how this structure contracts to an essentially conical structure around the closely spaced zeros/minima of the functions $|\zeta(s)|$, $|\zeta_N(s)'|$, and $|\zeta_N(s)''|$ for a larger value of $N$, and where the walls are a triple layer of these functions close together. For increasing $N$, the separations between these three layers become smaller and smaller, but the quotient $|\zeta_N(s)' - \zeta(s)| / |\zeta_N(s)'' - \zeta(s)|$ will nevertheless have some definite value even in the limit $N \rightarrow \infty$ as described in FAQ #7a.

---

## FAQ #8 (pages 8-9)

**"Can you please show me the details of how the rest of the expression on the left-hand side of (20) becomes equal to the middle expression (which is unity) in Remark 9.2 on pages 8-9 in your preprint."**

ANSWER: **Step1.** First I need to show that

$$|(s+2)(s+1)\,s\,\Gamma(s-3)| = |\Gamma(s)| \qquad (1)$$

for $s = ½ + i\,t$ as I say on top of page 9 in the preprint.

The following recurrence relation is valid for the gamma function (choose $n$ = 4 in Abramowitz and Stegun, eq 6.1.16)

$$(z+3)(z+2)(z+1)\,\Gamma(z+1) = \Gamma(4+z) \qquad (2)$$

Setting in particular $z$ = -7/2 + $i\,t$ in (2), and taking the absolute value of each factor, I get

$$\left|-\frac{1}{2}+i\,t\right|\left|-\frac{3}{2}+i\,t\right|\left|-\frac{5}{2}+i\,t\right|\left|\Gamma\!\left(-\frac{5}{2}+i\,t\right)\right| = \left|\Gamma\!\left(\frac{1}{2}+i\,t\right)\right| \qquad (3)$$

Without changing the values of the first three factors, I can replace the numbers within their absolute signs with their respective negative complex conjugates, and I then get

$$\left|\frac{1}{2}+i\,t\right|\left|\frac{3}{2}+i\,t\right|\left|\frac{5}{2}+i\,t\right|\left|\Gamma\left(-\frac{5}{2}+i\,t\right)\right|=\left|\Gamma\left(\frac{1}{2}+i\,t\right)\right| \qquad (4)$$

which becomes equal to (1) above if I set $s = ½ + i\,t$ there, and thus verifies the validity of (1).

For $s = ½ + i\,t$ the nontrivial factors in (20) can be calculated as follows.

**Step 2.** For $s = ½ + i\,t$ the cosine factor in the denominator in (20) becomes

$$\left|\cos\left(\frac{1}{2}\,s\,\pi\right)\right|=\left|\cos\left(\frac{1}{2}\left(\frac{1}{2}+i\,t\right)\pi\right)\right| \qquad (5)$$

Expand (5) in real and imaginary parts

$$\left|\cos\left(\frac{1}{2}\,s\,\pi\right)\right|=\left|\frac{1}{2}\sqrt{2}\,\cosh\left(\frac{1}{2}\pi\,t\right)-\frac{1}{2}\,i\sqrt{2}\,\sinh\left(\frac{1}{2}\pi\,t\right)\right| \qquad (6)$$

Form the absolute value

$$\left|\cos\left(\frac{1}{2}\,s\,\pi\right)\right|=\frac{1}{2}\sqrt{2\,\cosh\left(\frac{1}{2}\pi\,t\right)^2+2\,\sinh\left(\frac{1}{2}\pi\,t\right)^2} \qquad (7)$$

and simplify

$$\left|\cos\left(\frac{1}{2}\,s\,\pi\right)\right|=\frac{1}{2}\sqrt{2}\,\sqrt{\cosh(\pi\,t)} \qquad (8)$$

**Step 3.** For $s = ½ + i\,t$ the factor $|\,4^s - 8^s\,|$ in the numerator in (20) becomes

$$|\,4^s-8^s\,|=|\,4^s\,|\,|-1+2^s\,| \qquad (9)$$

$$|\,4^s-8^s\,|=|\,4^{(1/2+i\,t)}\,|\,|-1+2^{(1/2+i\,t)}\,| \qquad (10)$$

Expand the two factors on the right-hand side of (10) in real and imaginary parts

$$|\,4^s-8^s\,|=|\,2\cos(t\ln(4))+2\,i\sin(t\ln(4))\,|\,|-1+\sqrt{2}\cos(t\ln(2))+i\sqrt{2}\sin(t\ln(2))\,| \qquad (11)$$

Form the absolute values of the two factors on the right-hand side

$$|\,4^s-8^s\,|=2\sqrt{\cos(t\ln(4))^2+\sin(t\ln(4))^2}\,\sqrt{(-1+\sqrt{2}\cos(t\ln(2)))^2+2\sin(t\ln(2))^2} \qquad (12)$$

and simplify

$$|4^s - 8^s| = 2\sqrt{3 - 2\sqrt{2}\cos(t\ln(2))} \qquad (13)$$

**Step 4.** Finally the factor $|{-2}+2^s|$ in the denominator in (20) becomes for $s = ½ + i\,t$

$$|-2 + 2^s| = |-2 + 2^{(1/2 + i\,t)}| \qquad (14)$$

Expand (14) in real and imaginary parts

$$|-2 + 2^s| = |-2 + \sqrt{2}\cos(t\ln(2)) + i\sqrt{2}\sin(t\ln(2))| \qquad (15)$$

Form the absolute value

$$|-2 + 2^s| = \sqrt{(-2 + \sqrt{2}\cos(t\ln(2)))^2 + 2\sin(t\ln(2))^2} \qquad (16)$$

and simplify

$$|-2 + 2^s| = \sqrt{6 - 4\sqrt{2}\cos(t\ln(2))} \qquad (17)$$

**Step 5.** Inserting (8), (13), and (17) above into (20) in the preprint, and evaluating the trivial factors, we get

$$\frac{1}{2}\left|\frac{\pi^s\,(4^s - 8^s)}{\cos\left(\frac{1}{2}s\,\pi\right)(-2 + 2^s)\,\Gamma(s)}\right| = \frac{\sqrt{\frac{\pi}{\cosh(\pi\,t)}}}{\left|\Gamma\left(\frac{1}{2} + i\,t\right)\right|}$$

which thus verifies the equality of the first two members in Remark 9.2 in the preprint.

---

## FAQ #9 (page 9)

**"I think it would be more satisfactory to let the first integral in $I_N$ in Sect A1 on page 9 go from, say, -$\Lambda$ to $L$, so that the integral is indeed taken over a finite rectangle, and then afterwards let $\Lambda$ tend to infinity."**

ANSWER: Yes, but even though I agree on this point, I also think that the way it is done in my preprint is more visual. There are enough conceptual difficulties in the proof as it is, so I think allowing this particular description to be a little less abstract might help the reader somewhat at least in this part of the proof.

---

## FAQ #10 (page 10)

**"The expression for $I_1$ in Section A2 on page 10 in your preprint reads**

$$I_1 = (1 - e^{(2 i s \pi)})(1 - 2^{(1-s)})\,\Gamma(s)\,\zeta(s) - (1 - e^{(2 i s \pi)}) \int_{2\pi N}^{\infty} \frac{w^{(s-1)}}{e^w + 1}\, dw$$

**In your preprint you convert the last term on the right-hand side to the remainder $O(e^{-2\pi N})$. But for the factor in front of the second integral we have**

$$1 - e^{(2 i s \pi)} = (1 - e^{(2 i (\sigma + i t)\pi)}) = (1 - e^{(-2\pi t)}\, e^{(2 i \pi \sigma)})$$

**so when *t* tends to negative infinity this factor becomes infinite, and is thus hardly $O(e^{-2\pi N})$ as you state at the end of Section A2."**

ANSWER: Actually it is! The definition of a remainder O(f(*N*)) is that the remainder should be smaller in absolute terms than some constant (i e some number not containing *N*) times f(*N*). This is thus true for any given $s = \sigma + i\,t$. Even though the factor above tends to infinity with $-t$ as you point out, it is still a constant (albeit maybe a large one) from the point of view of *N*.

As *t* tends to negative infinity, both terms in $I_1$ (first term and remainder) tend to infinity. However, because of the exponential factor in the remainder,

$$\left| \int_{2N\pi}^{\infty} \frac{w^{(s-1)}}{e^w + 1}\, dw \right| < e^{(-2N\pi)}$$

the second term quickly becomes less and less important in $I_1$ compared to the first term as *N* gets larger, irrespective of the value of *t*, and can thus be neglected.

## FAQ #11 (page 11)

**"I get a little confused by all the explanations in Sect A3 in Appendix A in your preprint. Can't you give it in a more condensed mathematical form?"**

ANSWER: Define

$$\mathrm{h}(y, s, N) = \frac{\mathrm{e}^{(i\,s\,y)}}{\mathrm{e}^{(2\,\pi\,N\,(\cos(y) + i\,\sin(y)))} + 1} \tag{1}$$

and

$$\mathrm{g}(y, s) = \begin{cases} 0 & y < \frac{1}{2}\pi \\ \mathrm{e}^{(i\,s\,y)} & \frac{1}{2}\pi - y \le 0 \text{ and } y - \frac{3}{2}\pi \le 0 \\ 0 & \frac{3}{2}\pi < y \end{cases} \tag{2}$$

The integrals $I_2$ and $I_0$ in Sect A3 can then be written

$$I_2 = i\,2^s\,\pi^s\,N^s \int_0^{2\pi} \mathrm{h}(y, s, N)\,dy$$

and

$$I_0 = i\,2^s\,\pi^s\,N^s \int_0^{2\pi} \mathrm{g}(y, s)\,dy$$

We set

$$I_2 = I_0 + \Delta\,I_2\!\left(\frac{1}{2}\pi\right) + \Delta\,I_2\!\left(\frac{3}{2}\pi\right) \tag{3}$$

where the last two terms on the right-hand side are corrections to be determined, and of the form

$$\Delta\,I_2\!\left(\frac{1}{2}\pi\right) = \int_0^{2\pi} f_1(y, s, N)\,dy$$

$$\Delta\,I_2\!\left(\frac{3}{2}\pi\right) = \int_0^{2\pi} f_2(y, s, N)\,dy$$

where the two functions $f_1(y, s, N)$ and $f_2(y, s, N)$ are piecewise functions of type

$$f_1(y, s, N) = \{ \begin{matrix} \phi_1(y, s, N) & y \le \pi \\ 0 & \pi < y \end{matrix}$$

$$f_2(y, s, N) = \{ \begin{matrix} 0 & y \le \pi \\ \phi_2(y, s, N) & \pi < y \end{matrix}$$

and where the functions $\phi_1(y, s, N)$ and $\phi_2(y, s, N)$ are series expansions of the difference $i\,(2\pi N)^s\,[h(y, s, N) - g(y, s)]$ in the regions around $\pi/2$ and around $3\pi/2$, respectively [which regions with sufficient accuracy may be further confined to narrower regions of length $2\,\delta$ around these values, see (A4) in Appendix A in the preprint].

Inserting (1) and (2) into (3) then gives the final result for integral $I_2$ on the left-hand side, and where the integral in the first term on the right-hand side is easily calculated,

$$i\,2^s\,\pi^s\,N^s \int_0^{2\,\pi} \frac{e^{(i\,s\,y)}}{e^{(2\,\pi\,N\,(\cos(y) + i\,\sin(y)))} + 1}\,dy = i\,2^s\,\pi^s\,N^s \int_{1/2\,\pi}^{3/2\,\pi} e^{(i\,s\,y)}\,dy + \Delta\,I_2\left(\frac{1}{2}\,\pi\right) + \Delta\,I_2\left(\frac{3}{2}\,\pi\right)$$

---

## FAQ #11a (page 11)

**"I have difficulties to see the point with Section A3 in Appendix A in your preprint. Can you please explain what happens there?"**

ANSWER: OK, what happens in Section A3 is the following. I wish to calculate the integral $I_2$ in (A1). In order to do that I consider an approximating, piecewise function, which is zero for $0 \le y < \pi/2$ and for $3\pi/2 < y \le 2\pi$, whereas for $\pi/2 \le y \le 3\,\pi/2$ it is given by the exponential in the numerator in (A1). This piecewise function gives the integral $I_0$ in (A2).

In order to calculate $I_2$, I write

$$I_2 = I_0 + (I_2 - I_0).$$

The point with writing it this way is that the integrands in $I_2$ and $I_0$ differ essentially only in two rather narrow regions around $\pi/2$ and $3\pi/2$, respectively, as shown in the plots below. So I designate $(I_2 - I_0)$ around $\pi/2$ by $\Delta I_2(\pi/2)$, and $(I_2 - I_0)$ around $3\pi/2$ by $\Delta I_2(3\pi/2)$. I can thus write

$$I_2 = I_0 + \Delta I_2(\pi/2) + \Delta I_2(3\pi/2),$$

which can thus be regarded merely as a definition of $\Delta I_2(\pi/2)$ and $\Delta I_2(3\pi/2)$, provided that I also state that $\Delta I_2(\pi/2)$ designates the difference $I_2 - I_0$ for $y < \pi$, whereas $\Delta I_2(3\pi/2)$ designates this difference for $y > \pi$.

As I mentioned above, the point with this is thus that the integrands of $I_2$ and $I_0$ differ essentially only in two rather narrow regions around $\pi/2$ and $3\pi/2$, respectively. So I can make series expansions of the difference between the integrands of $I_2$ and $I_0$ in these two regions and calculate the integrals there. This is thus again nothing more than using the definition above and calculating

$\Delta I_2(\pi/2) = I_2 - I_0$ for $y < \pi$ by a series expansion around $y = \pi/2$, and

$\Delta I_2(3\pi/2) = I_2 - I_0$ for $y > \pi$ by a series expansion around $y = 3\pi/2$.

Around $\pi/2$ we thus end up with the first expression in Sect A3.2 on page 14 in the preprint, and around $3\pi/2$ we end up with the second expression in Sect A3.2, which added together and to the expression for $I_0$ give the final result on top of page 15 in the preprint.

Plot 1: Shows the integrand in $I_2 - I_0$ for small $N$ ($N = 5$), so that the localised deviations around $\pi/2$ and $3\pi/2$ are clearly visible (both real and imaginary parts are shown):

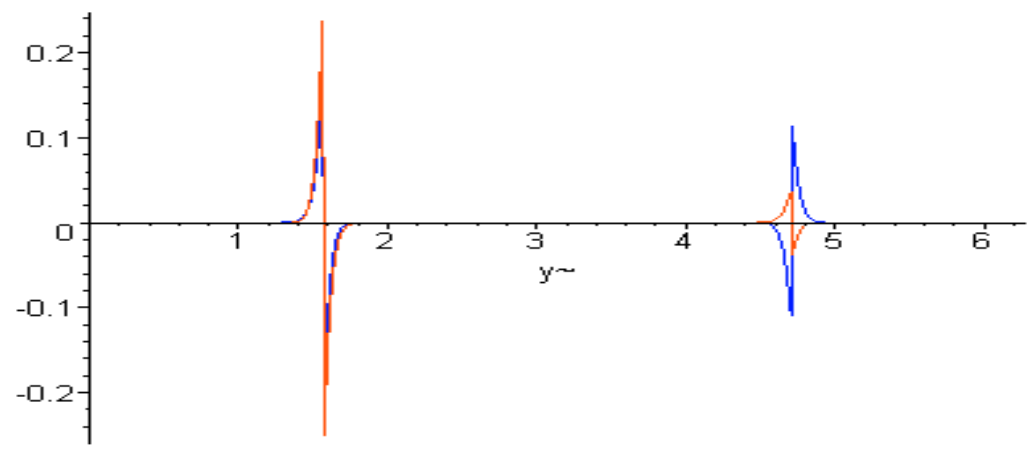

Plot 2: Shows the integrand in $I_2 - I_0$ for $N = 1000$ immediately around $\pi/2$. Note the different scale, and how thus the localised peaks get narrower for larger values of $N$:

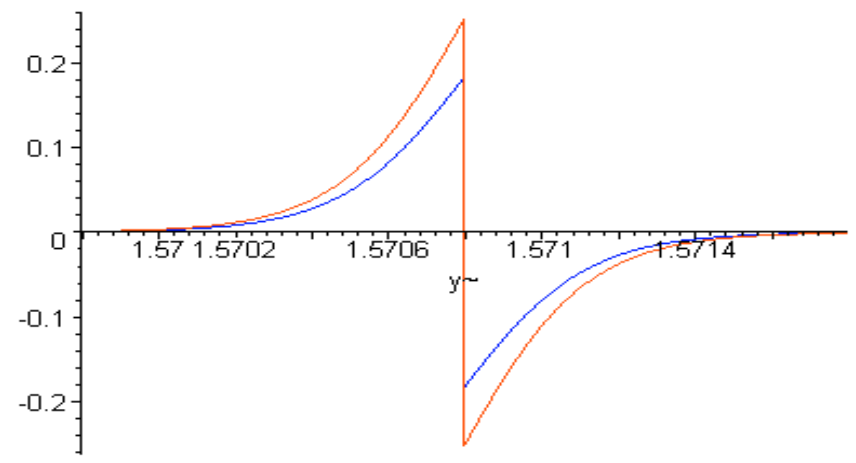

Plot 3: Corresponding plot of the integrand in $I_2 - I_0$ for $N = 1000$ immediately around $3\pi/2$:

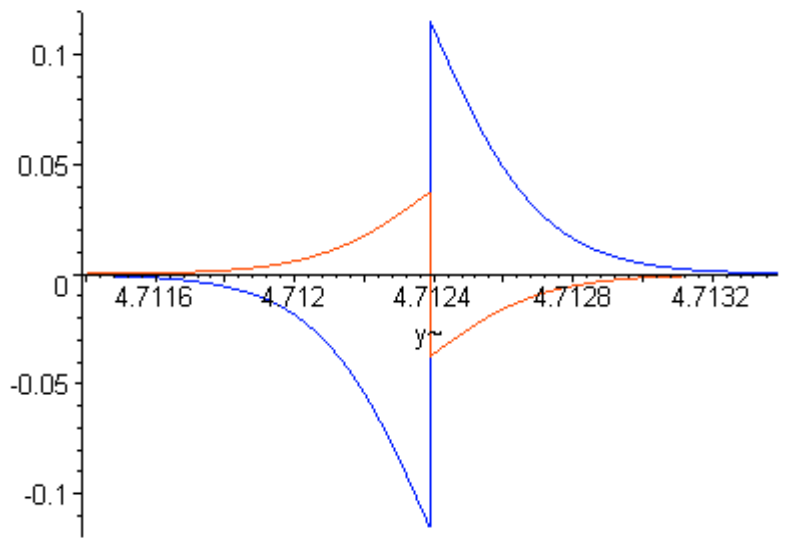

## FAQ #12 (page 11)

**"There seems to be some detailed analysis hidden in making the factor $e^{\frac{1}{2} i s \pi}$ appear outside the integral in the equation at the bottom of page 11 in your preprint."**

<u>ANSWER:</u> Straightforward Taylor expansion around $\pi/2$ of the exponential in the numerator in (A1) gives for the first few terms

$$e^{(i s y)} = e^{(1/2\, i s \pi)} + i\, e^{(1/2\, i s \pi)}\, s \left(y - \frac{1}{2}\pi\right) - \frac{1}{2} e^{(1/2\, i s \pi)}\, s^2 \left(y - \frac{1}{2}\pi\right)^2$$

Since the common factor $e^{\frac{1}{2} i s \pi}$ in all terms is independent of $y$, we can put it outside the integral and thus get the equation at the bottom of page 11.

---

## FAQ #13 (page 12)

**"Looking at the plots in FAQ #11a, I'm amazed that series expansions of just order two as you do on page 12 in your preprint could really be sufficient to describe the behavior of the peaks around $\pi/2$ and $3\pi/2$, respectively. Is that really possible?"**

<u>ANSWER:</u> Yes, because please remember that the series expansion is in $1/N$, not in the variable $\tau$ (or $y$), which is the variable in which the behavior in the plots above are described. As a function of $\tau$ the peaks are well described by the rather complicated expressions $a_1$ and $a_2$ at the top of page 13, which are the coefficients for the series expansion in $1/N$ that I make, and where order two in $1/N$ turns out to be sufficient.

---

## FAQ #14 (page 12)

**"I have a question on the last formula in Sect A3.1 in Appendix A in your preprint.**

$$I_2\left(\frac{1}{2}\pi, \delta\right) = i\, e^{(1/2\, i s \pi)}\, 2^s\, \pi^s\, N^s \int_{-N\delta}^{N\delta} \frac{\frac{1}{N} + \frac{i s \tau}{N^2} + O\left(\frac{1}{N^3}\right)}{e^{\left(-2\pi\tau - \frac{i \pi \tau^2}{N}\right)} + 1}\, d\tau$$

**What happens to the higher terms in the exponential in the denominator? I can understand how terms with powers less than 2 in $1/N$ become as shown in the formula, but what happens to the higher powers in $1/N$?"**

<u>ANSWER:</u> When you are used to working with Landau O:s, you do these things more or less automatically. There may be ways to see it simpler, but here is one sequence of steps to get to the formula in the preprint.

The higher powers are shown below in the form of remainders in the numerator and in the exponent in the denominator,

$$I_2\left(\frac{1}{2}\pi, \delta\right) = i\, e^{(1/2\, i\, s\, \pi)}\, 2^s\, \pi^s\, N^s \int_{-N\delta}^{N\delta} \frac{\frac{1}{N} + \frac{i\, s\, \tau}{N^2} + O\left(\frac{1}{N^3}\right)}{e^{\left(-2\pi\tau - \frac{i\pi\tau^2}{N} + O\left(\frac{1}{N^2}\right)\right)} + 1}\, d\tau$$

Separate out the remainder in the denominator as a separate factor,

$$I_2\left(\frac{1}{2}\pi, \delta\right) = i\, e^{(1/2\, i\, s\, \pi)}\, 2^s\, \pi^s\, N^s \int_{-N\delta}^{N\delta} \frac{\frac{1}{N} + \frac{i\, s\, \tau}{N^2} + O\left(\frac{1}{N^3}\right)}{e^{\left(-\frac{i\pi\tau^2}{N} - 2\pi\tau\right)}\, e^{O\left(\frac{1}{N^2}\right)} + 1}\, d\tau$$

Make a series expansion of this exponential of the remainder,

$$I_2\left(\frac{1}{2}\pi, \delta\right) = i\, e^{(1/2\, i\, s\, \pi)}\, 2^s\, \pi^s\, N^s \int_{-N\delta}^{N\delta} \frac{\frac{1}{N} + \frac{i\, s\, \tau}{N^2} + O\left(\frac{1}{N^3}\right)}{e^{\left(-\frac{i\pi\tau^2}{N} - 2\pi\tau\right)}\left(1 + O\left(\frac{1}{N^2}\right)\right) + 1}\, d\tau$$

Simplify the denominator,

$$I_2\left(\frac{1}{2}\pi, \delta\right) = i\, e^{(1/2\, i\, s\, \pi)}\, 2^s\, \pi^s\, N^s \int_{-N\delta}^{N\delta} \frac{\frac{1}{N} + \frac{i\, s\, \tau}{N^2} + O\left(\frac{1}{N^3}\right)}{e^{\left(-\frac{i\pi\tau^2}{N} - 2\pi\tau\right)} + 1 + O\left(\frac{1}{N^2}\right)}\, d\tau$$

Write the remainder in the denominator as a factor,

$$I_2\left(\frac{1}{2}\pi,\ \delta\right) = i\ e^{(1/2\ i\ s\ \pi)}\ 2^s\ \pi^s\ N^s \int_{-N\delta}^{N\delta} \frac{\frac{1}{N} + \frac{i\ s\ \tau}{N^2} + O\left(\frac{1}{N^3}\right)}{\left(e^{\left(-\frac{i\pi\tau^2}{N} - 2\pi\tau\right)} + 1\right)\left(1 + O\left(\frac{1}{N^2}\right)\right)}\, d\tau$$

Make a series expansion of 1/(this factor) so that it ends up in the numerator,

$$I_2\left(\frac{1}{2}\pi,\ \delta\right) = i\ e^{(1/2\ i\ s\ \pi)}\ 2^s\ \pi^s\ N^s \int_{-N\delta}^{N\delta} \frac{\left(1 + O\left(\frac{1}{N^2}\right)\right)\left(\frac{1}{N} + \frac{i\ s\ \tau}{N^2} + O\left(\frac{1}{N^3}\right)\right)}{e^{\left(-\frac{i\pi\tau^2}{N} - 2\pi\tau\right)} + 1}\, d\tau$$

Multiplying together the two factors in the numerator gives the final expression,

$$I_2\left(\frac{1}{2}\pi,\ \delta\right) = i\ e^{(1/2\ i\ s\ \pi)}\ 2^s\ \pi^s\ N^s \int_{-N\delta}^{N\delta} \frac{\frac{1}{N} + \frac{i\ s\ \tau}{N^2} + O\left(\frac{1}{N^3}\right)}{e^{\left(-\frac{i\pi\tau^2}{N} - 2\pi\tau\right)} + 1}\, d\tau$$

---

## FAQ #14a (page 12)

**"I am a little uneasy about the integrations in Sect A3 in the preprint. In the expression for $I_2(\frac{1}{2}\pi,\ \delta)$ in the first formula in Sect A3.1, the integration is performed over a small interval around $\pi/2$. But then after changing the integration variable in Sect A3.1 and the special choice in (A4), the integration of $\Delta I_2(\frac{1}{2}\pi,\ \delta)$ in the first formula in Sect A3.11 is now to be performed over an interval of the order of $N\delta = \ln(N)$, which thus tends to infinity with $N$. So one needs to be sure that the integrals of the series expansions really converge over this interval. Please then also remember that although the real part $Re(s)$ is assumed to be $< 1$, still the imaginary part $Im(s)$ must be permitted to be a very large positive or negative number."**

<u>ANSWER:</u> Consider here first the numerator of the integrand in the expression for $I_2(\frac{1}{2}\pi,\ \delta)$ in the first formula in Sect A3.1, i e (here given in explicit form to order three)

$$1 + i\ s\left(y - \frac{1}{2}\pi\right) - \frac{1}{2}s^2\left(y - \frac{1}{2}\pi\right)^2 - \frac{1}{6}i\ s^3\left(y - \frac{1}{2}\pi\right)^3 + \ldots$$

After changing integration variable as in Sect A3.1, the expression above becomes (to order three)

$$\frac{1}{N}+\frac{i\,s\,\tau}{N^2}-\frac{1}{2}\,\frac{s^2\,\tau^2}{N^3}-\frac{\frac{1}{6}\,i\,s^3\,\tau^3}{N^4}$$

where the common factor $1/N$ results from the change from $dy$ to $d\tau$. Since $\tau$ is at the most equal to $N\delta$, i e at the most equal to ln($N$) from (A4), we need to consider a remainder due to a sum of terms of type

$$\frac{s^p\,\tau^p}{N^{(p+1)}}=\frac{s^p\,\ln(N)^p}{N^{(p+1)}}$$

Thus no matter how large a given $Im(s)$ is, we can always choose an $N$ such that the remainder for some given order $p$ is arbitrarily small. Specifically, if we wish to consider a first order expansion, as is enough for the purpose of the preprint, we have (remember that there is a common factor $1/N$)

$$\frac{1}{N}+\frac{i\,s\,\tau}{N^2}+\mathrm{O}\left(\frac{1}{N^3}\right)$$

as given for the numerator in the preprint. The small deviation from a pure power in $N$ due to possible factors ln($N$) in the remainders as discussed above is denoted by a decimal point in the power of $N$ in the remainder in the final integral, as defined Sect A3.13 in the preprint.

Consider next the expansion in the exponential in the denominator of the integrand in the expression for $I_2(\frac{1}{2}\pi,\,\delta)$ in the above first formula in Sect A3.1. Given explicitly to order four, the expansion is as follows

$$2\,i\,\pi\,N-2\,\pi\,N\left(y-\frac{1}{2}\,\pi\right)-i\,\pi\,N\left(y-\frac{1}{2}\,\pi\right)^2+\frac{1}{3}\,\pi\,N\left(y-\frac{1}{2}\,\pi\right)^3+\frac{1}{12}\,i\,\pi\,N\left(y-\frac{1}{2}\,\pi\right)^4 + \ldots$$

After again changing the integration variable as in Sect A3.1, this expansion becomes

$$2\,i\,\pi\,N-2\,\pi\,\tau-\frac{i\,\pi\,\tau^2}{N}+\frac{\frac{1}{3}\,\pi\,\tau^3}{N^2}+\frac{\frac{1}{12}\,i\,\pi\,\tau^4}{N^3}+\ldots$$

Here again, $\tau$ is at the most equal to $N\delta$, i e at the most equal to ln($N$) from (A4). Thus we get the following expression for the exponential in the denominator (to first order in $1/N$, which is sufficient for the analysis in the preprint)

$$2\,i\,\pi\,N-2\,\pi\,\tau-\frac{i\,\pi\,\tau^2}{N}+\mathrm{O}\left(\frac{1}{N^2}\right)$$

Exponentiation of the first term gives a trivial factor unity as noted in the preprint. The remainder above [again including possible factors $\tau < N\,\delta = \ln(N)$ as above in the final integration] can be incorporated into the remainder in the numerator as discussed in detail in FAQ #14 above and in Sect A3.13 in the preprint.

As seen from the above calculations, there is thus no problem with the convergence of the integrals in Sect A3. This is as should be expected since, as discussed in FAQ #11a above, the integral $I_2$ is equal to the integral $I_0$ with the analytical expression given in (A3) plus a small correction.

---

## FAQ #15 (page 13)

**"I do not see how the integrals on the left-hand side of (A7) on page 13 in your preprint become the expression on the right-hand side with $\varepsilon_V(s)$ given in (A8)."**

ANSWER: This requires some calculations. We want to calculate the left-hand side of (A7)

$$Y = \int_{-N\delta}^{0} \frac{a_2}{N^2}\, d\tau + \int_{0}^{N\delta} \frac{a_2}{N^2} - \frac{i\,s\,\tau}{N^2}\, d\tau$$

Insert $a_2$ from page 13 in my preprint

$$Y = \int_{-N\delta}^{0} \frac{\dfrac{i\,s\,\tau}{e^{(-2\pi\tau)}+1} + \dfrac{i\,e^{(-2\pi\tau)}\,\pi\,\tau^2}{2\,(e^{(-2\pi\tau)}+1)}}{N^2}\, d\tau + \int_{0}^{N\delta} \frac{\dfrac{i\,s\,\tau}{e^{(-2\pi\tau)}+1} + \dfrac{i\,e^{(-2\pi\tau)}\,\pi\,\tau^2}{2\,(e^{(-2\pi\tau)}+1)}}{N^2} - \frac{i\,s\,\tau}{N^2}\, d\tau$$

Let Maple integrate

$$\begin{aligned}
Y = {} & \frac{1}{48} i\,(s\,\pi^2 + \pi^2 s\, e^{(2N\delta\pi)} - \pi^2 - \pi^2 e^{(2N\delta\pi)} - 24 s\,\pi^2 N^2 \delta^2 + 24 s \ln(1+e^{(2N\delta\pi)})\,\pi N\delta \\
& - 12\,\mathrm{dilog}(1+e^{(2N\delta\pi)})\, e^{(2N\delta\pi)} - 12\,\mathrm{dilog}(1+e^{(2N\delta\pi)}) + 24 s \ln(1+e^{(2N\delta\pi)})\,\pi N\delta\, e^{(2N\delta\pi)} \\
& - 24 \ln(1+e^{(2N\delta\pi)})\,\pi N\delta - 24\ln(1+e^{(2N\delta\pi)})\,\pi N\delta\, e^{(2N\delta\pi)} + 24\pi^2 N^2\delta^2 e^{(2N\delta\pi)} \\
& - 24 s\,\pi^2 N^2 \delta^2 e^{(2N\delta\pi)} + 12 s\,\mathrm{dilog}(1+e^{(2N\delta\pi)})\, e^{(2N\delta\pi)} + 12 s\,\mathrm{dilog}(1+e^{(2N\delta\pi)}))\Big/(N^2 \\
& (1+e^{(2N\delta\pi)})\,\pi^2) - \frac{1}{48} i\,(24 e^{(2N\delta\pi)} \ln(e^{(-2N\delta\pi)}+1)\,\delta N\pi - 24 e^{(2N\delta\pi)} s \ln(e^{(-2N\delta\pi)}+1)\,\delta N\pi \\
& + 24\pi^2 N^2\delta^2 - 24 s\ln(e^{(-2N\delta\pi)}+1)\,\pi N\delta - 12 e^{(2N\delta\pi)}\,\mathrm{dilog}((1+e^{(2N\delta\pi)})\, e^{(-2N\delta\pi)}) \\
& + 12 e^{(2N\delta\pi)} s\,\mathrm{dilog}((1+e^{(2N\delta\pi)})\, e^{(-2N\delta\pi)}) + 24\ln(e^{(-2N\delta\pi)}+1)\,\pi N\delta \\
& + 12 s\,\mathrm{dilog}((1+e^{(2N\delta\pi)})\, e^{(-2N\delta\pi)}) - 12\,\mathrm{dilog}((1+e^{(2N\delta\pi)})\, e^{(-2N\delta\pi)}) + s\,\pi^2 + \pi^2 s\, e^{(2N\delta\pi)} - \pi^2 \\
& - \pi^2 e^{(2N\delta\pi)})\Big/(N^2 (1+e^{(2N\delta\pi)})\,\pi^2)
\end{aligned}$$

and also simplify

$$
\begin{aligned}
\mathrm{Y} = & \frac{-1}{4} i\,(6 s \pi^2 N^2 \delta^2 + 4\ln(1+\mathrm{e}^{(2N\delta\pi)})\,\pi N\delta - 6\pi^2 N^2 \delta^2 \mathrm{e}^{(2N\delta\pi)} + \mathrm{dilog}(1+\mathrm{e}^{(2N\delta\pi)})\,\mathrm{e}^{(2N\delta\pi)} \\
& - s\,\mathrm{dilog}(1+\mathrm{e}^{(2N\delta\pi)}) + \mathrm{dilog}(1+\mathrm{e}^{(2N\delta\pi)}) + 4\ln(1+\mathrm{e}^{(2N\delta\pi)})\,\pi N\delta\,\mathrm{e}^{(2N\delta\pi)} \\
& + 6 s \pi^2 N^2 \delta^2 \mathrm{e}^{(2N\delta\pi)} - 4 s\ln(1+\mathrm{e}^{(2N\delta\pi)})\,\pi N\delta - 4 s \ln(1+\mathrm{e}^{(2N\delta\pi)})\,\pi N\delta\,\mathrm{e}^{(2N\delta\pi)} - 2\pi^2 N^2 \delta^2 \\
& + \mathrm{e}^{(2N\delta\pi)}\, s\,\mathrm{dilog}((1+\mathrm{e}^{(2N\delta\pi)})\,\mathrm{e}^{(-2N\delta\pi)}) - s\,\mathrm{dilog}(1+\mathrm{e}^{(2N\delta\pi)})\,\mathrm{e}^{(2N\delta\pi)} \\
& - \mathrm{e}^{(2N\delta\pi)}\,\mathrm{dilog}((1+\mathrm{e}^{(2N\delta\pi)})\,\mathrm{e}^{(-2N\delta\pi)}) - \mathrm{dilog}((1+\mathrm{e}^{(2N\delta\pi)})\,\mathrm{e}^{(-2N\delta\pi)}) \\
& + s\,\mathrm{dilog}((1+\mathrm{e}^{(2N\delta\pi)})\,\mathrm{e}^{(-2N\delta\pi)}))\ \big/\ (N^2\,(1+\mathrm{e}^{(2N\delta\pi)})\,\pi^2)
\end{aligned}
$$

This can be simplified further by defining

$$
\Phi = \frac{-i\,\delta^2}{1+\mathrm{e}^{(2N\delta\pi)}}
$$

since then we can write $\mathrm{Y} - \Phi$ as follows

$$
\mathrm{Y}-\Phi = \frac{\frac{1}{4} i\,(4\ln(1+\mathrm{e}^{(2N\delta\pi)})\,\pi N\delta - \mathrm{dilog}((1+\mathrm{e}^{(2N\delta\pi)})\,\mathrm{e}^{(-2N\delta\pi)}) + \mathrm{dilog}(1+\mathrm{e}^{(2N\delta\pi)}) - 6\pi^2 N^2 \delta^2)}{\pi^2 N^2}
$$

We want Y on the left-hand side, so we form

$$
\begin{aligned}
& -\frac{i N^2 (\mathrm{Y}-\Phi)}{s-1} - \frac{i N^2 \Phi}{s-1} = \\
& \quad \frac{1}{4}\,\frac{4\ln(1+\mathrm{e}^{(2N\delta\pi)})\,\pi N\delta - \mathrm{dilog}((1+\mathrm{e}^{(2N\delta\pi)})\,\mathrm{e}^{(-2N\delta\pi)}) + \mathrm{dilog}(1+\mathrm{e}^{(2N\delta\pi)}) - 6\pi^2 N^2\delta^2}{\pi^2} \\
& \quad - \frac{N^2 \delta^2}{(s-1)\,(1+\mathrm{e}^{(2N\delta\pi)})}
\end{aligned}
$$

because this simplifies to

$$
\begin{aligned}
\frac{-i N^2\,\mathrm{Y}}{s-1} = & \frac{\ln(1+\mathrm{e}^{(2N\delta\pi)})\,N\delta}{\pi} - \frac{1}{4}\,\frac{\mathrm{dilog}(\mathrm{e}^{(-2N\delta\pi)}+1)}{\pi^2} + \frac{\frac{1}{4}\,\mathrm{dilog}(1+\mathrm{e}^{(2N\delta\pi)})}{\pi^2} - \frac{3}{2} N^2 \delta^2 \\
& - \frac{N^2\delta^2}{(s-1)\,(1+\mathrm{e}^{(2N\delta\pi)})}
\end{aligned}
$$

Substitute $N\delta \to \nu$

$$
\frac{-i N^2\,\mathrm{Y}}{s-1} = \frac{1}{4}\,\frac{\mathrm{dilog}(\mathrm{e}^{(2\nu\pi)}+1)}{\pi^2} - \frac{1}{4}\,\frac{\mathrm{dilog}(\mathrm{e}^{(-2\nu\pi)}+1)}{\pi^2} + \frac{\ln(\mathrm{e}^{(2\nu\pi)}+1)\,\nu}{\pi} - \frac{3}{2}\nu^2 - \frac{\nu^2}{(s-1)\,(\mathrm{e}^{(2\nu\pi)}+1)}
$$

Defining

$$\varepsilon_\nu(s) = \frac{1}{4}\frac{\operatorname{dilog}(e^{(2\nu\pi)}+1)}{\pi^2} - \frac{1}{4}\frac{\operatorname{dilog}(e^{(-2\nu\pi)}+1)}{\pi^2} + \frac{\ln(e^{(2\nu\pi)}+1)\,\nu}{\pi} - \frac{3}{2}\nu^2 - \frac{\nu^2}{(s-1)(e^{(2\nu\pi)}+1)}$$

we can thus write Y and (A7) as

$$\int_{-N\delta}^{0} \frac{a_2}{N^2}\,d\tau + \int_{0}^{N\delta} \frac{a_2}{N^2} - \frac{i\,s\,\tau}{N^2}\,d\tau = \frac{i\,\varepsilon_\nu(s)\,(s-1)}{N^2}$$

as given in my preprint.

---

## FAQ #15a (page 13)

**"Your derivation in FAQ #15 above is not particularly transparent, and it relies completely on that Maple has done its job correctly, which is not very satisfactory. I think this is a weak point in your proof."**

ANSWER: No actually it isn't, because once we have derived the expression $\varepsilon_\nu(s)$ as in FAQ #15 above,

$$\varepsilon_\nu(s) = \frac{1}{4}\frac{\operatorname{dilog}(e^{(2\nu\pi)}+1)}{\pi^2} - \frac{1}{4}\frac{\operatorname{dilog}(e^{(-2\nu\pi)}+1)}{\pi^2} + \frac{\ln(e^{(2\nu\pi)}+1)\,\nu}{\pi} - \frac{3}{2}\nu^2 - \frac{\nu^2}{(s-1)(e^{(2\nu\pi)}+1)}$$

then we can check its correctness as follows.

(A7) in my preprint can be written

$$\int_{-\nu}^{0} \frac{a_2(\tau)}{N^2}\,d\tau + \int_{0}^{\nu} \frac{a_2(\tau)}{N^2} - \frac{i\,s\,\tau}{N^2}\,d\tau = \frac{i\,\varepsilon_\nu(s)\,(s-1)}{N^2} \qquad (A7)$$

where from page 13 in my preprint we have

$$a_2(\tau) = \frac{i\,s\,\tau}{e^{(-2\pi\tau)}+1} + \frac{i\,e^{(-2\pi\tau)}\,\pi\,\tau^2}{(e^{(-2\pi\tau)}+1)^2}$$

Now differentiate both sides in (A7) as given above with respect to ν. The left-hand side (*LHS*) then becomes

$$LHS = \frac{a_2(-\nu)}{N^2} + \frac{a_2(\nu)}{N^2} - \frac{i\,s\,\nu}{N^2}$$

or, after inserting $a_2$ from above,

$$LHS = \frac{-\frac{i\,s\,\nu}{e^{(2\,\nu\,\pi)}+1} + \frac{i\,e^{(2\,\nu\,\pi)}\,\pi\,\nu^2}{(e^{(2\,\nu\,\pi)}+1)^2}}{N^2} + \frac{\frac{i\,s\,\nu}{e^{(-2\,\nu\,\pi)}+1} + \frac{i\,e^{(-2\,\nu\,\pi)}\,\pi\,\nu^2}{(e^{(-2\,\nu\,\pi)}+1)^2}}{N^2} - \frac{i\,s\,\nu}{N^2}$$

which simplifies to

$$LHS = \frac{2\,i\,\nu\,(-s\,e^{(2\,\nu\,\pi)} - s + \nu\,\pi\,e^{(2\,\nu\,\pi)})}{N^2\,(e^{(2\,\nu\,\pi)}+1)^2}$$

We now similarly differentiate the right-hand side in (A7) as given above with respect to ν. For the factor $\varepsilon_\nu(s)$ we just calculated, we then get

$$\frac{\partial}{\partial\nu}\varepsilon_\nu(s) = -3\,\nu + \frac{2\,e^{(2\,\nu\,\pi)}\,\nu}{e^{(2\,\nu\,\pi)}+1} + \frac{\frac{1}{2}\ln(e^{(2\,\nu\,\pi)}+1)}{\pi} - \frac{1}{2}\,\frac{\ln(e^{(-2\,\nu\,\pi)}+1)}{\pi} - \frac{2\,\nu}{(s-1)\,(e^{(2\,\nu\,\pi)}+1)} + \frac{2\,\nu^2\,\pi\,e^{(2\,\nu\,\pi)}}{(s-1)\,(e^{(2\,\nu\,\pi)}+1)^2}$$

or

$$\frac{\partial}{\partial\nu}\varepsilon_\nu(s) = -3\,\nu + \frac{2\,e^{(2\,\nu\,\pi)}\,\nu}{e^{(2\,\nu\,\pi)}+1} + \frac{\ln\left(\frac{\sqrt{e^{(2\,\nu\,\pi)}+1}}{\sqrt{e^{(-2\,\nu\,\pi)}+1}}\right)}{\pi} - \frac{2\,\nu}{(s-1)\,(e^{(2\,\nu\,\pi)}+1)} + \frac{2\,\nu^2\,\pi\,e^{(2\,\nu\,\pi)}}{(s-1)\,(e^{(2\,\nu\,\pi)}+1)^2}$$

which since

$$\ln\left(\frac{\sqrt{e^{(2\,\nu\,\pi)}+1}}{\sqrt{e^{(-2\,\nu\,\pi)}+1}}\right) = \nu\,\pi$$

can be simplified to

$$\frac{\partial}{\partial\nu}\varepsilon_\nu(s) = -2\,\nu + \frac{2\,e^{(2\,\nu\,\pi)}\,\nu}{e^{(2\,\nu\,\pi)}+1} - \frac{2\,\nu}{(s-1)\,(e^{(2\,\nu\,\pi)}+1)} + \frac{2\,\nu^2\,\pi\,e^{(2\,\nu\,\pi)}}{(s-1)\,(e^{(2\,\nu\,\pi)}+1)^2}$$

After some further simplification, the derivative (*RHS*) of the right-hand side becomes

$$RHS = \frac{i\,(s-1)\left(\frac{\partial}{\partial \nu}\varepsilon_\nu(s)\right)}{N^2} = \frac{2\,i\,\nu\,(-s\,e^{(2\,\nu\,\pi)} - s + \nu\,\pi\,e^{(2\,\nu\,\pi)})}{N^2\,(e^{(2\,\nu\,\pi)}+1)^2}$$

Thus the left-hand side *LHS* and right-hand side *RHS* are equal, and the expression for the derivatives of both sides of (A7) are thus equal, so the left-hand side and the right-hand side in (A7) can hence differ by at most a constant, which we can demonstrate to be zero as follows. For $\nu = 0$ the left-hand side of (A7) becomes equal to 0, and so does also the above expression for $\varepsilon_\nu(s)$, and thus also the right-hand side of (A7).

Thus the left-hand side and the right-hand side of (A7) are equal, which thus proves that the expression derived for $\varepsilon_\nu(s)$ above in FAQ #15 is correct.

---

## FAQ #16 (page 13)

**"How do you approximate (A8) on page 13 in your preprint to the expression following it?"**

ANSWER: Set

$$e^{(2\,\nu\,\pi)} = X \qquad\qquad e^{(-2\,\nu\,\pi)} = \frac{1}{X}$$

Then (A8) can be written

$$\varepsilon_\nu(s) = -\frac{1}{4}\,\frac{\mathrm{dilog}\left(\frac{1}{X}+1\right)}{\pi^2} + \frac{\frac{1}{4}\,\mathrm{dilog}(X+1)}{\pi^2} + \frac{\ln(X+1)\,\nu}{\pi} - \frac{3}{2}\nu^2 - \frac{\nu^2}{(s-1)\,(X+1)}$$

The following series expansions are valid for $X >> 1$ (cf FAQ #16a below)

$$\mathrm{dilog}\left(\frac{1}{X}+1\right) = -\frac{1}{X} + \frac{\frac{1}{4}}{X^2} + O\left(\frac{1}{X^3}\right)$$

$$\mathrm{dilog}(X+1) = -\frac{1}{2}\ln(X)^2 - \frac{1}{6}\pi^2 + \frac{1}{X} - \frac{1}{4}\frac{1}{X^2} + O\left(\frac{1}{X^3}\right)$$

$$\ln(X+1) = \ln(X) + \frac{1}{X} - \frac{1}{2}\frac{1}{X^2} + O\left(\frac{1}{X^3}\right)$$

$$\frac{1}{X+1} = \frac{1}{X} - \frac{1}{X^2} + O\left(\frac{1}{X^3}\right)$$

Inserting these series expansion into the expression from (A8) above gives

$$\varepsilon_\nu(s) = \frac{1}{4}\,\frac{-\frac{1}{2}\ln(X)^2 - \frac{1}{6}\pi^2}{\pi^2} + \frac{\ln(X)\,\nu}{\pi} - \frac{3}{2}\nu^2 + \frac{\frac{1}{2}\frac{1}{\pi^2} + \frac{\nu}{\pi} - \frac{\nu^2}{s-1}}{X} + \frac{-\frac{1}{8}\frac{1}{\pi^2} - \frac{1}{2}\frac{\nu}{\pi} + \frac{\nu^2}{s-1}}{X^2} + O\left(\frac{1}{X^3}\right)$$

Substitute $X$ back to exponential

$$\varepsilon_\nu(s) = \frac{1}{4}\,\frac{-\frac{1}{2}\ln\left(e^{(2\,\nu\,\pi)}\right)^2 - \frac{1}{6}\pi^2}{\pi^2} + \frac{\ln\left(e^{(2\,\nu\,\pi)}\right)\nu}{\pi} - \frac{3}{2}\nu^2 + \frac{\frac{1}{2}\frac{1}{\pi^2} + \frac{\nu}{\pi} - \frac{\nu^2}{s-1}}{e^{(2\,\nu\,\pi)}} + \frac{-\frac{1}{8}\frac{1}{\pi^2} - \frac{1}{2}\frac{\nu}{\pi} + \frac{\nu^2}{s-1}}{\left(e^{(2\,\nu\,\pi)}\right)^2}$$
$$+ O\left(\frac{1}{\left(e^{(2\,\nu\,\pi)}\right)^3}\right)$$

Simplify

$$\varepsilon_\nu(s) = \frac{1}{4}\,\frac{-2\,\nu^2\,\pi^2 - \frac{1}{6}\pi^2}{\pi^2} + \frac{1}{2}\nu^2 + \frac{\frac{1}{2}\frac{1}{\pi^2} + \frac{\nu}{\pi} - \frac{\nu^2}{s-1}}{e^{(2\,\nu\,\pi)}} + \frac{-\frac{1}{8}\frac{1}{\pi^2} - \frac{1}{2}\frac{\nu}{\pi} + \frac{\nu^2}{s-1}}{\left(e^{(2\,\nu\,\pi)}\right)^2} + O\left(\frac{1}{\left(e^{(2\,\nu\,\pi)}\right)^3}\right)$$

Reduce order

$$\varepsilon_\nu(s) = \frac{1}{4}\,\frac{-2\,\nu^2\,\pi^2 - \frac{1}{6}\pi^2}{\pi^2} + \frac{1}{2}\nu^2 + \left(\frac{1}{2}\frac{1}{\pi^2} + \frac{\nu}{\pi} - \frac{\nu^2}{s-1}\right) e^{(-2\,\nu\,\pi)} + O\left(\nu^2\, e^{(-4\,\nu\,\pi)}\right)$$

Simplify the constant terms

$$\varepsilon_\nu(s) = -\frac{1}{24} + \left(\frac{1}{2}\frac{1}{\pi^2} + \frac{\nu}{\pi} - \frac{\nu^2}{s-1}\right) e^{(-2\,\nu\,\pi)} + O\left(\nu^2\, e^{(-4\,\nu\,\pi)}\right)$$

which is thus the final result on page 13.

---

## FAQ #16a (page 13)

**"Of the four series expansions you make in the beginning of FAQ #16 above, the expansions for dilog(1/*X*+1), ln(*X*+1), and 1/(*X*+1) seem straightforward. But how do you get to the series expansion of dilog(*X*+1)**

$$\operatorname{dilog}(X+1) = -\frac{1}{2}\ln(X)^2 - \frac{1}{6}\pi^2 + \frac{1}{X} - \frac{1}{4}\frac{1}{X^2} + O\left(\frac{1}{X^3}\right)$$

**as given above?"**

<u>ANSWER:</u> The dilog function has several interesting properties, especially in the form of sums of the type treated in Theorem 2 and Theorem 3 below.

For the proof of the series expansion above we also need one of the straightforward expansions, namely of dilog(1/*X*+1) as given in Theorem 1 below.

***<u>Theorem 1:</u>***

$$\mathrm{dilog}\left(\frac{1}{X}+1\right)=-\frac{1}{X}+\frac{\frac{1}{4}}{X^2}+\mathrm{O}\left(\frac{1}{X^3}\right) \tag{1}$$

*Proof of Theorem1:*

Make the transformation *X* -> 1/*x* on the left-hand side of (1)

$$\mathrm{dilog}\left(\frac{1}{X}+1\right)=\mathrm{dilog}(x+1) \tag{2}$$

From the definition of dilog(X)

$$\mathrm{dilog}(X)=-\int_1^X \frac{\ln(t)}{t-1}\,dt \tag{3}$$

we then have from (2)

$$\mathrm{dilog}\left(\frac{1}{X}+1\right)=-\int_1^{x+1} \frac{\ln(t)}{t-1}\,dt \tag{4}$$

For *x* = 0 the right-hand side of (4) vanishes (since upper and lower limits then are equal), so the constant term in the expansion of dilog(1/*X* +1) vanishes in agreement with (1).

For *x* = 0 the first derivative of the right-hand side of (4) becomes

$$\lim_{x\to 0}\frac{\partial}{\partial x}\left(-\int_1^{x+1}\frac{\ln(t)}{t-1}\,dt\right) = \lim_{x\to 0} -\frac{\ln(x+1)}{x} = -1 \tag{5}$$

in agreement with the coefficient for 1/*X* in the expansion in (1).

For *x* = 0 the coefficient involving the second derivative in the expansion of the right-hand side of (4) similarly becomes

$$\lim_{x\to 0}\frac{1}{2}\left(\frac{\partial^2}{\partial x^2}\left(-\int_1^{x+1}\frac{\ln(t)}{t-1}\,dt\right)\right) = \lim_{x\to 0} -\frac{1}{2}\frac{1}{(x+1)\,x}+\frac{\frac{1}{2}\ln(x+1)}{x^2} = \frac{1}{4} \tag{6}$$

in agreement with the coefficient for $1/X^2$ in the expansion in (1). Theorem 1 is thus proved.

*Theorem 2:*

$$\text{dilog}(X) + \text{dilog}\left(\frac{1}{X}\right) = -\frac{1}{2}\ln(X)^2 \qquad (7)$$

*Proof of Theorem 2:*

Definition of dilog(*X*):

$$\text{dilog}(X) = -\int_{1}^{X} \frac{\ln(t)}{t-1}\,dt \qquad (8)$$

*Corollary:*

$$\text{dilog}\left(\frac{1}{X}\right) = -\int_{1}^{\frac{1}{X}} \frac{\ln(t)}{t-1}\,dt \qquad (9)$$

Make variable transformation $t$ -> $1/t$ in (8)

$$\text{dilog}(X) = \int_{1}^{\frac{1}{X}} \frac{\ln\left(\frac{1}{t}\right)}{\left(\frac{1}{t}-1\right)t^2}\,dt \qquad (10)$$

i e

$$\text{dilog}(X) = -\int_{1}^{\frac{1}{X}} \frac{\ln\left(\frac{1}{t}\right)}{t\,(t-1)}\,dt \qquad (11)$$

Now form the sum in (7) using (11) and (9)

$$\text{dilog}(X) + \text{dilog}\left(\frac{1}{X}\right) = -\int_{1}^{\frac{1}{X}} \frac{\ln(t)}{t-1}\,dt - \int_{1}^{\frac{1}{X}} \frac{\ln\left(\frac{1}{t}\right)}{t\,(t-1)}\,dt \qquad (12)$$

i e

$$\text{dilog}(X) + \text{dilog}\left(\frac{1}{X}\right) = \int_{1}^{\frac{1}{X}} -\frac{\ln(t)}{t-1} + \frac{\ln(t)}{t\,(t-1)}\,dt \qquad (13)$$

or

$$\text{dilog}(X) + \text{dilog}\left(\frac{1}{X}\right) = -\int_{1}^{\frac{1}{X}} \frac{\ln(t)}{t}\, dt \qquad (14)$$

(14) integrates to

$$\text{dilog}(X) + \text{dilog}\left(\frac{1}{X}\right) = -\frac{1}{2}\ln(X)^2 \qquad (15)$$

which thus proves Theorem 2.

*<u>Theorem 3</u>:*

$$\text{dilog}\left(\frac{1}{x} + 1\right) + \text{dilog}(x + 1) = -\frac{1}{2}\ln(x)^2 - \frac{1}{6}\pi^2 \qquad (16)$$

*Proof of Theorem 3:*

Using the definition of dilog($X$) in (3) we can evaluate the left-hand side of (16) as follows

$$\text{dilog}\left(\frac{1}{x} + 1\right) + \text{dilog}(x + 1) = -\int_{1}^{\frac{1}{x}+1} \frac{\ln(t)}{t-1}\, dt - \int_{1}^{x+1} \frac{\ln(t)}{t-1}\, dt \qquad (17)$$

Now consider the following function

$$\text{dilog}\left(\frac{1}{x} + 1\right) + \text{dilog}(x + 1) + \frac{1}{2}\ln(x)^2 = -\int_{1}^{\frac{1}{x}+1} \frac{\ln(t)}{t-1}\, dt - \int_{1}^{x+1} \frac{\ln(t)}{t-1}\, dt + \frac{1}{2}\ln(x)^2 \qquad (18)$$

and differentiate both sides of (18) with respect to $x$

$$\frac{\partial}{\partial x}\left(\text{dilog}\left(\frac{1}{x} + 1\right) + \text{dilog}(x + 1) + \frac{1}{2}\ln(x)^2\right) = \frac{\partial}{\partial x}\left(-\int_{1}^{\frac{1}{x}+1} \frac{\ln(t)}{t-1}\, dt - \int_{1}^{x+1} \frac{\ln(t)}{t-1}\, dt + \frac{1}{2}\ln(x)^2\right) \qquad (19)$$

which evaluates to

$$\frac{\partial}{\partial x}\left(\text{dilog}\left(\frac{1}{x} + 1\right) + \text{dilog}(x + 1) + \frac{1}{2}\ln(x)^2\right) = \frac{\left(\frac{1}{x} - \frac{x+1}{x^2}\right)\ln\left(\frac{x+1}{x}\right)}{1 - \frac{x+1}{x}} - \frac{\ln(x+1)}{x} + \frac{\ln(x)}{x} \qquad (20)$$

and simplifies to

$$\frac{\partial}{\partial x}\left(\mathrm{dilog}\left(\frac{x+1}{x}\right)+\mathrm{dilog}(x+1)+\frac{1}{2}\ln(x)^2\right)=0 \qquad (21)$$

Since the derivative in (21) thus vanishes identically for all $x$, the function on the left-hand side of (18) must be a constant,

$$\mathrm{dilog}\left(\frac{1}{x}+1\right)+\mathrm{dilog}(x+1)+\frac{1}{2}\ln(x)^2=C \qquad (22)$$

The constant $C$ can be determined by, e g, considering the case $x = 1$, when (22) becomes

$$2\ \mathrm{dilog}(2)=C \qquad (23)$$

or

$$-\frac{1}{6}\pi^2=C \qquad (24)$$

From (22) and (24) we then get

$$\mathrm{dilog}\left(\frac{1}{x}+1\right)+\mathrm{dilog}(x+1)=-\frac{1}{2}\ln(x)^2-\frac{1}{6}\pi^2 \qquad (25)$$

This hence proves Theorem 3, which is thus a relationship of a similar type as that given in Theorem 2 above.

*Proof of FAQ #16a:*

Setting $x = 1/X$ in Theorem 3 gives

$$\mathrm{dilog}(X+1)+\mathrm{dilog}\left(\frac{1}{X}+1\right)=-\frac{1}{2}\ln\left(\frac{1}{X}\right)^2-\frac{1}{6}\pi^2 \qquad (26)$$

$$\mathrm{dilog}(X+1)+\mathrm{dilog}\left(\frac{1}{X}+1\right)=-\frac{1}{2}\ln(X)^2-\frac{1}{6}\pi^2 \qquad (27)$$

Inserting dilog($1/X$+1) from Theorem 1 into (27) then gives the final result

$$\mathrm{dilog}(X+1)=-\frac{1}{2}\ln(X)^2-\frac{1}{6}\pi^2+\frac{1}{X}-\frac{1}{4}\frac{1}{X^2}+\mathrm{O}\left(\frac{1}{X^3}\right) \qquad (28)$$

which thus proves FAQ #16a.

---

## FAQ #17 (page 14)

**"In the middle of page 14 in your preprint you say that you want to estimate the error you make in $\Delta I_2(\pi/2)$ when you neglect the rest of the integral outside the interval $-N\delta \leq \tau \leq N\delta$, and you then give what is obviously only the final result of a calculation. Can you please give the complete calculation of the error."**

<u>ANSWER:</u> Define

$$\mathrm{h}(y, s, N) = \frac{\mathrm{e}^{(i\,s\,y)}}{\mathrm{e}^{(2\,\pi\,N\,(\cos(y) + i\,\sin(y)))} + 1} \qquad \textit{(a1)}$$

and

$$\mathrm{g}(y, s) = \begin{cases} 0 & y < \frac{1}{2}\pi \\ \mathrm{e}^{(i\,s\,y)} & \frac{1}{2}\pi - y \le 0 \text{ and } y - \frac{3}{2}\pi \le 0 \\ 0 & \frac{3}{2}\pi < y \end{cases} \qquad \textit{(a2)}$$

The integrals (A1) and (A2) in Appendix A in the preprint can then be written

$$I_2 = i\,(2\,\pi\,N)^s \int_0^{2\,\pi} \mathrm{h}(y, s, N)\,dy \qquad \textit{(A1)}$$

and

$$I_0 = i\,(2\,\pi\,N)^s \int_0^{2\,\pi} \mathrm{g}(y, s)\,dy \qquad \textit{(A2)}$$

The functions h(*y*, *s*, *N*) and g(*y*, *s*) differ appreciably only in two narrow regions around $y = \pi/2$ and $y = 3\pi/2$, respectively. Below I consider only the case around $y = \pi/2$ (the case around $y = 3\pi/2$ is treated similarly).

We want to estimate the error that I make when instead of integrating the difference between (A1) and (A2) over the whole range around $\pi/2$,

$$\Delta\, I_2\!\left(\frac{1}{2}\pi\right) = i\,(2\,\pi\,N)^s \int_0^{\pi} \mathrm{h}(y, s, N) - \mathrm{g}(y, s)\,dy \qquad \textit{(1)}$$

I only integrate it over the particular, much narrower domain that I use in (A4) in my preprint,

$$\Delta\, I_2\!\left(\frac{1}{2}\pi\right) = i\,(2\,\pi\,N)^s \int_{1/2\,\pi - \delta}^{1/2\,\pi + \delta} \mathrm{h}(y, s, N) - \mathrm{g}(y, s)\,dy \qquad \textit{(2)}$$

The error I make when I use my expression (2) above compared to the correct (1) is thus

$$\Delta\!\left(\Delta\, I_2\!\left(\frac{1}{2}\pi\right)\right) = i\,(2\,\pi\,N)^s \left(\int_0^{\pi} \mathrm{h}(y, s, N) - \mathrm{g}(y, s)\,dy - \int_{1/2\,\pi - \delta}^{1/2\,\pi + \delta} \mathrm{h}(y, s, N) - \mathrm{g}(y, s)\,dy\right) \qquad \textit{(3)}$$

i e

$$\Delta\left(\Delta I_2\left(\frac{1}{2}\pi\right)\right) = i\,(2\,\pi\,N)^s\left(\int_0^{1/2\,\pi-\delta} \mathrm{h}(y,s,N)-\mathrm{g}(y,s)\,dy + \int_{1/2\,\pi+\delta}^{\pi} \mathrm{h}(y,s,N)-\mathrm{g}(y,s)\,dy\right) \quad (4)$$

or since, according to its definition, $g(y,s) = 0$ for $y < \pi/2$,

$$\Delta\left(\Delta I_2\left(\frac{1}{2}\pi\right)\right) = i\,(2\,\pi\,N)^s\left(\int_0^{1/2\,\pi-\delta} \mathrm{h}(y,s,N)\,dy + \int_{1/2\,\pi+\delta}^{\pi} \mathrm{h}(y,s,N)-\mathrm{g}(y,s)\,dy\right) \quad (5)$$

Inserting (a1) from above, the first integrand in (5) above becomes

$$\mathrm{h}(y,s,N) = \frac{\mathrm{e}^{(i\,s\,y)}}{\mathrm{e}^{(2\,\pi\,N\,(\cos(y)+i\,\sin(y)))}+1} \quad (6)$$

For $y < \pi/2$ the exponential in the denominator becomes very large for large $N$. I can rewrite (6) as follows

$$\mathrm{h}(y,s,N) = \frac{\mathrm{e}^{(i\,s\,y)}\,\mathrm{e}^{(-2\,\pi\,N\,(\cos(y)+i\,\sin(y)))}}{1+\mathrm{e}^{(-2\,\pi\,N\,(\cos(y)+i\,\sin(y)))}} \quad (7)$$

or after series expansion of the denominator (where the exponential now is << 1),

$$\mathrm{h}(y,s,N) = \mathrm{e}^{(i\,s\,y)}\,\mathrm{e}^{(-2\,\pi\,N\,(\cos(y)+i\,\sin(y)))}\,(1-\mathrm{e}^{(-2\,\pi\,N\,(\cos(y)+i\,\sin(y)))}+\,...) \quad (8)$$

$$= \mathrm{e}^{(i\,s\,y)}\,(\mathrm{e}^{(-2\,\pi\,N\,(\cos(y)+i\,\sin(y)))}-\mathrm{e}^{(-4\,\pi\,N\,(\cos(y)+i\,\sin(y)))}+\,...) \quad (9)$$

$$= \mathrm{e}^{(i\,s\,y)}\,(\mathrm{e}^{(-2\,\pi\,N\,(\cos(y)+i\,\sin(y)))}+\mathrm{O}(\mathrm{e}^{(-4\,\pi\,N\,\cos(y))})) \quad (10)$$

I will later show that the correction to the integral $I_2$ corresponding to the first exponential term above is negligible in the present context. Since the remainder in (10) above is the square of this correction, it is even smaller and can thus be neglected in the following calculations. The first integrand in (5) can thus be written

$$\mathrm{h}(y,s,N) = \mathrm{e}^{(-i\,s\,y)}\,\mathrm{e}^{(-2\,\pi\,N\,(\cos(y)+i\,\sin(y)))} \quad (11)$$

Expressed in (a1) and (a2) above, the second integrand in (5) above similarly becomes

$$\mathrm{h}(y,s,N)-\mathrm{g}(y,s) = \frac{\mathrm{e}^{(i\,s\,y)}}{\mathrm{e}^{(2\,\pi\,N\,(\cos(y)+i\,\sin(y)))}+1}-\mathrm{e}^{(i\,s\,y)} \quad (12)$$

$$= -\frac{e^{(i\,s\,y)}\,e^{(2\,\pi\,N\,(\cos(y)+i\,\sin(y)))}}{e^{(2\,\pi\,N\,(\cos(y)+i\,\sin(y)))}+1} \qquad (13)$$

$$= -\frac{e^{(i\,s\,y)}}{1+e^{(-2\,\pi\,N\,(\cos(y)+i\,\sin(y)))}} \qquad (14)$$

Here $y > \pi/2$ and thus $cos(y) < 0$, so again the exponential in the denominator becomes very large for large $N$. As above I can rewrite (14) as follows [or revert to (13)]

$$= -\frac{e^{(i\,s\,y)}\,e^{(2\,\pi\,N\,(\cos(y)+i\,\sin(y)))}}{e^{(2\,\pi\,N\,(\cos(y)+i\,\sin(y)))}+1} \qquad (15)$$

where the exponential in the denominator now is << 1 for $N$ >> 1. After series expansion of the denominator, I thus get

$$= -e^{(i\,s\,y)}\,e^{(2\,\pi\,N\,(\cos(y)+i\,\sin(y)))}\,(1 - e^{(2\,\pi\,N\,(\cos(y)+i\,\sin(y)))} + \ldots) \qquad (16)$$

$$= -e^{(i\,s\,y)}\,(e^{(2\,\pi\,N\,(\cos(y)+i\,\sin(y)))} - e^{(4\,\pi\,N\,(\cos(y)+i\,\sin(y)))} + \ldots) \qquad (17)$$

$$= -e^{(i\,s\,y)}\,(e^{(2\,\pi\,N\,(\cos(y)+i\,\sin(y)))} + O(e^{(4\,\pi\,N\,\cos(y))})) \qquad (18)$$

Again I will later show that the correction to the integral $I_2$ corresponding to the first exponential term above is negligible in the present context. Since the remainder in (18) above is the square of this correction, it is even smaller and can thus in this case too be neglected in the following calculations. The second integrand in (5) can thus be written

$$h(y, s, N) - g(y, s) = -e^{(i\,s\,y)}\,e^{(2\,\pi\,N\,(\cos(y)+i\,\sin(y)))} \qquad (19)$$

We now want to estimate the magnitude of the error I make when I use my integral (2) instead of the correct integral (1). Inserting the results from (11) and (19), we thus now calculate the absolute value of the error in (5) above

$$\left|\Delta\left(\Delta I_2\left(\frac{1}{2}\pi\right)\right)\right| =$$

$$(2\,\pi\,N)^{\sigma}\left|\int_0^{1/2\,\pi-\delta} e^{(-i\,s\,y)}\,e^{(-2\,\pi\,N(\cos(y)+i\,\sin(y)))}\,dy - \int_{1/2\,\pi+\delta}^{\pi} e^{(i\,s\,y)}\,e^{(2\,\pi\,N(\cos(y)+i\,\sin(y)))}\,dy\right| \qquad (20)$$

The first exponential in the two integrals in (20) will have its maximum value when $Im(s)$ is negative. Since $y$ is at most $2\pi$, (20) can thus be estimated as follows

$$\left|\Delta\left(\Delta I_2\left(\frac{1}{2}\pi\right)\right)\right| < (2\pi N)^{\sigma} e^{(2\pi|\Im(s)|)}\left(\left|\int_0^{1/2\pi-\delta} e^{(-2\pi N\cos(y))}\,dy\right| + \left|\int_{1/2\pi+\delta}^{\pi} e^{(2\pi N\cos(y))}\,dy\right|\right) \quad (21)$$

$$< (2\pi N)^{\sigma} e^{(2\pi|\Im(s)|)}\left(\left|\int_0^{1/2\pi-\delta} e^{\left(-2\pi N\left(1-\frac{2y}{\pi}\right)\right)}\,dy\right| + \left|\int_{1/2\pi+\delta}^{\pi} e^{\left(2\pi N\left(1-\frac{2y}{\pi}\right)\right)}\,dy\right|\right) \quad (22)$$

where the last inequality uses a linear approximation through (0,1) and (π,-1) of the cosine, which underestimates the function and thus overestimates the contribution from the negative exponentials, and thus overestimates also the error calculated on the right-hand side of (26) below.

Evaluating (22) we get (for $\delta << 1$)

$$\left|\Delta\left(\Delta I_2\left(\frac{1}{2}\pi\right)\right)\right| < \frac{1}{2}\frac{(2\pi N)^{\sigma} e^{(2\pi|\Im(s)|)}\left|e^{(-4N\delta)} - e^{(-2\pi N)}\right|}{N} < \frac{1}{2}\frac{(2\pi N)^{\sigma} e^{(2\pi|\Im(s)|)} e^{(-4N\delta)}}{N} \quad (23)$$

Setting (for $N$>1)

$$N\delta = \ln(N) \quad (24)$$

(23) gives

$$\left|\Delta\left(\Delta I_2\left(\frac{1}{2}\pi\right)\right)\right| < \frac{1}{2}\frac{(2\pi N)^{\sigma} e^{(2\pi|\Im(s)|)}}{N^5} \quad (25)$$

i e

$$\Delta\left(\Delta I_2\left(\frac{1}{2}\pi\right)\right) < O\left(\frac{1}{N^4}\right) \quad (26)$$

With the assignment of $N\delta = \nu = \ln(N)$ as in (24) above [(A4) in the preprint], the error introduced when I perform the integration in (A5) in the preprint only over the limited domain $-N\delta < \tau < N\delta$, instead of over the complete domain, thus lies well within the remainder in (A5).

---

## FAQ #18 (page 14)

**"In all your calculations on the integrals in Sect A3 in your preprint you have forgotten to discuss the factor $N^s$ in front of the integrals. This factor tends to infinity with $N$ for values of s inside the critical strip. Shouldn't this factor also be taken into account?"**

ANSWER: The term with this factor $N^s$ is part of the pair discussed in Remark 5.2 in the preprint, and is necessary to cancel the corresponding divergence in $S_N$ in (6) according to Cauchy's theorem. So it makes sense to keep this factor outside the calculations of the integrals themselves as is done in, e g, FAQ #17 above.

---

## FAQ #19 (page 14)

**"I'm not sure I can repeat the calculation of $\Delta I_2$ in Sect A3.2 on page 14 in your preprint for $3\pi/2$ correctly. Can you give it?"**

ANSWER: Here are the corresponding steps for $3\pi/2$:

Around $3\pi/2$ the expression corresponding to the one on bottom of page 11 is

$$I_2\left(\frac{3}{2}\pi, \delta\right) = i\, e^{(3/2\, i\, s\, \pi)}\, 2^s\, \pi^s\, N^s \int_{3/2\,\pi-\delta}^{3/2\,\pi+\delta} \frac{1 + i\, s\left(y - \frac{3}{2}\pi\right) - \frac{1}{2} s^2 \left(y - \frac{3}{2}\pi\right)^2 + \ldots}{e^{(2\,\pi\, N\,(-i + y - 3/2\,\pi + 1/2\, i\,(y - 3/2\,\pi)^2 - 1/6\,(y - 3/2\,\pi)^3 + \ldots))} + 1}\, dy$$

Changing integration variable

$$N\left(y - \frac{3}{2}\pi\right) = \tau$$

we get

$$I_2\left(\frac{3}{2}\pi, \delta\right) = i\, e^{(3/2\, i\, s\, \pi)}\, 2^s\, \pi^s\, N^s \int_{-N\delta}^{N\delta} \frac{\frac{1}{N} + \frac{i\, s\, \tau}{N^2} + O\left(\frac{1}{N^3}\right)}{e^{\left(\frac{i\,\pi\,\tau^2}{N} + 2\,\pi\,\tau\right)} + 1}\, d\tau$$

Form the correction, i e the difference between the above integral and corresponding integral over the approximating piecewise function,

$$\Delta I_2\left(\frac{3}{2}\pi\right) = i\, e^{(3/2\, i\, s\, \pi)}\, 2^s\, \pi^s\, N^s \left( \int_{-N\delta}^{N\delta} \frac{\frac{1}{N} + \frac{i\, s\, \tau}{N^2} + O\left(\frac{1}{N^3}\right)}{e^{\left(\frac{i\,\pi\,\tau^2}{N} + 2\,\pi\,\tau\right)} + 1}\, d\tau - \int_{-N\delta}^{0} \frac{1}{N} + \frac{i\, s\, \tau}{N^2} + O\left(\frac{1}{N^3}\right) d\tau \right)$$

Expand the first integrand as a series in $1/N$

$$\Delta I_2\left(\frac{3}{2}\pi\right) = i\, e^{(3/2\, i\, s\, \pi)}\, 2^s\, \pi^s\, N^s \left( \int_{-N\delta}^{N\delta} \frac{a_1}{N} + \frac{a_2}{N^2} + O\left(\frac{1}{N^3}\right) d\tau - \int_{-N\delta}^{0} \frac{1}{N} + \frac{i\, s\, \tau}{N^2} + O\left(\frac{1}{N^3}\right) d\tau \right)$$

where

$$a_1 = \frac{1}{e^{(2\pi\tau)} + 1}$$

$$a_2 = \frac{i\, s\, \tau}{e^{(2\pi\tau)} + 1} - \frac{i\, e^{(2\pi\tau)}\, \pi\, \tau^2}{(e^{(2\pi\tau)} + 1)^2}$$

Integrate for each power of $N$

$$\int_{-N\delta}^{0} \frac{a_1}{N} - \frac{1}{N}\, d\tau + \int_{0}^{N\delta} \frac{a_1}{N}\, d\tau = 0$$

$$\int_{-N\delta}^{0} \frac{a_2}{N^2} - \frac{i\, s\, \tau}{N^2}\, d\tau + \int_{0}^{N\delta} \frac{a_2}{N^2}\, d\tau = \frac{-i\, \varepsilon_\nu(s)\,(s-1)}{N^2}$$

where

$$\varepsilon_\nu(s) = \frac{1}{4}\, \frac{\mathrm{Li}(e^{(2\nu\pi)} + 1) - \mathrm{Li}(e^{(-2\nu\pi)} + 1)}{\pi^2} + \frac{\ln(e^{(2\nu\pi)} + 1)\, \nu}{\pi} - \frac{3}{2}\nu^2 + \frac{\nu^2}{(e^{(2\nu\pi)} + 1)(1 - s)}$$

which thus calculates to the same result as for $\pi/2$, and thus again simplifies to

$$\varepsilon_\nu(s) = -\frac{1}{24} + O(\nu^2\, e^{(-2\nu\pi)})$$

Inserting the integrated results above into the integral, we thus get the final expression on the bottom of page 14,

$$\Delta I_2\left(\frac{3}{2}\pi\right) = i\, e^{(3/2\, i\, s\, \pi)}\, 2^s\, \pi^s\, N^s \left( -\frac{i\, \varepsilon_\nu(s)\,(s-1)}{N^2} + O\left(\frac{\ln(N)}{N^3}\right) \right)$$

## FAQ #20 (pages 17-18)

**"As I understand it, on pages 17-18 in the preprint you want to calculate differences of type $\zeta_N(s)$ - $\zeta(s)$ between some approximation of the zeta-function and the zeta-function itself, because then you can form quotients from which you can extract information about zeros of the zeta-function. But why do you need to choose expressions as in (13) and (14) in the preprint, in which you have vanishing remainders? Wouldn't any remainder within $O(1/N^{3-\sigma})$ or $O(1/N^{2+\sigma})$, respectively, around $\zeta(s)$ do, if you only avoid the very particular remainders that give the exact zeta-function?"**

ANSWER: Yes, in principle you are right, any approximate function within the proper remainder around $\zeta(s)$ would do, but the two functions in (13) and (14) with vanishing remainders have one important property, namely that for them the series expansions on top of page 17 get the simple form given there.

For any other choice of remainders we would have to calculate what form they would get in the series expansions. In principle we can use any approximate function within $O(1/N^{3-\sigma})$ or $O(1/N^{2+\sigma})$, respectively, around $\zeta(s)$ but the particular choice of functions in (13) and (14) with vanishing remainders makes the proof much simpler.

---

## FAQ #20a (page 17)

**"I cannot easily see how you get the series expansions on top of page 17 from the expressions on page 16. Can you show me?"**

ANSWER: In the derivation of the expressions on top on page 17 in my preprint from the expressions in the middle of page 16, I use two alternative approaches (which naturally give the same result). For illustration, I below use one of the approaches to calculate the first expression, and the other approach to calculate the second expression.

FIRST EXPRESSION (PRIME)

EqB1 below is equivalent to the equation immediately after Remark B.1 in the preprint; it differs from that equation only in that a factor of 2 in numerator and denominator are not yet cancelled.

$$EqB1 := \zeta_N(s)' - \zeta_{N+1}(s)' =$$

$$-\frac{1}{2}\,\frac{\pi^s\left(2^s N^s\left(1-\frac{1}{24}\frac{s(s-1)}{N^2}\right)-2^s(N+1)^s\left(1-\frac{1}{24}\frac{s(s-1)}{(N+1)^2}\right)+2\,s\,(1+2\,N)^{(s-1)}\right)}{\cos\left(\frac{1}{2}s\,\pi\right)(-1+2^{(1-s)})\,\Gamma(s+1)}$$

In TMP1 below some substitutions have for simplicity been made in EqB1

$$TMP2 := (N+1)^s\left(1-\frac{1}{24}\frac{s(s-1)}{(N+1)^2}\right) = tmp2$$

$$TMP3 := (1+2N)^{(s-1)} = tmp3$$

$$TMP1 := \zeta_N(s)' - \zeta_{N+1}(s)' = -\frac{1}{2}\frac{\pi^s\left(2^s N^s\left(1-\frac{1}{24}\frac{s(s-1)}{N^2}\right) - 2^s\,tmp2 + 2\,s\,tmp3\right)}{\cos\left(\frac{1}{2}s\pi\right)(-1+2^{(1-s)})\,\Gamma(s+1)}$$

Rewrite TMP2 above as follows, then expand it as a Taylor series in $1/N$

$$N^s\left(1+\frac{1}{N}\right)^s\left(1-\frac{1}{24}\frac{s(s-1)}{(N+1)^2}\right)$$

$$TMP4 := tmp2 = N^s\left(1+\frac{s}{N}+\frac{\frac{11}{24}s(s-1)}{N^2}+\frac{\frac{1}{12}s(s-1)-\frac{1}{24}s^2(s-1)+\frac{1}{6}s(s-1)(s-2)}{N^3}\right.$$
$$+\frac{\frac{1}{12}s^2(s-1)+\frac{1}{24}s(s-1)(s-2)(s-3)-\frac{1}{48}s^2(s-1)^2-\frac{1}{8}s(s-1)}{N^4}+\left(-\frac{1}{8}s^2(s-1)\right.$$
$$\left.+\frac{1}{6}s(s-1)+\frac{1}{24}s^2(s-1)^2-\frac{1}{144}s^2(s-1)^2(s-2)+\frac{1}{120}s(s-1)(s-2)(s-3)(s-4)\right)\Big/N^5$$
$$\left.+\mathrm{O}\left(\frac{1}{N^6}\right)\right)$$

Similarly rewrite TMP3 above as follows, then expand it as a Taylor series in $1/N$, then simplify the factors

$$N^{(s-1)}\left(\frac{1}{N}+2\right)^{(s-1)}$$

$$tmp3 = N^s\left(\frac{2^{(s-1)}}{N}+\frac{2^{(s-1)}\left(\frac{1}{2}s-\frac{1}{2}\right)}{N^2}+\frac{\frac{1}{4}2^{(s-1)}\left(\frac{1}{2}s-\frac{1}{2}\right)(s-2)}{N^3}+\frac{\frac{1}{24}2^{(s-1)}\left(\frac{1}{2}s-\frac{1}{2}\right)(s-2)(s-3)}{N^4}\right.$$
$$\left.+\frac{\frac{1}{192}2^{(s-1)}\left(\frac{1}{2}s-\frac{1}{2}\right)(s-2)(s-3)(s-4)}{N^5}+\mathrm{O}\left(\frac{1}{N^6}\right)\right)$$

$$TMP5 := tmp3 = N^s\left(\frac{2^{(s-1)}}{N} + \frac{\frac{1}{2}2^{(s-1)}(s-1)}{N^2} + \frac{\frac{1}{8}2^{(s-1)}(s-1)(s-2)}{N^3} + \frac{\frac{1}{48}2^{(s-1)}(s-1)(s-2)(s-3)}{N^4} + \frac{\frac{1}{384}2^{(s-1)}(s-1)(s-2)(s-3)(s-4)}{N^5} + O\left(\frac{1}{N^6}\right)\right)$$

Now reinsert into TMP1 the expressions for tmp2 and tmp3 we have just calculated in TMP4 and TMP5

$$EqB2 := \zeta_N(s)' - \zeta_{N+1}(s)' = -\frac{1}{2}\pi^s\left(2^s N^s\left(1 - \frac{1}{24}\frac{s(s-1)}{N^2}\right) - 2^s N^s\left(1 + \frac{s}{N} + \frac{\frac{11}{24}s(s-1)}{N^2} + \frac{\frac{1}{12}s(s-1) - \frac{1}{24}s^2(s-1) + \frac{1}{6}s(s-1)(s-2)}{N^3} + \frac{\frac{1}{12}s^2(s-1) + \frac{1}{24}s(s-1)(s-2)(s-3) - \frac{1}{48}s^2(s-1)^2 - \frac{1}{8}s(s-1)}{N^4} + \left(-\frac{1}{8}s^2(s-1) + \frac{1}{6}s(s-1) + \frac{1}{24}s^2(s-1)^2 - \frac{1}{144}s^2(s-1)^2(s-2) + \frac{1}{120}s(s-1)(s-2)(s-3)(s-4)\right)\Big/ N^5 + O\left(\frac{1}{N^6}\right)\right) + 2\,s\,N^s\left(\frac{2^{(s-1)}}{N} + \frac{\frac{1}{2}2^{(s-1)}(s-1)}{N^2} + \frac{\frac{1}{8}2^{(s-1)}(s-1)(s-2)}{N^3} + \frac{\frac{1}{48}2^{(s-1)}(s-1)(s-2)(s-3)}{N^4} + \frac{\frac{1}{384}2^{(s-1)}(s-1)(s-2)(s-3)(s-4)}{N^5} + O\left(\frac{1}{N^6}\right)\right)\right)\Big/\left(\cos\left(\frac{1}{2}s\pi\right)(-1 + 2^{(1-s)})\Gamma(s+1)\right)$$

Simplify the above equation, factorise it (without remainders), and simplify again. This gives the final result as given in the preprint.

$$\zeta_N(s)' - \zeta_{N+1}(s)' = \frac{1}{11520}\pi^s\left(11520\,2^s N^{(s+5)}\,s\,O\left(\frac{1}{N^6}\right) - 5760\,4^s N^{(s+5)}\,O\left(\frac{1}{N^6}\right) + 7\,4^s N^s s^5 - 70\,4^s N^s s^4 + 245\,4^s N^s s^3 - 350\,4^s N^s s^2 + 168\,4^s N^s s\right)\Big/\left(N^5\cos\left(\frac{1}{2}s\pi\right)\Gamma(s+1)(2^s - 2)\right)$$

$$\zeta_N(s)' - \zeta_{N+1}(s)' = \frac{7}{11520}\,\frac{\pi^s\,4^s\,N^s\,s\,(s-1)(s-2)(s-3)(s-4)}{N^5\cos\left(\frac{1}{2}s\pi\right)\Gamma(s+1)(2^s-2)}$$

SECOND EXPRESSION (DOUBLE PRIME)

EqB3 below is equivalent to the second equation immediately after Remark B.1 in the preprint; it differs from that equation only in that a factor of 2 in numerator and denominator are not yet cancelled.

$$EqB3 := \zeta_N(s)'' - \zeta_{N+1}(s)'' =$$
$$\frac{1}{2}\,\frac{-2\,N^{(1-s)}\left(1-\frac{1}{24}\frac{s\,(s-1)}{N^2}\right)+2\,(N+1)^{(1-s)}\left(1-\frac{1}{24}\frac{s\,(s-1)}{(N+1)^2}\right)-2^{(s+1)}\,(1-s)\,(1+2\,N)^{(-s)}}{(1-s)\,(-1+2^s)}$$

Define tmp6 as follows

$$TMP6 := tmp6 = 1-\frac{1}{24}\frac{s\,(s-1)}{N^2}-\left(1+\frac{1}{N}\right)^{(1-s)}\left(1-\frac{1}{24}\frac{s\,(s-1)}{(N+1)^2}\right)+\frac{2^s\,(1-s)\left(\frac{1}{N}+2\right)^{(-s)}}{N}$$

in which case we can write EqB3 as follows

$$EqB4 := \zeta_N(s)'' - \zeta_{N+1}(s)'' = -\frac{N^{(1-s)}\,tmp6}{(1-s)\,(-1+2^s)}$$

Rewrite TMP6 with $N = 1/n$, then expand expressions, calculate lead term in series expansion, factorise lead term, and then finally insert tmp6 into EqB4 above, which gives the final result as given in the preprint

$$tmp6 = 1-\frac{1}{24}\,s\,(s-1)\,n^2-(1+n)^{(1-s)}\left(1-\frac{1}{24}\frac{s\,(s-1)}{\left(\frac{1}{n}+1\right)^2}\right)+2^s\,(1-s)\,(n+2)^{(-s)}\,n$$

$$tmp6 = 1-\frac{1}{24}\,s^2\,n^2+\frac{1}{24}\,s\,n^2-(1+n)^{(-1-s)}-2\,(1+n)^{(-1-s)}\,n-(1+n)^{(-1-s)}\,n^2$$
$$+\frac{1}{24}\,(1+n)^{(-1-s)}\,s^2\,n^2-\frac{1}{24}\,(1+n)^{(-1-s)}\,s\,n^2+n\left(\frac{1}{2}\,n+1\right)^{(-s)}-n\left(\frac{1}{2}\,n+1\right)^{(-s)}\,s$$

$$tmp6 = \left(\frac{7}{960}\,s+\frac{7}{1152}\,s^2-\frac{7}{1152}\,s^3-\frac{7}{1152}\,s^4-\frac{7}{5760}\,s^5\right)n^5$$

$$tmp6 = -\frac{7}{5760} s\,(s-1)\,(s+3)\,(s+2)\,(s+1)\,n^5$$

$$\zeta_N(s)'' - \zeta_{N+1}(s)'' = -\frac{7}{5760}\frac{N^{(-4-s)}\, s\,(s+3)\,(s+2)\,(s+1)}{-1+2^s}$$

---

## FAQ #21 (page 18)

**"I have difficulties understanding how you get to the second quotient involving *N*+ *k* in *LHS* on page 18 in your preprint from the first quotient involving *N*+ 1. Can you show me the detailed calculation?"**

ANSWER: It is simplest to start from the following expressions on page 17 in the preprint, which I here call (1a) and (1b),

$$\frac{\zeta_N(s)' - \zeta_{N+1}(s)'}{N^{(s-4)} - (N+1)^{(s-4)}} = -\frac{7}{11520}\frac{\pi^s\, 4^s}{\cos\left(\frac{1}{2} s\,\pi\right)(-2+2^s)\,\Gamma(s-3)} + \mathrm{O}\left(\frac{1}{N}\right) \quad (1a)$$

$$\frac{\zeta_N(s)'' - \zeta_{N+1}(s)''}{N^{(-3-s)} - (N+1)^{(-3-s)}} = -\frac{7}{5760}\frac{(s+2)\,(s+1)\,s}{-1+2^s} + \mathrm{O}\left(\frac{1}{N}\right) \quad (1b)$$

Setting

$$C' = -\frac{7}{11520}\frac{\pi^s\, 4^s}{\cos\left(\frac{1}{2} s\,\pi\right)(-2+2^s)\,\Gamma(s-3)} + \mathrm{O}\left(\frac{1}{N}\right) \quad (2a)$$

$$C'' = -\frac{7}{5760}\frac{(s+2)\,(s+1)\,s}{-1+2^s} + \mathrm{O}\left(\frac{1}{N}\right) \quad (2b)$$

(1a) and (1b) can be written as

$$\frac{\zeta_N(s)' - \zeta_{N+1}(s)'}{N^{(s-4)} - (N+1)^{(s-4)}} = C' \quad (3a)$$

$$\frac{\zeta_N(s)'' - \zeta_{N+1}(s)''}{N^{(-3-s)} - (N+1)^{(-3-s)}} = C'' \quad (3b)$$

or

$$\zeta_N(s)' - \zeta_{N+1}(s)' = C'\,(N^{(s-4)} - (N+1)^{(s-4)}) \quad (4a)$$

$$\zeta_N(s)'' - \zeta_{N+1}(s)'' = C''\,(N^{(-3-s)} - (N+1)^{(-3-s)}) \quad (4b)$$

Now use (4a) to form a sequence of equations by successively substituting $N$ -> $N + 1$, $N$ -> $N + 2$, $N$ -> $N + 3$, $N$ -> $N + 4$, … $N$ -> $N + k - 2$, $N$ -> $N + k - 1$, and then summing them, whereby all terms cancel except the first and the last ones on the left- and right-hand sides, giving the result in (5a) below [within O(1/$N$)],

$$\zeta_N(s)' - \zeta_{N+1}(s)' = C'\,(N^{(s-4)} - (N+1)^{(s-4)})$$

$$\zeta_{N+1}(s)' - \zeta_{N+2}(s)' = C'\,((N+1)^{(s-4)} - (N+2)^{(s-4)})$$

$$\zeta_{N+2}(s)' - \zeta_{N+3}(s)' = C'\,((N+2)^{(s-4)} - (N+3)^{(s-4)})$$

$$\zeta_{N+3}(s)' - \zeta_{N+4}(s)' = C'\,((N+3)^{(s-4)} - (N+4)^{(s-4)})$$

.
.
.

$$\zeta_{N+k-2}(s)' - \zeta_{N+k-1}(s)' = C'((N+k-2)^{(s-4)} - (N+k-1)^{(s-4)})$$

$$\zeta_{N+k-1}(s)' - \zeta_{N+k}(s)' = C'\,((N+k-1)^{(s-4)} - (N+k)^{(s-4)})$$

---

$$\zeta_N(s)' - \zeta_{N+k}(s)' = C'\,(N^{(s-4)} - (N+k)^{(s-4)}) \qquad (5a)$$

Similarly use (4b) to form a sequence of equations by again successively substituting $N$ -> $N + 1$, $N$ -> $N + 2$, $N$ -> $N + 3$, $N$ -> $N + 4$, … $N$ -> $N + k - 2$, $N$ -> $N + k - 1$, and then summing them, whereby again all terms cancel except the first and the last ones on both sides, giving the result in (5b) below [within O(1/$N$)],

$$\zeta_N(s)'' - \zeta_{N+1}(s)'' = C''\,(N^{(-3-s)} - (N+1)^{(-3-s)})$$

$$\zeta_{N+1}(s)'' - \zeta_{N+2}(s)'' = C''\,((N+1)^{(-3-s)} - (N+2)^{(-3-s)})$$

$$\zeta_{N+2}(s)'' - \zeta_{N+3}(s)'' = C''\,((N+2)^{(-3-s)} - (N+3)^{(-3-s)})$$

$$\zeta_{N+3}(s)'' - \zeta_{N+4}(s)'' = C''\,((N+3)^{(-3-s)} - (N+4)^{(-3-s)})$$

.
.
.

$$\zeta_{N+k-2}(s)'' - \zeta_{N+k-1}(s)'' = C''((N+k-2)^{(-3-s)} - (N+k-1)^{(-3-s)})$$

$$\zeta_{N+k-1}(s)'' - \zeta_{N+k}(s)'' = C''\,((N+k-1)^{(-3-s)} - (N+k)^{(-3-s)})$$

---

$$\zeta_N(s)'' - \zeta_{N+k}(s)'' = C''\,(N^{(-3-s)} - (N+k)^{(-3-s)}) \qquad (5b)$$

Now solve *C'* and *C"* from (5a) and (5b), respectively,

$$C' = \frac{\zeta_N(s)' - \zeta_{N+k}(s)'}{N^{(s-4)} - (N+k)^{(s-4)}} \qquad (6a)$$

$$C'' = \frac{\zeta_N(s)'' - \zeta_{N+k}(s)''}{N^{(-3-s)} - (N+k)^{(-3-s)}} \qquad (6b)$$

Combining (3a) with (6a), and (3b) with (6b), and also evaluating (6a) and (6b) in the limit $k$ -> ∞, we then get within O(1/*N*),

$$\frac{\zeta_N(s)' - \zeta_{N+1}(s)'}{N^{(s-4)} - (N+1)^{(s-4)}} = \frac{\zeta_N(s)' - \zeta_{N+k}(s)'}{N^{(s-4)} - (N+k)^{(s-4)}} = \frac{\zeta_N(s)' - \zeta(s)}{N^{(s-4)}} \qquad (7a)$$

$$\frac{\zeta_N(s)'' - \zeta_{N+1}(s)''}{N^{(-3-s)} - (N+1)^{(-3-s)}} = \frac{\zeta_N(s)'' - \zeta_{N+k}(s)''}{N^{(-3-s)} - (N+k)^{(-3-s)}} = \frac{\zeta_N(s)'' - \zeta(s)}{N^{(-3-s)}} \qquad (7b)$$

which thus agree with, respectively, the numerators and the denominators in *LHS* on page 18 in the preprint.

---

## FAQ #21a (page 18)

**"In your derivations in FAQ #21 above, you treat the quantities *C'* and *C"* as constants. But this is correct only up to the last step when you calculate the rightmost equality in (7a) and (7b). According to equations (2a) and (2b) above, *C'* and *C"* are actually functions of *N* due to the remainders O(1/*N*), although these do vanish if *N* tends to infinity. However, in the rightmost equalities in (7a) and (7b) you are studying finite *N*, and you add *k* terms containing these remainders, and then let *k* tend to infinity in order to get the third equalities in (7a) and (7b). But then the remainder contributions obviously add up to a sum of the order of *k* times 1/*N*, which might tend to infinity with *k*."**

ANSWER: No actually they don't. One has to do the summation you describe carefully and take into account that there are indeed $k$ terms as you say, but also that their sum is a sum over $n$ of remainders smaller than of type $1/(N+n)^5$ as I will show below. There is then no problem when we sum them over $n$ with $n$ going from 0 to $k$, even if we let $k$ tend to infinity.

First we need to calculate the remainders explicitly. As derived above in FAQ #20a, the leading terms $T_{5,N}'$ and $T_{5,N}''$ in the expansions of $\zeta_N(s)' - \zeta_{N+1}(s)'$ and $\zeta_N(s)'' - \zeta_{N+1}(s)''$ are given on top of page 17 in the preprint. The remainder is equal to the next term in each expansion, and can be obtained by calculating the leading terms in, respectively, $\zeta_N(s)' - \zeta_{N+1}(s)' - T_{5,N}'$ and $\zeta_N(s)'' - \zeta_{N+1}(s)'' - T_{5,N}''$. We then get the two remainders in explicit form as follows

$$\zeta_N(s)' - \zeta_{N+1}(s)' - T_{5,N}' = -\frac{7}{23040}\frac{\pi^s N^{(s-6)} 2^s}{\cos\left(\frac{1}{2}s\pi\right)(-1+2^{(1-s)})\Gamma(s-5)} \qquad (1)$$

$$\zeta_N(s)'' - \zeta_{N+1}(s)'' - T_{5,N}'' = \frac{7}{11520}\frac{N^{(-5-s)} s(s+4)(s+3)(s+2)(s+1)}{-1+2^s} \qquad (2)$$

Thus the remainder terms O(1/*N*) on the right-hand sides of (1a), (2a), and (1b), (2b) in FAQ #21 become, respectively,

$$\mathrm{O}\left(\frac{1}{N}\right)' = \frac{c'}{N} \qquad (3)$$

$$\mathrm{O}\left(\frac{1}{N}\right)'' = \frac{c''}{N} \qquad (4)$$

where

$$c' = -\frac{7}{23040}\frac{\pi^s 2^s}{\cos\left(\frac{1}{2}s\pi\right)(-1+2^{(1-s)})\Gamma(s-5)} \qquad (5)$$

$$c'' = \frac{7}{11520}\frac{s(s+4)(s+3)(s+2)(s+1)}{-1+2^s} \qquad (6)$$

The "constants" $C'$ and $C''$ in FAQ #21 should thus rightly be written as

$$C' = C_0' + \frac{c'}{N} \qquad (7)$$

$$C'' = C_0'' + \frac{c''}{N} \qquad (8)$$

where $C_0'$ and $C_0''$ are true constants, for which the cancelling of terms in the schemes leading to (5a) and (5b) in FAQ #21 is exactly true.

Inserting the expression in (7) for $C'$ in the scheme leading to (5a) in FAQ #21, we then have for the equations around index $N+n$,

$$\zeta_{N+n-2}(s)' - \zeta_{N+n-1}(s)' = \left(C_0' + \frac{c'}{N+n-2}\right)((N+n-2)^{(s-4)} - (N+n-1)^{(s-4)})$$

$$\zeta_{N+n-1}(s)' - \zeta_{N+n}(s)' = \left(C_0' + \frac{c'}{N+n-1}\right)((N+n-1)^{(s-4)} - (N+n)^{(s-4)})$$

$$\zeta_{N+n}(s)' - \zeta_{N+n+1}(s)' = \left(C_0' + \frac{c'}{N+n}\right)((N+n)^{(s-4)} - (N+n+1)^{(s-4)})$$

$$\zeta_{N+n+1}(s)' - \zeta_{N+n+2}(s)' = \left(C_0' + \frac{c'}{N+n+1}\right)((N+n+1)^{(s-4)} - (N+n+2)^{(s-4)})$$

and similarly by inserting the expression in (8) for $C''$ in the scheme leading to (5b) in FAQ #21, we have for the equations around index $N+n$ in that case,

$$\zeta_{N+n-2}(s)'' - \zeta_{N+n-1}(s)'' = \left( C_0'' + \frac{c''}{N+n-2} \right) ((N+n-2)^{(-s-3)} - (N+n-1)^{(-s-3)})$$

$$\zeta_{N+n-1}(s)'' - \zeta_{N+n}(s)'' = \left( C_0'' + \frac{c''}{N+n-1} \right) ((N+n-1)^{(-s-3)} - (N+n)^{(-s-3)})$$

$$\zeta_{N+n}(s)'' - \zeta_{N+n+1}(s)'' = \left( C_0'' + \frac{c''}{N+n} \right) ((N+n)^{(-s-3)} - (N+n+1)^{(-s-3)})$$

$$\zeta_{N+n+1}(s)'' - \zeta_{N+n+2}(s)'' = \left( C_0'' + \frac{c''}{N+n+1} \right) ((N+n+1)^{(-s-3)} - (N+n+2)^{(-s-3)})$$

So instead of vanishing pairs of terms according to the schemes in FAQ #21 of type

$$T'_{N,n} \;=\; C'\,((N+n)^{(s-4)} - (\,N+n)^{(s-4)}) \;=\; 0 \tag{9}$$

$$T''_{N,n} \;=\; C''\,((N+n)^{(-s-3)} - (\,N+n)^{(-s-3)}) \;=\; 0 \tag{10}$$

we thus actually have nonvanishing contributions from the remainders equal to, respectively,

$$T'_{N,n} \;=\; c' \left( \frac{(N+n)^{(s-4)}}{N+n} - \frac{(N+n)^{(s-4)}}{N+n-1} \right) \;=\; -\frac{c'\,(N+n)^{(s-5)}}{N+n-1} \tag{11}$$

$$T''_{N,n} \;=\; c'' \left( \frac{(N+n)^{(-s-3)}}{N+n} - \frac{(N+n)^{(-s-3)}}{N+n-1} \right) \;=\; -\frac{c''\,(N+n)^{(-s-4)}}{N+n-1} \tag{12}$$

Since we are considering very large $N$ and/or $n$, we can disregard the term –1 in the denominators, and thus get

$$T'_{N,n} \;=\; -c'\,(N+n)^{(s-6)} \tag{13}$$

$$T''_{N,n} \;=\; -c''\,(N+n)^{(-5-s)} \tag{14}$$

When summed over $n$ with $n$ from 0 to $k$, the extra contributions in (13) and (14) give corrections to the sums in (5a) and (5b) in FAQ #21 as follows

$$\sum_{n=0}^{k} T'_{N,n} = \sum_{n=0}^{k} (-c'\,(N+n)^{(s-6)}) \tag{15}$$

$$\sum_{n=0}^{k} T''_{N,n} = \sum_{n=0}^{k} (-c''\,(N+n)^{(-5-s)}) \tag{16}$$

The sums in (15) and (16) are smaller in absolute value than some finite constant times the following sum, which we can thus use to estimate the order of magnitude of the sums in (15) and (16)

$$T_{N,\,k} = \sum_{n=0}^{k} \frac{1}{(N+n)^5} \qquad (17)$$

When $k$ tends to infinity, the sum in (17) evaluates to

$$T_{N,\,\infty} = -\frac{1}{24}\,\Psi(4, N) \qquad (18)$$

where the polygamma function $\Psi(4, N)$ above is the fourth derivative of the digamma function (see Abramowitz and Stegun, eqs 6.31 and 6.4.1). It gives finite values for finite $N$, and tends to zero when $N$ tends to infinity.

The contributions to the resulting sums in (5a) and (5b) in FAQ #21 from the remainders in (2a) and (2b) in FAQ #21 are thus negligible even if there are infinitely many of them.

---

## FAQ #21b (page 18)

**"In Remark B.2 on page 18 in the preprint, you say that when $N$–> oo, the limit of (20) must be either 0, 1 or infinity. It would seem to me that there is nothing that prevents the limit from being different from 1; there is nothing that forces the limit to be 0, 1, or infinity; there is nothing that compels the limit even to exist."**

ANSWER: I have understood that your argument is based on the following limited analysis (cf also FAQ #7).

From Cauchy's theorem one can derive the following two alternative expressions for the zeta-function in Appendix A [by inserting respectively (12a), (A9) into (9), and (12b), (A9) into (11)]

$$\zeta(s) = \frac{\pi^s \left( -2^{(s-1)} N^s \left( 1 - \frac{1}{24}\frac{s\,(s-1)}{N^2} + \mathrm{O}\left(\frac{1}{N^{3.}}\right)\right) + \left( \sum_{n=1}^{N} (2\,n-1)^{(s-1)} \right) s \right)}{\cos\left(\frac{1}{2}\,s\,\pi\right)\,(-1+2^{(1-s)})\,\Gamma(s+1)} \qquad (a1)$$

$$\zeta(s) = \frac{-N^{(1-s)}\left( 1 - \frac{1}{24}\frac{s\,(s-1)}{N^2} + \mathrm{O}\left(\frac{1}{N^{3.}}\right)\right) + 2^s \left( \sum_{n=1}^{N} (2\,n-1)^{(-s)} \right)(1-s)}{(-1+2^s)\,(1-s)} \qquad (a2)$$

where the decimal point after the powers in the remainders designate a remainder such as $\mathrm{O}(\ln(N)^2/N^3)$, which is greater than $\mathrm{O}(1/N^3)$ but (for sufficiently large $N$) is smaller than any remainder $O(1/N^{3-|\varepsilon|})$, no matter how small ε may be.

The following approximations of (a1) and (a2), corresponding to (13) and (14) in the preprint, are obtained by setting the remainders in (a1) and (a2) equal to zero.

$$\zeta_N(s)' = \frac{\pi^s \left( -2^{(s-1)} N^s \left( 1 - \frac{1}{24} \frac{s\,(s-1)}{N^2} \right) + \left( \sum_{n=1}^{N} (2\,n-1)^{(s-1)} \right) s \right)}{\cos\left(\frac{1}{2} s\,\pi\right) (-1 + 2^{(1-s)})\, \Gamma(s+1)} \tag{b1}$$

$$\zeta_N(s)'' = \frac{-N^{(1-s)} \left( 1 - \frac{1}{24} \frac{s\,(s-1)}{N^2} \right) + 2^s \left( \sum_{n=1}^{N} (2\,n-1)^{(-s)} \right) (1-s)}{(-1+2^s)\,(1-s)} \tag{b2}$$

Now use (a1), (b1) and (a2), (b2), respectively, to calculate the following differences,

$$\zeta_N(s)' - \zeta(s) = \frac{\pi^s\, 2^{(s-1)} N^s\, O'(N^{(-3.)})}{\cos\left(\frac{1}{2} s\,\pi\right) (-1 + 2^{(1-s)})\, \Gamma(s+1)} \tag{c1}$$

$$\zeta_N(s)'' - \zeta(s) = -\frac{N^{(1-s)}\, O''(N^{(-3.)})}{(s-1)\,(-1+2^s)} \tag{c2}$$

Calculate their quotient

$$\frac{\zeta_N(s)' - \zeta(s)}{\zeta_N(s)'' - \zeta(s)} = -\frac{\pi^s\, 2^{(s-1)} N^s\, O'(N^{(-3.)})\,(s-1)\,(-1+2^s)}{\cos\left(\frac{1}{2} s\,\pi\right) (-1 + 2^{(1-s)})\, \Gamma(s+1)\, N^{(1-s)}\, O''(N^{(-3.)})} \tag{d}$$

Simplifying (d) and setting

$$\frac{O'(N^{(-3.)})}{O''(N^{(-3.)})} = K(s) + O(N^{(-1.)}) \tag{e}$$

gives

$$\frac{\zeta_N(s)' - \zeta(s)}{\zeta_N(s)'' - \zeta(s)} = \frac{1}{2} \frac{(s-1)\, N^{(2\,s-1)}\, K(s)\, \pi^s\, (-4^s + 8^s)}{\cos\left(\frac{1}{2} s\,\pi\right) \Gamma(s+1)\,(-2+2^s)} + O(N^{(2\,\sigma-2.)}) \tag{f}$$

This expression has the following properties for finite, nonvanishing K(*s*) in the limit when *N* tends to infinity – if the limit exits at all: For σ < ½ the right-hand side tends to zero, whereas for σ > ½ the right-hand side tends to infinity. For σ = ½ the limit of the right-hand side depends on K(*s*), and could thus in principle be any number (and also depend on *t*). And nothing of course even forces K(*s*) to be finite and nonvanishing as stated above.

This can thus be taken as the basis for a statement that my assertion that the limit in (f) is either exactly 1, or 0, or infinity is not correct, as you do. However, my assertion is based on the much more detailed calculation of the quotient in Appendix B in the preprint, giving as final result the expression (18),

$$\frac{\zeta_N(s)' - \zeta(s)}{\zeta_N(s)'' - \zeta(s)} = \frac{1}{2}\frac{N^{(2s-1)}\pi^s(-4^s+8^s)}{\cos\left(\frac{1}{2}s\pi\right)(-2+2^s)\Gamma(s-3)(s+2)(s+1)s} + O(N^{(2\sigma-2)}) \quad (18)$$

which is similar to (f) above, but in the derivation of which K(*s*) is explicitly calculated, so that the quotient indeed has the properties I assert, and which are given in Appendix B.

If we wish, we can calculate K(*s*) in (f) from (18) as follows. Setting the right-hand sides of (f) and (18) equal and solving for K(*s*), we get

$$K(s) = \frac{(s-2)(s-3)}{(s+2)(s+1)}$$

which for $s = ½ + i\,t$ gives the following absolute value

$$\left|K\left(\frac{1}{2}+i\,t\right)\right| = 1$$

The result of this more detailed analysis thus proves that the quotient in (18) indeed must be either 0, 1 or infinity in the limit $N \rightarrow \infty$, as I show in the preprint.

---

## FAQ #22 (page 18)

**"It is always good to check theoretical derivations by numerical examples whenever possible. If the results come out right, then this doesn't necessarily prove anything of course. But if they do come out wrong, then that's a clear sign that the derivation is most probably flawed somewhere. Have you checked the final expression (18) on page 18 in your preprint in this way?"**

ANSWER: Yes, indeed I have. Inserting (12a) and (A9) into (16), we get

$$\zeta_N(s)' - \zeta(s) = \frac{\pi^s \left( -2^{(s-1)} N^s \left( 1 - \frac{1}{24} \frac{s\,(s-1)}{N^2} \right) + \left( \sum_{n=1}^{N} (2\,n-1)^{(s-1)} \right) s \right)}{\cos\left(\frac{1}{2} s\,\pi\right) (-1 + 2^{(1-s)})\,\Gamma(s+1)} - \zeta(s)$$

Similarly inserting (12b) and (A9) into (17), we get

$$\zeta_N(s)'' - \zeta(s) = \frac{-N^{(1-s)} \left( 1 - \frac{1}{24} \frac{s\,(s-1)}{N^2} \right) + 2^s \left( \sum_{n=1}^{N} (2\,n-1)^{(-s)} \right) (1-s)}{(-1+2^s)\,(1-s)} - \zeta(s)$$

The quotient of the above two expressions should thus be compared to (18), i e

$$\frac{\zeta_N(s)' - \zeta(s)}{\zeta_N(s)'' - \zeta(s)} = \frac{1}{2} \frac{N^{(2s-1)}\,\pi^s\,(-4^s + 8^s)}{\cos\left(\frac{1}{2} s\,\pi\right) (-2+2^s)\,\Gamma(s-3)\,(s+2)\,(s+1)\,s} + O(N^{(2\sigma-2)}) \qquad (18)$$

For, say, $N = 10^5$ and $s = 0.7 + 30\,i$, we have

$\zeta_N(s)' - \zeta(s)$ = -0.373511085416521 $10^{-15}$ - 0.174631217030210 $10^{-15}$ $i$
$\zeta_N(s)'' - \zeta(s)$ = 0.333100167313817 $10^{-17}$ - 0.338795678198075 $10^{-17}$ $i$

So now we can calculate (18):

Left - hand side: -28.906181537 - 81.826470700 $i$
Right-hand side: -28.906181640 - 81.826470663 $i$ + O(0.001)

---

Difference: 0.103 $10^{-6}$ - 0.037 $10^{-6}$ $i$

The difference is thus well within the accuracy defined by the remainder. This numerical example thus offers some support that (18) may be correct.

---

## FAQ #23 (page 20)

**"Your discussion in Remark C.II on page 20 in the preprint for the case ½ < σ < 1 is very condensed. Can you explain what happens there in somewhat more detail ?"**

ANSWER: The expression for $\Lambda_N$ in Theorem C.II is given in (ι) in the range $0 < \sigma \leq ½$ by using (18) from Sect. 8 in the preprint,

$$\Lambda_N = \left| \frac{\zeta_N(s_0)'}{\zeta_N(s_0)''} \right| = \lim_{s \to s_0} \left| \frac{1}{2} \frac{N^{(2s-1)}\,\pi^s\,(-4^s+8^s)}{\cos\left(\frac{1}{2} s\,\pi\right) (-2+2^s)\,\Gamma(s-3)\,(s+2)\,(s+1)\,s} + O(N^{(2\sigma-2)}) \right| \qquad (\iota)$$

For ½ < σ < 1, the expression on the right-hand side in (ι) above obviously does not exist uniformly with respect to *N*, *i e* we cannot use it to support (λ) in Theorem C.II in this range. However, to cover this case we can proceed as follows.

Instead of (9) and (11) in Sect. 5 in the preprint, we use (10) and a new equation (11a), in which we have made the substitution $s \to 1-s$ in (11),

$$\zeta(1-s) = \frac{\pi^{(1-s)}\left(-2^{(-s)}N^{(1-s)}\,\mathrm{E}_N(1-s) + \left(\sum_{n=1}^{N}(2n-1)^{(-s)}\right)(1-s)\right)}{\sin\left(\frac{1}{2}s\pi\right)(-1+2^{s})\,\Gamma(2-s)} \qquad (10)$$

$$\zeta(1-s) = \frac{-N^{s}\,\mathrm{E}_N(s) + 2^{(1-s)}\left(\sum_{n=1}^{N}(2n-1)^{(s-1)}\right)s}{(-1+2^{(1-s)})\,s} \qquad (11a)$$

Note that if we make the substitution $s \to 1-s$, then these equations become identical to (9) and (11) (as they should). This thus means that we can perform all the calculations from Sect. 5 through (18) in Sect. 8 in exactly the same way as in the preprint, but now with the substitution $s \to 1-s$. Instead of (18), we then get (18′) below (note that here the remainder is then $O(N^{-2s})$), which is thus still finite in the critical strip $0 < Re(s) < 1$),

$$\frac{\zeta_N(1-s)' - \zeta(1-s)}{\zeta_N(1-s)'' - \zeta(1-s)} =$$

$$\frac{1}{2}\frac{N^{(1-2s)}\pi^{(1-s)}(-4^{(1-s)}+8^{(1-s)})}{\cos\left(\frac{1}{2}(1-s)\pi\right)(-2+2^{(1-s)})\,\Gamma(-2-s)(3-s)(2-s)(1-s)} + O(N^{(-2\sigma)}) \qquad (18')$$

Analogously to (ι) in Theorem C.II in Appendix C in the preprint, we thus get from (18′) at a zero $s_0$ of the zeta-function, *i e* when $\zeta(1-s_0) = 0$,

$$\Lambda_N{}' = \left|\frac{\zeta_N(1-s_0)'}{\zeta_N(1-s_0)''}\right| =$$

$$\lim_{s \to s_0}\left|\frac{1}{2}\frac{N^{(1-2s)}\pi^{(1-s)}(-4^{(1-s)}+8^{(1-s)})}{\cos\left(\frac{1}{2}(1-s)\pi\right)(-2+2^{(1-s)})\,\Gamma(-2-s)(3-s)(2-s)(1-s)} + O(N^{(-2\sigma)})\right| \qquad (\iota')$$

Here $\Lambda_N'$ in (ι′) exists uniformly with respect to *N* in the range ½ ≤ σ < 1. Continuing the proof as in Theorem C.II, we thus next consider the left-most part of (ι′),

$$\Lambda_N{}' = \left| \frac{\zeta_N(1 - s_0)'}{\zeta_N(1 - s_0)''} \right|$$

A variable transformation like $s$ –> 1 – $s$ above can of course always be made, but the crucial point here is whether the left-hand part of (ι') above can then be related to the corresponding left-hand part of (ι). As remarked in Sect. 6 in the preprint, the functions $\zeta_N(s)'$ and $\zeta_N(s)''$ do not obey the functional equation, so despite the fact that the numerator and denominator in $\Lambda_N'$ above differ from the numerator and denominator in $\Lambda_N$ only by a constant factor in $s$, the quotient $\Lambda_N'$ above is not necessarily equal to $\Lambda_N$. However, this problem is solved by the fact that in the limit $N \rightarrow \infty$, the functions $\zeta_N(1 - s)'$ and $\zeta_N(1 - s)''$ are both equal to $\zeta(1 - s)$ since, as discussed in Sects. 6 and 7 in the preprint, they differ from the zeta-function only by their remainders, which are then zero. In this limit, the quotients $\Lambda_N'$ and $\Lambda_N$ are thus indeed equal, so the limit

$$\lim_{N \to \infty} \Lambda_N{}' = \lim_{N \to \infty} \Lambda_N = \lim_{N \to \infty} \left| \frac{\zeta_N(s_0)'}{\zeta_N(s_0)''} \right| \qquad (\kappa)$$

thus exists, and according to (λ) in Theorem C.II has the value

$$\lim_{N \to \infty} \left| \frac{\zeta_N(s_0)'}{\zeta_N(s_0)''} \right| = 1 \qquad (\lambda)$$

Hence the double limit in (μ) in Theorem C.II exists and has the value in (ν). Theorem C.II in Sect. C in the preprint is thus proven for all σ within the critical strip 0 < σ < 1.

---

14 October 2013

Arne Bergstrom